\documentclass[10pt,a4paper]{article}
\usepackage{fullpage}   
\usepackage[pdftex]{graphicx} %\usepackage[dvips]{graphicx} \usepackage{psfrag} % incompatibile con dvipdfm
\usepackage{amsmath, amsfonts, amssymb, amsbsy}
\usepackage{latexsym} \usepackage{stmaryrd} 
\usepackage[textsize=footnotesize]{todonotes}
\usepackage[labelfont=bf]{caption}
\usepackage{hyperref}
\hypersetup{
    colorlinks,%
    citecolor=black,%
    filecolor=black,%
    linkcolor=black,%
    urlcolor=black
}
\usepackage{subcaption}
\usepackage{multirow}
\usepackage{verbatim}
\usepackage{graphicx}%
\usepackage{multirow}%
\usepackage{amsthm}%
\usepackage{mathrsfs}%
\usepackage[title]{appendix}%
\usepackage{xcolor}%
\usepackage{textcomp}%
\usepackage{manyfoot}%
\usepackage{booktabs}%
\usepackage{algorithm}%
\usepackage{algorithmicx}%
\usepackage{algpseudocode}%
\usepackage{listings}%
\usepackage{arydshln}
\usepackage{empheq}
\usepackage[scr=boondox]{mathalfa}
\usepackage{subcaption}

\begin{document}

%%%%%%%%%%%%%%%%%%%%%%%%%%%%%%%%%%%%%%%%%%%%%% NUMERAZIONE
\newtheorem{theorem}{Theorem}[section]
\newtheorem{proposition}[theorem]{Proposition}
\newtheorem{corollary}[theorem]{Corollary}
\newtheorem{lemma}[theorem]{Lemma}
\newtheorem{definition}[theorem]{Definition}
\newtheorem{remark}[theorem]{Remark}
\newtheorem{example}[theorem]{Example}

%%%%%%%%%%%%%%%%%%%%%%%%%%%%%%%%%%%%%%%%%%%%%%% MACRO TESTO
\renewcommand{\proof}{\noindent \textbf{Proof. }}
\renewcommand{\qed}{ \hfill {\vrule width 6pt height 6pt depth 0pt} \medskip }
\newcommand{\separe}{\medskip \centerline{\tt -------------------------------------------- } \medskip}
%\newcommand{\doit}[1]{\hfill {\small \tt [#1]}}
%\newcommand{\note}[1]{{\footnotesize #1}}

%%%%%%%%%%%%%%%%%%%%%%%%%%%%%%%%%%%%%%%%%%%%%%%% MACRO MATEMATICA
\newcommand{\eps}{\varepsilon} \renewcommand{\det}{\mathrm{det}} \newcommand{\argmin}{ \mathrm{argmin} \,}
\def\interior{\mathaccent'27} 
\newcommand{\weakto}{ \rightharpoonup}  \newcommand{\weakstarto}{\stackrel{*}{\rightharpoonup}}
\newcommand{\R}{\mathbb{R}}
\newcommand{\F}{\mathcal{F}}\newcommand{\wF}{\tilde{\mathcal{F}}}
\newcommand{\I}{\mathcal{I}}\newcommand{\J}{\mathcal{J}}
\renewcommand{\H}{\mathcal{H}}
\newcommand{\U}{\mathcal{U}}
\newcommand{\V}{\mathcal{V}}
\newcommand{\K}{\mathcal{K}}
\newcommand{\E}{\mathcal{E}}
\renewcommand{\R}{\mathbb{R}} 
\renewcommand{\c}{\gamma} 
\renewcommand{\v}{\rho} 
\newcommand{\bu}{\boldsymbol{u}}
\newcommand{\bg}{\boldsymbol{g}}
\newcommand{\bn}{\boldsymbol{n}}
\newcommand{\stress}{\boldsymbol{\sigma}} \newcommand{\strain}{\boldsymbol{\epsilon}} 
\newcommand{\jump}[1]{\llbracket #1\hspace{1pt}\rrbracket}

\def\Xint#1{\mathchoice
 {\XXint\displaystyle\textstyle{#1}}%
 {\XXint\textstyle\scriptstyle{#1}}%
 {\XXint\scriptstyle\scriptscriptstyle{#1}}%
 {\XXint\scriptscriptstyle\scriptscriptstyle{#1}}%
 \!\int}
\def\XXint#1#2#3{{\setbox0=\hbox{$#1{#2#3}{\int}$}
 \vcenter{\hbox{$#2#3$}}\kern-.5\wd0}}

\def\ddashint{\Xint=} % ----------------------------- questi sono i comandi 
\def\dashint{\Xint-}

\newcommand{\bb}{\color{black}} 
\newcommand{\bl}{\color{black}} 
\newcommand{\rr}{\color{black}} 
\newcommand{\cc}{\color{black}} 
\newcommand{\vv}{\color{black}} 
\newcommand{\gv}{\color{black}} 
\newcommand{\mg}{\colorblack} 
\definecolor{burgundy}{rgb}{0.5, 0.0, 0.13}
\newcommand{\ap}{\color{black}}
\newcommand{\aprev}[1]{{\color{black} #1}}
\newcommand{\NOTE}[1]{\todo[inline, color=yellow!20]{\tt #1}}

%%%%%%%%%%%%%%%%%%%%%%%%%%%%%%%%%%%%%%%%% MACRO LUIGI

\newcommand\uu{\mathbf{u}}
\newcommand\dd{\textrm{d}}
\newcommand\ww{\textrm{w}}
%\newcommand\strain{\bm{\varepsilon}}

%=======================================================================================
%% BEGIN ADDED definitions, commands 
% comments
\definecolor{frenchblue}{rgb}{0.0, 0.45, 0.73}
% format
\newcommand\mtx[1]{\textrm{#1}}
\newcommand\hh{\mtx{h}}
\newcommand\un{_\mtx{n}}
% displ
%\newcommand\uu{\mathbf{u}}
\newcommand\huu{\hat{\mathbf{u}}}
\newcommand\uuh{\mathbf{u}_\hh}
\newcommand\duu{\delta\mathbf{u}}
\newcommand\duuh{\delta\mathbf{u}_\hh}
\newcommand\dtuu{\dot{\mathbf{u}}}
% phase-field
%\newcommand\dd{\textrm{d}}
\newcommand\hdd{\hat{\mathbf{d}}}
\newcommand\ddh{\textrm{d}_\hh}
\newcommand\ddd{\delta\textrm{d}}
\newcommand\dddh{\delta\textrm{d}_\hh}
\newcommand\dtdd{\dot{\textrm{d}}}
% stress, strain
%\newcommand\strain{\bm{\varepsilon}}
%\newcommand\stress{\bm{\sigma}}
% parametric coord
\newcommand\bxi{\bm{\xi}}
% shape functions, discretization
\newcommand\Rsf{R_{\textbf{i},\textbf{p}}}
\newcommand\Rsfj[1]{R_{\textbf{j},\textbf{p}}}
\newcommand\Rd{R^{\mtx{d}}_{\textbf{i},\textbf{p}}}
\newcommand\RRu{\mathbf{R}^{\mathbf{u}}_{\textbf{i},\textbf{p}}}
\newcommand\dRsf[1]{\frac{\partial R_{\textbf{i},\textbf{p}}}{\partial x_{#1}}}
\newcommand\ddRsf[1]{\frac{\partial^2 R_{\textbf{i},\textbf{p}}}{\partial x^2_{#1}}}
\newcommand\BBu{\mathbf{B}^{\mathbf{u}}_{\textbf{i},\textbf{p}}}
\newcommand\BBd{\mathbf{B}^{\mtx{d}}_{\textbf{i},\textbf{p}}}
\newcommand\BBdj[1]{\mathbf{B}^{\mtx{d}}_{\textbf{#1},\textbf{p}}}
\newcommand\CCd{\mathbf{C}^{\mtx{d}}_{\textbf{i},\textbf{p}}}
\newcommand\CCdj[1]{\mathbf{C}^{\mtx{d}}_{\textbf{#1},\textbf{p}}}
% math other
\newcommand\eye{\mathbf{I}}
\newcommand\one{\mathbf{1}}
\newcommand{\oo}{\mathscr{o}}
% text
\newcommand\ie{i.e., }
\newcommand\eg{e.g., }
\newcommand\ii{\textbf{i}}
\newcommand\jj{\textbf{j}}
% algorithms
\renewcommand{\algorithmicrequire}{\textbf{Input:}}
\renewcommand{\algorithmicensure}{\textbf{Output:}}
% french numbering for Matteo's addition
\newcommand*\latinnumeral[1]
{\ifcase#1\unskip\or\unskip\or bis \or ter \fi}
\newcommand*\latintags
{\xdef\startingequationnumber{\the\numexpr\value{equation}+1\relax}%
	\setcounter{equation}{0}%
	\aftergroup\resetarabicequations
	\def\theequation
	{\startingequationnumber~\latinnumeral{\value{equation}}}}
\newcommand*\resetarabicequations
{\setcounter{equation}{\numexpr\startingequationnumber\relax}}

%%%%%%%%%%%%%%%%%%%%%%%%%%%%%%%%%%%%%%%%%%%%% ARTICLE 

% !TeX root = AT1IV.tex

%%%%%%%%%%%%%%%%%%%%%%%%%%%%%%%%%%%%%%%%%%%%%%%%%%%% COVER 
%\phantom{.}

\vspace{0.5cm}
\noindent {\LARGE {\bf \boldmath{$AT_1$} fourth-order isogeometric phase-field modeling \\[9pt] of brittle fracture}}

\vspace{24pt}

\begin{small}
\noindent {\bf L.~Greco, E.~Maggiorelli, M.~Negri, A.~Patton,   A.~Reali}

\bigskip
\noindent {\bf Abstract.} 
%% Text of abstract
A crucial aspect in phase-field modeling, based on the variational formulation of brittle fracture, is the accurate representation of how the fracture surface energy is dissipated during the fracture process in the energy competition within a minimization problem. In general, the family of $AT_1$ functionals showcases a well-defined elastic limit and narrow transition regions before crack onset, as opposed to $AT_2$ models. On the other hand, high-order functionals provide similar accuracy as low-order ones but allow for larger mesh sizes in their discretization, remarkably reducing the computational cost. In this work, we aim to combine both these advantages and propose a novel $AT_1$ fourth-order phase-field model for brittle fracture within an isogeometric framework, which provides a straightforward discretization of the high-order term in the crack surface density functional. For the {introduced} $AT_1$ functional, we first prove a $\Gamma$-convergence result (in both the continuum and {discretized} isogeometric setting) based on a careful study of the optimal transition profile, which ultimately provides the explicit correction factor for the toughness and the exact size of the transition region. Fracture irreversibility is modeled by monotonicity of the damage variable and is conveniently enforced using the {Projected Successive Over-Relaxation} algorithm. %Mesh size is choosen as a function of the size of the transition region. 
Our numerical results indicate that the proposed fourth-order $AT_1$ model {is more accurate than the considered} lower-order $AT_1$ and $AT_2$ models; this allows to employ larger mesh sizes, entailing a lower computational cost. 
%\NOTE{AP: Non so se come incipit dell'abstract e' ok parlare solo dell'energia dissipata (anche se ho cercato di dare l'idea di competizione delle energie) in quanto per modelli AT si intende quadratic degradation function e linear o quadratic local part nella dissipazione a seconda di AT1/AT2; e' anche vero che fissata la degradation function quadratica poi quello che cambia effettivamente e' la dissipazione. Che ne pensate?}
\bigskip
\\\\\noindent {\bf AMS Subject Classification.} 
\end{small}

\section*{Introduction}

Numerical simulation has the potential to serve as a decision-making tool in engineering, for example, at the design stage of structural elements of mechanical, civil, or aviation systems, ultimately reducing the reliance on expensive and time-consuming experimental tests to detect potential failure due to fracture. This requirement has led to developing many theoretical models and numerical methods.

The reference work on which fracture mechanics was founded is that of Griffith \cite{Griffith1921}, based on the ``competition'' between stored elastic and dissipated fracture  energy. More precisely, Griffith's criterion asserts that fracture propagation occurs whenever the {\it energy release} $G$ (the configurational variation of elastic energy with respect to the crack surface) reaches a critical value $G_c$ (a material parameter, called {\it toughness}). In this approach, the physical crack is a (possibly branched) path in 2{\sc d} or a surface in 3{\sc d}. Alternatively, on the basis of \cite{Ambrosio1990999}, Bourdin et al.~\cite{BOURDIN2000797} introduced the phase-field approach for fracture: %Several comments are in order from the mathematical, mechanical and numerical point of view. 
in this setting, the crack is represented by means of a scalar variable $v$ taking values in $[0,1]$, with $v=1$ and $v=0$ corresponding respectively to fracture and sound material.  {\it Phase-field energies} usually take the form 
$$
	\F_\eps ({\boldsymbol{u}}, v) = \E_\eps ({\boldsymbol{u}},v) + G_c \, \K_\eps (v) , 
$$
where $\boldsymbol{u}$ is the displacement, $\eps$ is the internal length, $\E_\eps$ is the elastic energy, while $\K_\eps$ is the crack surface energy. There are nowadays many choices for both $\K_\eps$ and $\E_\eps$. 
For an extensive overview of phase-field work, we refer the interested reader to the following review works and books: \cite{ambati2015review,LI2023109419,ZHUANG2022108234,BourdFrancMar08,Wick20}.

In general, low-order energies take the form 
$$
	\K_\eps (v) = \int_\Omega \left(\eps^{-1} \phi (v) + \eps | \nabla v |^2 \right) dx  ,
$$
while their high-order counterparts \cite{BORDEN2014100,HeschSchussDittmannFrankeWeinberg_CMAME16} read
$$
	\K_\eps (v) = \int_\Omega \left(\eps^{-1} \phi (v) + \eps | \nabla v |^2 + \eps^3 | \Delta v |^2 \right) dx . 
$$
In the literature, the former are called second-order and the latter fourth-order; making reference to the associated Euler-Lagrange equations. 

There are several choices for the function $\phi$, see, e.g.,~\cite{Wu_JMPS17} and the above references. Here, we mention only the most common: $\phi(v) = v$ and $\phi(v) = v^2$, corresponding to the so-called $AT_1$ and $AT_2$ functionals. From the mechanical point of view, $AT_1$ functionals have in general, better properties than $AT_2$ functionals since they result in a clear linear elastic regime before the onset of fracture. Indeed, as observed by Pham et al.~\cite{pham2011gradient}, the linearity of  $\phi$ introduces an analytical elastic limit of order $\sqrt{ G_c \mu / \eps }$ under which fracture does not occur. This consideration also makes this model interesting in the field of plasticity (see, e.g.,~Marengo et al.~\cite{MARENGO2023115992}).
On the contrary, $AT_2$ promotes damage at arbitrarily small values of stress. Moreover, for given values of internal length and mesh size, the damage profile around the crack is usually narrower with $AT_1$. {It follows that numerical simulations carried out using the $AT_1$ model are closer to representing the physics behavior of the specimen and, therefore, a preferable modeling choice.} Most often, simulations based on high-order functionals are more accurate approximating the dissipated energy, see, e.g.,~\cite{greco2024higher}, and allow for larger mesh sizes, leading to considerably lower computational costs despite the extra effort related to second-order derivatives. As far as the elastic energy $\E_\eps$, it is nowadays standard to employ energy splits to provide the correct mechanical behavior under tension, compression, and shear. Generally, elastic energies take the form
%{\tt Since we focus on the approximation of surface energies, for the elastic energy we take into account the classic volumetric-deviatoric split \cite{+} 
$$
	\E_\eps ({\boldsymbol{u}} ,v) = \int_\Omega  \left(\psi_\eps (v, \nabla v ) W_+ ( \strain ) + W_- ( \strain ) \right) dx ,
%	\qquad 
%	W_+ (\strain) : = \bb \mu \bl | \strain_d |^2 +  \bb \kappa \bl | \strain_v^+ |^2,
%	\qquad
%	W_- (\strain)  := \bb \kappa \bl | \strain_v^- |^2,
$$
where $\psi_\eps$ is the degradation function and $\strain$ is the strain.  In the literature, there are actually several choices for both the terms $W_\pm$ (e.g., volumetric-deviatoric \cite{Amor2009JMPS}, spectral \cite{MieheWelschingerHofacker10}{, and the recent star-convex \cite{Vicentini2024} decomposition}) and the degradation function $\psi_\eps$ \cite{Wu_JMPS17}, with different mechanical outcomes, in particular under cyclic loadings. 

The phase-field approach relies on the solid basis of {\it $\Gamma$-convergence} \cite{DalMaso93,Braides98} of the phase-field energy to the {\it sharp crack energy}, as the internal length $\eps$ vanishes. In this direction, besides \cite{Ambrosio1990999}, are nowadays available several interesting results, see e.g.~\cite{Braides98,CCF,N_CMAME20,BurgerEspositoZeppieri_MMS15}. 

{\it Evolutions} are usually defined by means of time discrete incremental problems: at each time instant an {\it equilibrium configuration} of the system is computed, not necessarily being a global minimizer, subject to an irreversibility constraint. Numerically, different schemes are employed to find an equilibrium configuration, e.g., staggered \cite{BOURDIN2000797}, monolithic \cite{FarrellMaurini_IJNME17}, and active set \cite{HWW15}. In this respect, note that the energy $\F_\eps$ is not (jointly) convex but only separately convex; therefore, in general, there are many equilibria and the solution to the incremental problem could be non-unique. Studying the time-continuous limit of the time-discrete evolutions \cite{KRZ,KneesNegri_15,Maggiorelli24,MaggiorelliNegri24} reveals a peculiar feature of this model, actually shared by any rate-independent evolution for non-convex energies. There are two regimes of propagation: stable (or steady-state) and unstable (or catastrophic) \cite{Maggiorelli24,MaggiorelliNegri24}. In the former, the evolution is continuous in time, and, noteworthy, it satisfies a phase-field version of Griffith's criterion. In the latter, a discontinuous instantaneous propagation occurs and Griffith's criterion is not always satisfied; in this case, dynamical or other rate-dependent models should be preferred. 
 
\bigskip
 
In this work, we study a high-order $AT_1$ functional of the form $\F_\eps ({\boldsymbol{u}} , v) = \E_\eps ({\boldsymbol{u}}, v) + G_c \K_\eps(v)$ {(see \S \ref{s.evolution})}. The elastic energy features the volumetric-deviatoric split \cite{Amor2009JMPS}
$$
	\E_\eps ({\boldsymbol{u}} ,v) = \int_\Omega  \left(\psi_\eps (v) W_+ ( \strain ) + W_- ( \strain ) \right) dx, 
	\qquad 
	W_+ (\strain) = \mu | \strain_d |^2 +  \kappa | \strain_v^+ |^2,
	\qquad
	W_- (\strain)  = \kappa | \strain_v^- |^2,
$$%
with a quadratic degradation function, i.e.,~$\psi_\eps (v) =  ( v-1)^2 + \eta_\eps$ for $\eta_\eps > 0$. The crack surface energy reads 
$$
	\K_\eps (v) = \frac{1}{c_\rho} \int_\Omega \left(\eps^{-1} v+ \eps | \nabla v |^2  + \rho \, \eps^3 | \Delta v |^2 \right) dx  ,
$$
where $\rho >0$ is a parameter weighting the effects of the high-order term and $c_\rho$ is the normalizing constant that depends on $\rho$. We underline that the present model is different from the one in \cite{nagaraja2023deterministic}, where a high-order AT1 model is used to analyze anisotropy problems without considering the phase-field gradient term in the crack surface energy, as well as the technical tuning of $c_\rho$.

First of all, we provide a $\Gamma$-convergence result {in \S \ref{s.gammaConv}}. To this end, a main crucial point is the characterization of $c_\rho$ (as a function of $\rho$)  in such a way that the toughness is recovered exactly; in other terms, the factor $c_\rho$ should be chosen in such a way that the sharp crack $\Gamma$-limit reads, in 2{\sc d}, 
$$
	\int_{\Omega \setminus J_u} W ( \strain) \, dx + G_c \mathcal{H}^1 ( J_u) ,
$$
where $\mathcal{H}^1$ is the 1-dimensional Hausdorff measure of the displacement jump $J_u$ (loosely speaking, the length of the crack).

Technically, previous results in the literature, see, e.g., \cite{Braides98,N_CMAME20,BurgerEspositoZeppieri_MMS15}, do not apply here due to the combination of the high-order term and the constraint $v \in [0,1]$. The value $c_\rho$ is \bl usually given by  %generically characterized as 
$c_\rho = 2 \K (w_*)$, where $w_*$ is the solution of the constrained optimal profile problem in $\mathbb{R}_+$. \bl 
However, this characterization is not useful \bl for explicitly computing the optimal profile $w_*$.
Hence, we developed a novel line of proof \bl in \S \ref{convergence} based on approximations with unconstrained local optimal profiles, which ultimately yield an explicit form for $w_*$. As a by product, we show that $w_*$ is compactly supported in an interval $[0,R_*]$, where $R_*$ is explicitly computed in terms of the coefficients  of $\mathcal{K}_\epsilon$. Note that, on the contrary,  \bl the optimal profile of $AT_2$ is an exponential, supported in the whole $\mathbb{R}_+$. This technical point, has a couple of important benefits, for $AT_1$ versus $AT_2$: the transition region is narrower and the value of toughness is better approximated. As a result, numerical simulations are very accurate at relatively large values of the internal length and mesh size. 

\bigskip

After the mathematical part, the paper provides a detailed study of the numerical performance of the proposed $AT_1$ fourth-order model in \S \ref{sec:numerical test}. 
Dealing with high-order energies, we employ an isogeometric approach, as in, e.g., \cite{BORDEN2014100,HeschSchussDittmannFrankeWeinberg_CMAME16,greco2024higher}. 
Isogeometric Analysis (IgA) was initially developed by Hughes et al. \cite{HUGHES20054135} to extend and enhance the capabilities of the finite element (FE) method in the field of geometric modeling, resulting in a technique also appearing to be preferable to standard FEs in many applications on the basis of per-degree-of-freedom accuracy. In particular, IgA displays a unique blend of features that can be harnessed for tackling challenges associated with the modeling of higher-order differential operators, among which higher-order accuracy, robustness, geometric flexibility, and, in particular, $C^1$ and higher-order continuity. 
The phase-field method in combination with an isogeometric discretization proved to be successful in studying fracture in a wide variety of problems involving, e.g., the modeling of structural elements in statics such as solid shells \cite{AMBATI2016351} and Kirchhoff-Love shells \cite{PROSERPIO2021114019,PROSERPIO2020113363}, dynamics \cite{LI2023104783}, but also different types of materials, e.g., piezoceramics \cite{KIRAN2022108738}, piezoelectric composites \cite{KIRAN2023109181}, porous functionally graded structures \cite{Nguyen2023}, and rock-like materials \cite{LI2022108298}.
Furthermore, this technique has been applied within an IgA framework to the Cahn–Hillard equation \cite{GOMEZ20084333} and the isothermal Navier-Stokes–Korteweg equations \cite{GOMEZ20101828} and extended to a wide variety of areas of science and engineering, including the modeling of shape memory alloys \cite{DHOTE201548}, liquid–vapor flows with surfactants \cite{BUENO2016797}, but also biomedical applications such as tumor growth \cite{lorenzo2016,LORENZO2017515,lorenzo2019}, due to its capability to capture the interface implicitly without solving a moving boundary problem. 

In this study, we consider various benchmarks to highlight different properties of the newly proposed $AT_1$ isogeometric higher-order model compared to with classical $AT$-models from the literature.
For $AT_1$, we link the mesh size $h$ to the finite support of the optimal profile, i.e.,~to the value $R_*$, and just adopt the same choice of $h$ for $AT_2$ (since $R_*$ is not well defined). In this respect, it is important to remark that $R_*$ actually depends on the coefficient $\rho$, appearing in the surface energy. 
In all numerical simulations, irreversibility is consistently enforced by applying the Projected Successive Over-Relaxation algorithm \cite{MARENGO2021114137}, which seamlessly integrates with the isogeometric discretization of the high-order functional under examination \cite{greco2024higher}.
At first, we employ a pure tension setting to check the accuracy of the elastic limit. In this case, we compare only second- and fourth-order $AT_1$ functionals since $AT_2$ does not introduce a clear threshold for the stress before the onset of fracture. Our numerical results show that both functionals are extremely accurate, with a slightly better result for the fourth-order one.
Next, we consider the Double Cantilever Beam (DCB) and the Single Edge Notch (SEN) tests under tension. In these Mode I benchmarks, the crack geometry is a straight line, which allows to focus on the accuracy of the toughness. We compare fourth- and second-order $AT_1$ and $AT_2$ functionals. In essence, the fourth-order $AT_1$ functional performs better than the other $AT$-functionals, yielding very low errors on the toughness, at fixed mesh size, as well as higher convergence rates, with respect to the mesh size 
(see Table \ref{Tab: errors DCB AT} and \ref{Tab: errors tensile AT}).
Then, we consider a SEN under shear. In terms of toughness, the picture is very similar to that of the previous cases: the error is lower and the convergence rate higher. However, in terms of crack paths, fourth-order and second-order functionals produce quite different results (see Figure \ref{crack patter SEN shear novel AT}). We investigated this point in more detail in \S \ref{sec: sensitivity study} by changing the weight $\rho$ in front of the higher-order term in the surface integral. From this point of view, the picture is quite complex: for a fixed mesh size, smaller values of $\rho$ are better in terms of elastic limit, while larger values are better in terms of toughness. Finally, in \S \ref{trade} we study the trade-off between accuracy and computational cost, and, more precisely, we investigate first the actual savings in terms of control points for a fixed level of accuracy and, then, the role of the mesh as $R_*$ varies for different values of $\rho$. 
Finally, we draw our conclusions in \S \ref{sec:conclusions}.

\bl

\section{Phase-field energy and governing equations\label{s.evolution}} 

In this section, we introduce the phase-field energy and the evolution law in the time- and space-continuous setting. 

\subsection{Energy} 
%\NOTE{AP: Ho messo in bold il campo di spostamento $\bu$, dbcs $\bg$, normale $\bn$, variazione $\boldsymbol{\xi}$ (vettoriale) per la sezione 1.}
The reference configuration $\Omega$ is an open, bounded, and Lipschitz domain in $\mathbb{R}^2$. The set of admissible displacements is given by $\U= \{  \bu \in H^1(\Omega;\R^2) : \bu = \bg \bl  \text{ on }  \partial_D \Omega  \} $ where $\partial_D \Omega$ is a (relatively open) subset of the boundary $\partial \Omega$ and $\bg \in W^{1,\infty} (\Omega, \R^2)$ is the boundary displacement at a given time. The phase-field density $v$ belongs to the convex set 
$\V = H^1_{\mathrm{div}} (\Omega; [0,1])= \{ v \in H^1(\Omega) : v \in [0,1] , \, \Delta v  \in L^2(\Omega) \}$. We consider the following $AT_1$ energy functional: 
\begin{gather*}
	\F_\eps ( \bu, v) = \E_\eps ( \bu,v) + G_c \, \K_\eps (v) , \\ 
         \E_\eps (\bu,v)  =  \int_\Omega W_\eps (  v , \strain(\bu) )  \, dx ,  
\qquad
	\K_\eps (v) =  \frac{1}{c_\v} \int_\Omega \left(\eps^{-1} v + \eps | \nabla v |^2 + \v \, \eps^3 | \Delta v |^2 \right) dx ,
\end{gather*}
where $\strain (\bu) = \tfrac12 ( \nabla\bu + \nabla\bu^T)$ is the linear strain, $\eps > 0$ is the internal length, and $\v > 0$ is a parameter. 
%\NOTE{AP: Ok $\strain (\bu) = \tfrac12 ( D\bu + D\bu^T)$ oppure $\strain (\bu) = \tfrac12 ( \nabla\bu + \nabla\bu^T)$? Per essere uniformi con la notazione nell'energia del phase-field.}
Clearly, for $\v=0$, the energy boils down to the lower order $AT_1$ functional (see Appendix \ref{optimalAT1secondo}). The elastic energy density $W_\eps$ takes into account the volumetric-deviatoric split \cite{Amor2009JMPS,CCF} and is of the form 
$$
	W ( v , \strain ) :=  \psi_\eps (v) W_+ ( \strain ) + W_- ( \strain ),
	\qquad 
	W_+ (\strain) : =  \mu  | \strain_d |^2 +  \kappa  | \strain_v^+ |^2,
	\qquad
	W_- (\strain)  :=  \kappa  | \strain_v^- |^2 ,
$$where $\psi_\eps (v) =  ( v-1)^2 + \eta_\eps $ (for $\eta_\eps >0$) is the degradation function, $\strain_v^+ : = \tfrac12 \mathrm{tr}^+ ( \strain) \boldsymbol{I}$ and $\strain_v^- : = \tfrac12 \mathrm{tr}^- ( \strain) \boldsymbol{I}$  are respectively the \aprev{positive and negative volumetric parts} of the strain, \aprev{whereas} $\strain_{d} := \strain - \strain_v$ is \aprev{its} deviatoric part. %The function $\kappa$ appearing in the facture energy $\K$ may take the form $(1-v)^2$ or $[1-v]_+$. %respectively for $AT_2$ and $AT_1$. 
The phase-field energy $\F_\eps$ provides an approximation of the sharp crack energy 
$$
	\mathcal{F} ( \bu ) = \int_{\Omega \setminus J_u } W ( \strain (\bu) ) \, dx + G_c \mathcal{H}^1 ( J_u) , 
$$
where the displacement $\bu$ belongs to the space $SBD^2 (\Omega)$ and satisfies the non-interpenetration condition $\llbracket \bu \,\rrbracket \cdot \hat{\bn} = ( \bu^+ - \bu^- ) \cdot \hat{\bn} \ge 0$ on the crack set $J_u$. % (which represents the crack). 

In the sequel, we will prove that this convergence property holds in terms of $\Gamma$-convergence \cite{Braides98,DalMaso93} (as the internal length $\eps$ vanishes) at least in a reference one-dimensional case. Note that the scaling factor $c_\v$ in the energy $\F_\eps$ is needed to get the correct toughness value $G_c$ in the limit energy $\F$. We will indeed provide an explicit formula for $c_\v$, as a function of the parameter $\v$. %which is non trivial  due to the combination of the constraint $v \in [0,1]$ and the second order term $\Delta v$.  

\subsection{Evolution law}

We consider a quasi-static (rate-independent) evolution, driven by a boundary condition of the form $\bu = t \bg $ where $t \in [0,T]$ is a ``pseudo-time''  variable. As we consider linear elasticity, the energy at each time $t \in [0,T]$ can be written as 
$\F_\eps (t, \bu , v) = t^2 \E_\eps (\bu,v) + G_c \K_\eps (v)$ for $\bu \in \U = \{  \bu \in H^1(\Omega;\R^2) : \bu = \bg \text{ on }  \partial_D \Omega  \} $  and $v \in \V$. In the evolution we will further introduce the irreversibility constraint $\dot{v} \ge 0$. 

It is well known that the energy $\F_\eps$ is only ``separately convex'' i.e., while $\F_\eps ( t, \bu , \cdot) : \V \to [0,+\infty)$ and $\F_\eps ( t, \cdot , v) : \U \to [0,+\infty)$ are convex, the energy $\F_\eps (t, \cdot, \cdot ) : \U \times \V \to [0,+\infty)$ in general is not (jointly) convex. However, separate convexity is enough to characterize equilibrium configurations in terms of partial derivatives of the energy. More precisely, taking into account the irreversibility constraint, we will say that a configuration $(\bu , v)$ is in equilibrium at time $t \in [0,T]$ if 
$$
	\partial_{\bu} \F_\eps ( t, \bu , v ) [ \boldsymbol{\xi} ] = 0    , \quad  \partial_v \F_\eps ( t, \bu , v) [ \zeta ] \ge 0 
$$
for every admissible variation $\boldsymbol{\xi} = \boldsymbol{w} - \bu$ for $\boldsymbol{w} \in \U$ and $\zeta = z - v$ for $z\in \V$ with $z \ge v$; in other terms, by irreversibility the equilibrium of the phase-field is unilateral. The partial derivatives of the energy read 
\begin{gather*}
	\partial_{\bu} \F_\eps ( t, \bu , v ) [ \boldsymbol{\xi} ] = \int_{\Omega} \partial_{\strain} W_\eps ( v , \strain (\bu) ) : \strain ( \boldsymbol{\xi}) \, dx ,  \\
	\partial_v \F_\eps ( t, \bu ,v) [ \zeta ] = \int_{\Omega} \partial_v W_\eps  ( v , \strain (\bu) ) \, \zeta \, dx + \frac{G_c}{c_\v}  \int_{\Omega} \left(\eps^{-1} \zeta - 2 \eps \, \nabla v \cdot \nabla \zeta + 2 \v \, \eps^{3} \Delta v \, \Delta \zeta\right) dx .
\end{gather*}
Note that $\partial_{\strain} W_\eps ( v , \strain (\bu) )$ gives the phase-field stress 
$ \stress ( v , \bu) =  2 \psi_\eps (v) ( \mu \strain_d (\bu) + \kappa \strain^+_v (\bu) ) - 2 \kappa \strain_v^- (\bu) $, 
while $\partial_v W_\eps ( v , \strain(\bu))$ gives the crack driving force. 

Since the energy $\F_\eps$ is non-convex, in the quasi-static setting, we should expect both stable (or steady-state) and unstable (or catastrophic) propagation regimes, see, e.g.,~\cite{Maggiorelli24}. For the sake of simplicity, we confine ourselves to the former since the latter is delicate and should be better governed by dynamics or other rate-dependent effects. In the quasi-static setting, the stable regime is governed by the following system:
\begin{equation} \label{e.evol}
	\begin{cases}
	\partial_{\bu} \F_\eps (t, \bu (t) , v (t) ) [ \boldsymbol{\xi} ] = 0  \\
	\partial_v \F_\eps ( t, \bu (t) , v(t)) [ \zeta ] \ge 0 , \quad \partial_v \F_\eps ( t, \bu (t) , v(t)) [ \dot{v} (t) ] = 0 , \quad \dot{v} (t) \ge 0 ,
	\end{cases}
\end{equation} 
where $\boldsymbol{\xi} = \boldsymbol{w} - \bu (t) $ for $\boldsymbol{w} \in \U$ and $\zeta = z - v(t)$  for $z \in \V$ with  $z \ge v(t)$.  

In practice this evolution is obtained by solving incremental problems in a time-discrete setting. Let us consider the discrete times $t_k = k T / n_{\mathrm{ts}}$ where $n_{\mathrm{ts}}$ is the number of time steps and $k=0,...,n_{\mathrm{ts}}$. 
Starting from an equilibrium configuration $(\bu (t_{k-1}) , v (t_{k-1}) )$ at time $t_{k-1}$ we employ the staggered scheme 
$$
	\begin{cases}
		\bu_{m+1} \in \mathrm{argmin} \{  \F (t_k , \bu , v_m ) : \bu \in \U \} \\
		v_{m+1} \in \mathrm{argmin} \{  \F (t_k , \bu_m , v ) : v \in \V , \, v \ge v (t_{k-1}) \} \\
	\end{cases}
$$
to provide (technically, up to subsequences) an equilibrium configuration $( \bu (t_k ) , v(t_k))$ which satisfies the discrete counterpart of \eqref{e.evol}, i.e., 
$$
	\begin{cases}
	\partial_{\bu} \F (t, \bu (t_k) , v (t_k) ) [ \boldsymbol{\xi} ] = 0  \\
	\partial_v \F ( t, \bu (t_k) , v(t_k)) [ \zeta ] \ge 0 , \quad  \partial_v \F ( t, \bu (t_k) , v(t_k)) [ v (t_k) - v (t_{k-1} )  ] = 0 , \quad v (t_k) \ge v ( t_{k-1} ) .
	\end{cases}
$$
Here, we do not enter into the mathematical details about these evolutions. For second-order functionals, we refer the reader to \cite{MaggiorelliNegri24} for a comprehensive analysis of  \eqref{e.evol} including Griffith's criterion in the phase-field context.

%the technical details on the convergence of the time discrete evolution, which actually holds independently of the specific incremental problem (staggered or monolithic). \bl 
%We do not enter into the mathematical details about the convergence of the evolution, which will be the subject of future investigations. 

% !Tex root=AT1IV.tex

\section{Statement of the {\boldmath{$\Gamma$}}-convergence result}\label{s.gammaConv}

\subsection{One-dimensional setting \label{ss.1D}}

For the sake of simplicity, let us consider an interval $(a,b)$. Let $\U = H^1 (a,b)$  be the space of admissible displacements and $ \V =  H^2 ( (a,b)  ; [0,1]) $ the set of admissible phase-field functions. For $\eps>0$ and $\eta_\eps = o (\eps)$, let $\F_{\eps}  : L^1 (a,b) \times L^1(a,b) \to [0,+\infty]$ be the $AT_1$ functional defined by 
\begin{align} \label{e.Feps}
	\F_{\eps} ( u , v) = \begin{cases} 
	{\displaystyle \int_{(a,b)} \left(\kappa  \psi_\eps (v) | u' |_+^2 + \kappa |u'|_-^2  \right) dx +  \frac{G_c}{c_\v} \int_{(a,b)}  \left(\eps^{-1} v + \eps | v' |^2 + \v \, \eps^3
 %c \eps^3
 | v'' |^2  \right) dx }  & (u,v) \in \U \times \V  , \\ 	+\infty	& \text{otherwise,}
\end{cases}
\end{align}
with $\psi_\eps (v) = (v-1)^2 + \eta_\eps $.  In the sequel, we will provide the explicit value of $c_\v$ (as a function \aprev{of} $\rho$), while the parameter $\v>0$ will be of the form $1/ \c^2$ (this choice will be convenient to compute $c_\v$ in \S \ref{s.explicit}). In order to simplify the notation, we omit the dependence of the energy on the parameter $\v$. In \S \ref{convergence}, we will prove the convergence of \eqref{e.Feps} to the functional $\F : L^1 (a,b) \times L^1(a,b) \to [0,+\infty]$ defined as follows:
\begin{align} \label{e.Flim}
	\F ( u , v)  = \begin{cases} 
	 {\displaystyle \int_{(a,b) \setminus J_u}  \kappa  | u' |^2 \, dx + %c_{\c}\hspace{0.3pt} 
 G_c  \# J_u } & \text{if $u \in SBV^2 ( a,b)$ with $\jump{u} > 0$ and $v =0$ a.e.~in $(a,b)$}, \\
	+\infty	& \text{otherwise.}
	\end{cases}
\end{align}
Note that the energy $\F$ is finite if $u' \in L^2(a,b)$, the number of discontinuities is finite, and the non-interpenetration condition $\llbracket u \rrbracket > 0$ is satisfied in the jump points. 
The optimality constant $c_\v$ will be characterized in \S \ref{optimal}. 

Before stating our main convergence result, it is necessary to define the \textit{optimal profile problem}, which will ultimately provide the constant $c_\v$ for the $\Gamma$-convergence result. 
First, let us introduce the transition energy with unit internal length, $\K: \mathcal{W} \to \R$ given by 
$$
	\K( w) := \int_{\R_+} \left(w + | w' |^2 +  \c^{-2} | w''|^2 \right) dx , 
$$
where $\mathcal{W} = \{ w \in L^1 (\mathbb{R}_+) :  w' \in H^1  (\mathbb{R}_+) ,\, w(0)=1,\,w'(0)=0 \}$. %, i.e.~the space where $\K$ is finite. 
%\bb In order to simplify the notation, we omit de dependence of the eergy on the parameter $\gamma$.
\bl  
We now provide a result of existence and uniqueness of the optimal profile that will be subsequently characterized in \S \ref{optimal}.
\begin{proposition} \label{p.wop} There exists a unique 
\begin{equation} \label{e.optprof} 
	w_* \in \argmin \{ \K (w) : w \in \mathcal{W} \text{ such that }  \, 0 \le w \le 1 \} . 
\end{equation} 
\end{proposition}
\proof Take the  \textit{minimizing sequence} $\{\hat w_n \}_{n\in\mathbb N} \subset \mathcal{W}  $ such that $\K(\hat w_n)\to \inf\{\K(w) : w \in \mathcal{W} \text{ such that }  0 \le w \le 1\}$. 
Observe that such a sequence $\{\hat w_n \}_{n\in\mathbb N}$ is bounded in $ H^2(\R_+)$.  Indeed, being $0\leq w_n\leq 1$, $\|w_n\|^2_{H^2(\R_+)}=\int_{\R_+}\left(|w_n|^2+|w_n'|^2+|w_n''|^2\right)dx\leq C\int_{\R_+}\left(|w_n|+|w_n'|^2+\c^{-2}|w_n''|^2\right)dx<+\infty$ by definition of $\mathcal W$.
 Therefore, it exists a (non-relabeled) subsequence of $\{\hat w_n\}_{n\in\mathbb N}$ that weakly converges in $H^2(\R_+)$ to a certain $w_*\in H^2(\R_+)$. 
The set $\{w \in \mathcal{W}\,:\, 0\leq w\leq1\}$ is weakly closed (being convex and strongly closed) in $H^2(\R_+)$, whence $w_*$ belongs to this set. 
Now, since the functional $\K$ is weakly-lower semicontinuous (being strictly convex), 
we obtain that $w_*$ is indeed the unique minimum in definition \eqref{e.optprof}. 
\qed
\\In \S \ref{convergence} we will prove the following $\Gamma$-convergence result. 

\begin{theorem} \label{t.Gamma} Let $c_\v = 2 \K (w_*) $ and let $\eta_\eps = o (\eps)$. Then, $\F_{\eps}$ $\Gamma$-converges to $\F$ as $\eps \to 0^+$ with respect to the (strong) topology of $L^1 ( a,b ) \times L^1( a,b)$.
\end{theorem}
For numerical purposes it is fundamental to know the support of the optimal profile $w_*$ (in order to calibrate the mesh size) and the explicit value of the constant $c_\v$ (to estimate the effective toughness). To this purpose, in \S \ref{optimal} we will prove the following result. 

\begin{theorem} \label{t.param} The support of $w_*$ is the finite interval $[0,R_*]$ where $R_*$ is the unique solution of the non-linear equation
\begin{equation} \label{e.Rstar}
	\gamma R_* (1+\cosh(\gamma R_*) ) = 2(\c+1)\sinh(\gamma R_*) .
\end{equation} 
% %Moreover,  $w_*|_{[0,R_*]}=w_{R_*}$, where, %set $r_*= \c R_*$, $w_{R_*}(x)=z_{r_*}(\c x)$ and $z_{r_*}$ is the unique solution of \eqref{EL_r} with $r= \c R_*=:r_*$.
In the interval $[0,R_*]$ the optimal profile takes the form 
\begin{equation} \label{e.wstar}
	w_* (x) = a_3 (r_*)  \bigl(\cosh(\c x)-1\bigr)+ a_4 (r_*) \bigl(\sinh(\c x)-\c x\bigr)+\tfrac{1}{4}x^2+1, 
\end{equation}
where $r_* = \gamma R_*$ while $a_3$ and $a_4$ are given in \eqref{ci} and \eqref{cici}. The optimal constant $c_\v$ takes the form
\begin{equation} \label{e.cgamma}
	c_\v =\frac{2(1+\c)}{\gamma R_* }+\frac{(1+2\c)R_*}{2 \gamma} -\frac{R_*^3}{24} . 
\end{equation}
\end{theorem}

\begin{figure}
\centering
\includegraphics[width=.5\textwidth]{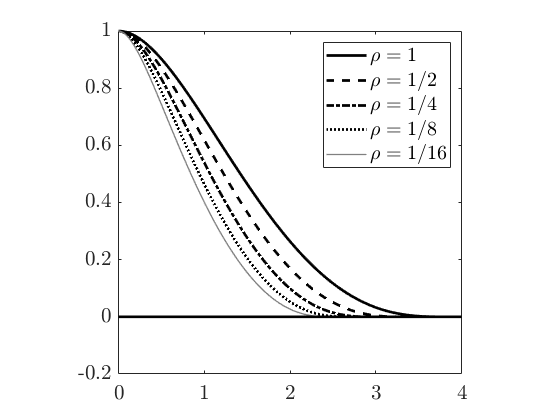}
\caption{The optimal profiles $w_*$ as a function of $\v$.}\label{fig:c}
\end{figure}

\begin{remark} 
As the optimal profile $w_*$ has compact support, it turns out that the phase-field ``interfaces'' of $AT_1$ are narrower than those of $AT_2$; indeed the optimal profile of $AT_2$ is supported in the whole $\mathbb{R}_+$  (see, for instance, \cite{N_CMAME20}). %{\color{red}This property will also allow to ... in the numerical experiments. } 
Moreover, the support of $w_*$ increases with $\v= 1 / \c^2$ (see Figure \ref{fig:c}). \end{remark}
In our numerical simulations, the weight $\v=1/\c^2$ appearing in the phase-field energy $\F_{\eps}$ takes values of the form $2^i$ for $i\in\{-4, ..., 4\}.$ %$\c\in\{1,\frac12, \frac14, \frac18,\frac{1}{16}\}$. \bl
From \eqref{e.Rstar}, by a fixed point algorithm, we obtain (with arbitrary precision) the explicit values of $R_*$ from which we can compute $c_\v$ from \eqref{e.cgamma}. Such values,  as a function of $\v=1/ \c^2$, are reported in \aprev{Table \ref{tab1}}. 
 
\begin{table}[h]  \centering\caption{\aprev{Tabulated $c_\v$ and $R_*$ values for $\v=\{\frac{1}{16},\frac{1}{8},\frac{1}{4},\frac{1}{2},1,2,4,8,16\}$}.}
		\vspace{-5pt}
\label{tab1}       
\begin{tabular}{p{0.5cm}p{1cm}p{1cm}p{1cm}p{1cm}p{1cm}p{1cm}p{1cm}p{1cm}p{1cm}}
\hline\noalign{\smallskip}
$\v$ & $\frac{1}{16}$  &  $\frac{1}{8}$  & $\frac{1}{4}$ & $\frac{1}{2}$    &$1$ & $2$ &     $4$ & $8$ & $16$\\
\noalign{\smallskip}
$R_*$ & $2.4998$  &  $2.7045$  & $ 2.9847$ & $3.3554$    &$3.8300$ & $4.4230$ &     $5.1515$ & $6.0387$ & $  7.1265$\\
\noalign{\smallskip}
$c_\v$ & $3.1615
$  & $3.3593$ & $3.6281$ & $ 3.9852$    &$4.4485$ & $5.0369$ &     $ 5.7715$ & $6.6714$ & $  7.7022$\\
\noalign{\smallskip} \hline
\end{tabular}
\end{table}
%\NOTE{AP: Inserito caption per uniformare con tabelle successive.}

%\NOTE{AP: Se possibile, toglierei il titolo in Figure 1 (la descrizione e' gia' nella caption) e anche se mi piacciono molto i colori, giocherei di spessori o markers per distinguere le curve (non so se c'e' qualche request dalla rivista per fare in modo che anche persone con daltonismo possano vederle...o anche semplicemente perche' solitamente si stampa in b/w e le curve poi risultano difficilmente distinguibili).}

\subsection{Isogeometric setting}

In this section, we provide the statement of the $\Gamma$-convergence result in the isogeometric setting. Consider a (physical) mesh $\mathcal{T}_h$ of size $h$ in the interval $[a,b]$ and let $\mathcal{S}_h \subset H^2 (a,b)$ be the space of splines on $\mathcal{T}_h$. We denote by $\U_h = \mathcal{S}_h$ the space for the displacement and by $\V_h$ the set for the phase field, i.e., $\V_h = \{ v_h \in \mathcal{S}_h : v_h \in [0,1] \}$. As $\mathcal{S}_h \subset H^2(a,b)$, the isogeometric functional $\F_{\eps,h}$ is simply the restriction of the functional $\F_\eps$ to the spline spaces $\mathcal{S}_h$. More precisely, $\F_{\eps,h}  : L^1 (a,b) \times L^1(a,b) \to [0,+\infty]$ is given by 
\begin{align} \label{e.Fepsh}
	\F_{\eps,h} ( u_h , v_h ) = \int_{(a,b)} \left(\kappa \psi_\eps (v_h) | u'_h |_+^2 +  \kappa |u'_h|_-^2  \right) dx +  \frac{G_c}{c_\v} \int_{(a,b)}  \left(\eps^{-1} v_h + \eps | v'_h |^2 + \v \, \eps^3 | v''_h  |^2  \right) dx  , 
\end{align}
for $(u_h,v_h) \in \U_h \times \V_h$ and $\F_{\eps,h}  ( u_h , v_h ) = +\infty$ otherwise. The convergence result is stated in the next Theorem.

\begin{theorem} \label{t.gammah} Let $c_\v = 2 \K (w_*) $, $\eta_\eps = o (\eps)$, and $h = o(\eps)$, then $\F_{\eps,h}$ $\Gamma$-converges to $\F $ as $\eps \to 0^+$ with respect to the (strong) topology of $L^1 ( a,b ) \times L^1( a,b)$.
\end{theorem}

We do not enter into the technical details and do not provide a proof, which follows from Theorem \ref{t.Gamma} together with the arguments of \cite{N_CMAME20}. %For the sake of simplicity we consider again the one-dimensional setting of \S \ref{ss.1D}. 
However, we highlight that the mesh size $h$ is much smaller than the internal length $\eps$; this is fundamental in order to obtain the correct $\Gamma$-limit and, in particular, the right value of toughness, choosing $h$ of order $\eps$ would result in a case similar to that of \cite{N_NFAO99}, where the toughness depends on the orientation of the mesh. 

%\subsection{2-dimensional setting}

%In this section we state the $\Gamma$-convergence result in the bidimensional setting. 

\section{Optimal profile: characterization of {\boldmath $w_*$}  and \boldmath{$c_\v$}  \label{optimal}}

\subsection{Auxiliary unconstrained problems} \label{uncost}
\begin{figure}[h] \label{f.Rstar}
\centering
\includegraphics[width=.5\textwidth]{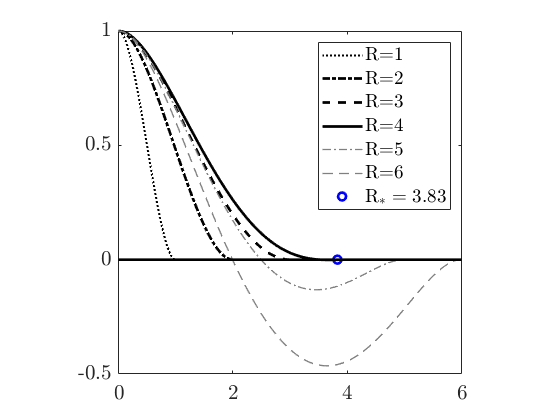} %dekstop 30 30 !
\caption{Solutions $w_R$ of \eqref{minR} as a function of $R$ for fixed $\v= 1$ and the value $R_*$ \aprev{corresponding to the maximum strictly positive value of the interval limit such that the profile $0 \le w_R \le 1$ (namely, $0 \le w_R \le 1$ for $R \le R_*$, while $w_R$ takes negative values for $R >R_*$).}}
\label{fig:R}
\end{figure}

Let us now introduce 
the localized energies with unit internal length, $\K_{R} : \mathcal{W}_R \to \R$ given by 
$$
	\K_{R}( w) := \int_{(0,R)} \left(w + | w' |^2 +  \c^{-2} | w''|^2 \right) dx,
$$
where $\mathcal{W}_R =\{w\in H^2(0,R)\,: \, w(0)=1,\,w'(0)=0,\, w(R)=0,\,w'(R)=0\}$. 
%, i.e., the space where $\K_{R}$ is finite. %When $R=+\infty$, we will simply \bb denote \bl  the energy and its space of definition $\K$ and $\mathcal{W}$  respectively, \bb as before. \bl 
%\bb To simplify the notation, we neglect the dependence of the energy on $\gamma$. \bl
Let us now consider the  minimization problem 
\begin{equation} \label{minR}
	w_R \in \argmin \{ \K_{R} (w) : w \in \mathcal{W}_{R} \} .
\end{equation} 
%\NOTE{AP: Stessa considerazione per la Figura 2 relativamente a colori, markers.}
We remark on the fact that, in contrast to \eqref{e.optprof}, in \eqref{minR}, it is not required that $0 \le w \le 1$. In general, solutions $w_R$ do not satisfy the constraint $0 \le w_R \le 1$ (see Figure \ref{fig:R}). However, solutions taking values in $[0,1]$ will be fundamental to characterize the optimal profile and will be called \textit{admissible solutions}; indeed, we will see that the support of the optimal profile $w_*$ is the finite interval $[0,R_*],$ where $R_* = \max \{ R > 0 : 0 \le w_R \le 1 \}$.
\\
Arguing as in Proposition \ref{p.wop}, it follows that the minimizer in \eqref{minR} is unique; moreover, it is characterized by 
\begin{equation}\label{EL_R}\begin{cases} 1-2 w_R''+\frac{2}{\c^2} w_R^{(4)}=0\qquad\text{on $(0,R)$} \\w_R(0)=1\\w_R'(0)=0 \\w_R(R)=0\\w_R'(R)= 0. \end{cases}\end{equation}
The solution of this ODE has the form 
\begin{equation}\label{wR}
w_R(x)=b_1+b_2x+b_3\cosh(\c x)+b_4\sinh(\c x)+\tfrac{1}{4}x^2.
\end{equation}
 Indeed $\frac{1}{4}x^2$ is a particular solution while the characteristic polynomial of the associated homogeneous equation has roots $\lambda_{1,2}=0$ with multiplicity two and $\lambda_{3,4}=\pm \c$ (an explicit calculation of the coefficients $b_i$ is provided in the sequel).

\subsection[Characterization of $w_*$]{Characterization of {\boldmath $w_*$}}

%In the following, we will characterize the value of $r_*$, and as a consequence that of $R_*=\frac{r_*}{\c}$. This will be crucial to provide the explicit form of $w_*$ and the explicit value of $c_\v(\c)$ to be used in numerical experiments. 

The next propositions establish the link between the (global constrained) profile $w_*$ and the local unconstrained profiles $w_R$ introduced in \S \ref{uncost}. %that solve \eqref{minR}. 

\begin{proposition} \label{p.eq} It exists $R_* > 0$ such that the support of $w_*$ is $[0,R_*]$ (the value of $R_*$ will be characterized in Proposition \ref{p.R*}). 
Moreover, $w_*|_{[0,R_*]}(x)= w_{R_*}(x)$, where $w_{R_*}$ is the unique solution of \eqref{EL_R} with $R=R_*$. 
% %Moreover,  $w_*|_{[0,R_*]}=w_{R_*}$, where, %set $r_*= \c R_*$, $w_{R_*}(x)=z_{r_*}(\c x)$ and $z_{r_*}$ is the unique solution of \eqref{EL_r} with $r= \c R_*=:r_*$.
%Moreover, set $r_*= \c R_*$, then $w_*|_{[0,R_*]}(x)=z_{r_*}(\c x)$, where $z_{r_*}$ is the unique solution of \eqref{EL_r} with $r=r_*$.
\end{proposition}
\bl 

\proof 
Notice that $w_*(x)<1$ for all $x>0$. Indeed, if it existed  $R'>0$ such that $w_*(R')=1$, then 
$w'_*(R') =0$  ($R'$ is a maximum point, since $0\leq  w_* \leq 1$, and, by Sobolev embedding, $w_*$ is of class $C^1$), hence 
$\hat w(x):=w_*(x+R')$ is an admissible competitor for $w_*$, indeed it belongs to $\mathcal W$, and satisfies  $0\leq \hat w\leq 1$. On the other hand, $\K(\hat w)<\K(w_*)$, which is in contradiction with \eqref{e.optprof}. 

We define $R_*:=\inf\{x\,:\,w_*(x)=0\}$ 
and observe that $R_*<+\infty$. Indeed, if this was not the case, $0<w_*(x)<1$ for all $x\in \R_+$ and thus $d\K(w_*)[\xi]=0$ for all $\xi\in C^\infty_0(\R_+)$, that means $w_*$ is the solution of:
\begin{equation}\label{inf}\begin{cases} 1-2w''+\frac{2}{\c^2} w^{(4)}=0\qquad\text{on $(0,+\infty)$} \\w(0)=1\\w'(0)=0.\end{cases}
\end{equation}
The solution of this ODE has the form $w(x)=b_1+b_2x+b_3\cosh(\c x)+b_4\sinh(\c x)+\frac{1}{4}x^2$. %Indeed $\frac{1}{4}x^2$ is a particular solution while the characteristic polynomial of the associated homogeneous equation has roots $\lambda_{1,2}=0$ with multiplicity two and $\lambda_{3,4}=\pm \c$. 
Such a function does not belong to $ L^1(\R_+)$ for any choice of the coefficients $b_i$ (as in any case it tends either to $+\infty$ or $-\infty$  for $x \to +\infty$), which is in contradiction with \eqref{e.optprof}.

Note that $w_*|_{[R_*,+\infty)}=0$, indeed, being $w_* (R_*) =0$ and $0 \le w_* \le 1$, $R_*$ is a minimum point for $w_*$ and thus $w_*'(R_*)  = 0$. As a consequence, the function 
$$
	\bar{w} (x ) = \begin{cases}  w_* (x) & x < R_* \\ 0 & x \ge R_*  \end{cases}
$$
belongs to the space $\mathcal{W}$ and satisfies $0 \le \bar{w} \le 1$. Hence, by minimality it coincides with the unique minimizer $w_*$ of $\mathcal{K}$.  
 We have hence proven the first part of the proposition.

Now, the function $w_*|_{(0,R_*)}$  takes values in $(0,1)$, hence $d\K(w_*)[\xi]=0$ for every $\xi\in C^\infty_0(0,R)$;  therefore, it solves %the Euler-Lagrange equation associated to $\K_{(0,R_*)}$ with suitable boundary conditions:
\begin{equation}\label{EL_R*}\begin{cases} 1-2w''+\frac{2}{\c^2} w^{(4)}=0\qquad\text{on $(0,R_*)$} \\w(0)=1\\w'(0)=0 \\w(R_*)=0\\w'(R_*)=0\end{cases}\end{equation}
where the last condition comes from the fact that $R_*$ is a minimum point. Hence, $w_*|_{[0,R_*]}(x)= w_{R_*}(x)$, where $w_{R_*}$ is the unique solution of  \eqref{minR} \bl for $R=R_*$. \qed 
%By a simple change of variable, the solution of the ODE \eqref{EL_R*} is equal to $z_{r_*}(\c x)$, where $z_{r_*}$ is the solution of \eqref{EL_r} with $r=r_*=\gamma R_*$. \qed 

\begin{proposition} \label{p.sup} It holds that $R_*=\max\{R \in (0,+\infty) :  w_R \in [0,1] \}$.
\end{proposition} 

\proof %
As proven in Proposition \ref{p.eq}, the function $w_*$ solves \eqref{EL_R*}  and by definition takes values in $(0,1)$. Therefore, $R_*$ belongs to the set $\{R \in (0,+\infty) \,:\,  w_R \in [0,1] \}$; hence, the latter set is non-empty and $R_* \leq \sup\{R \in (0,+\infty) \,:\, w_R \in [0,1] \}$.

To prove the inverse inequality, assume by contradiction that there exists $\bar{R} > R_*$ such that $ 0 \le w_{\bar{R}} \le 1$. 
By definition $w_{\bar{R}}$ 
%\eqref{EL_r} with $r=\bar r$ and equivalently, set $\bar{R}=\c \bar{r}$, $w_{\bar{R}}$ 
minimizes $\K_{\bar{R}}$ over the set  $\mathcal{W}_{\bar{R}}$ and it  takes values in $[0,1]$, so 
$$w_{\bar{R}}\in \argmin\{\K_{\bar{R}}(w) : w\in \mathcal W_{\bar{R}} \ \text{such that} \ w \in [0,1]  \}. $$Now, by Proposition \ref{p.eq}, $w_*|_{[R_*,+\infty)}=0$ 
and 
$w_*'(R_*)=0$; 
from this, since $R_* < \bar{R}$, it follows that $w_*(\bar{R})=w_*'(\bar{R})=0$. Therefore, by definition \eqref{e.optprof} of $w_*$, it holds that 
$$w_*|_{[0,\bar{R}]}\in \argmin\{\K_{\bar{R}}(w) : w \in \mathcal{W}_{\bar{R}} \ \text{such that} \  w \in[0,1] \}. $$  
By uniqueness of the minimizer, $w_*|_{[0,\bar{R}]}(x)=w_{\bar{R}}(x)
$.  This is a contradiction, since $w_* (x) = 0$ for $x >R_*$ while $w_{\bar{R}}(x)$ has the form  \eqref{wR} with $R=\bar{R}$. 
\qed

\begin{remark} \label{rem} \normalfont
  In Propositions \ref{p.eq} and \ref{p.sup}, it was shown that the solution of \eqref{e.optprof} is precisely given by  the \textit{admissible solution} $w_R$ with the largest possible support.  
\end{remark}

\subsection[Computation of $R_*$  and $c_\v$]{Computation of {\boldmath $R_*$}  and {\boldmath $c_\v$} \label{s.explicit}} 
In this section, we will provide several explicit calculations and, to this end, it is convenient to consider first the \textit{general auxiliary problem} of finding $z_r$ such that:
\begin{equation}\label{EL_r}\begin{cases} 1-2\c^2 z_r''+2\c^2z_r^{(4)}=0\qquad\text{on $(0,r)$} \\z_r(0)=1\\z_r'(0)=0 \\z_r(r)=0\\z_r'(r)=0.\end{cases}\end{equation}
%\bb In general, solutions $z_r$ of the above ODE do not satisfy the constraint $0 \le z_r \le 1$ (see Figure \ref{fig:R}). Solutions taking values in $[0,1]$ will be fundamental to characterize the optimal profile and will be called \textit{admissible solutions}. Indeed, we will see that the support of $w_*$ in the finite interval $[0,R_*]$ where $R_* = r_* / \gamma$ and $r_* = \max \{ r > 0 : 0 \le z_r \le 1 \}$. \bl
%
 %The values assumed by the solutions of \eqref{EL_r} depend on $r$ and we refer to 
%    the $z_r$ that solve \eqref{EL_r} and take values between $0$ and $1$ as \textit{admissible functions}.% In fact, % set $R:=\frac{r}{\c}$,  for such functions, %$w_{r/\c}(x):=z_r(\c x)$ 
  %  $w_R$ belongs to the set %whose minimum is looked for in \eqref{e.optprof}.
%featured in \eqref{e.optprof}.
%
%\noindent
%In Figure \ref{fig:R}, the solutions of \eqref{EL_r} (with $c=\frac{1}{\c^2}$ fixed) as $r$ varies are shown.
%
%It can be observed that, as pointed out in Remark \ref{rem}, $z_{r_*}$ is indeed the admissible solution with greatest possible support.
The solution of the ODE \eqref{EL_r} has the form
$$ z_r(y)=a_1(r)+a_2(r)y+a_3(r)\cosh(y)+a_4(r)\sinh(y)+\frac{y^2}{4\c^2} . $$ 
Indeed $y^2/ 4\c^2$ is a particular solution, while the characteristic polynomial of the associated homogeneous equation has roots $\lambda_{1,2}=0$ with multiplicity two and $\lambda_{3,4}=\pm 1$. Thus, the boundary conditions read:
\begin{equation*}
\begin{cases}
     a_1(r)+a_3(r)=1\\
    a_2(r)+a_4(r)=0\\
    a_1(r)+a_2(r)r+a_3(r)\cosh(r)+a_4(r)\sinh(r)+\frac{r^2}{4\c^2}=0\\
    a_2(r)+a_3(r)\sinh(r)+a_4(r)\cosh(r)+\frac{r}{2\c^2}=0 , 
\end{cases}
\end{equation*}
that can be rewritten as the linear system $\mathbf{Ax=b}$, with:
\begin{equation*}
    \mathbf{A}=\begin{pmatrix}
        \cosh(r)-1\quad &\sinh(r)-r\\\sinh(r) &\cosh(r)-1
    \end{pmatrix},\quad \mathbf{x}=\begin{pmatrix}
        a_3(r) \\a_4(r)
    \end{pmatrix},\quad \mathbf{b}=-\begin{pmatrix}
        1+\frac{r^2}{4\c^2}\\\frac{r}{2\c^2}
    \end{pmatrix}
\end{equation*} %and where 
together with the conditions $a_1(r)=1-a_3(r)$ and $a_2(r)=-a_4(r)$.
From direct calculations it follows that: 
\begin{align} \label{ci}
    &a_3(r)=\frac{(A_{22} b_1-A_{12} b_2)}{\det(A)}=-\frac{(r^2+4\c^2)\cosh(r)+r^2-4\c^2-2r\sinh(r)}{4\c^2(2-2\cosh(r)+r\sinh(r))}=1-a_1(r),\\ 
    &a_4(r)=\frac{(-A_{21} b_1+ A_{11} b_2)}{\det(A)}=\frac{(r^2+4\c^2)\sinh(r)+2r-2r\cosh(r)}{4\c^2(2-2\cosh(r)+r\sinh(r))}=-a_2(r). \label{cici}
\end{align} 
Given these coefficients, we can write $z_r$ in the following condensed form:\begin{equation} \label{zr}\displaystyle z_r(y)=a_3(r)\bigl(\cosh(y)-1\bigr)+a_4(r)\bigl(\sinh(y)-y)+\frac{y^2}{4\c^2}+1.\end{equation}
Clearly, the function $z_r$ is the unique minimizer of the energy 
\begin{equation}
	\mathcal{J}_r (z) := \int_{(0,r)} \left(z + \gamma^2 | z'|^2 + \gamma^2 | z''|^2 \right) dy 
\end{equation}
over the set $\mathcal{Z}_r = \{  z \in H^2(0,r) : z(0)=1,\, z'(0) = 0,\, z (r) = 0,\, z'(r) = 0 \}$. %\bb As before, we neglect the dependence of the energy on $\gamma$. \bl 
Finally, note that in general $z_r$ does not satisfy the admissibility constraint $z_r \in [0,1]$; however, solutions with $z_r \in [0,1]$ {will play an important role in the characterization of the optimal profile}. 
%\NOTE{\textbf{Luigi commenta}: riprendendo il commento di Ale dell'ultima revisione, c'era una parte commentata che ho scommentato ed evidenziato il magenta, può essere la parte mancante?}
\begin{remark} \normalfont 
Setting $r = \gamma R$, by a simple change of variable, we get $w_R(x) =  z_r (\c x) $, where $z_r$ is the unique solution of \eqref{EL_r}, from which we also get an explicit expression for $w_R$:
%and 
\begin{equation} \label{wRz}\displaystyle w_R(x)=z_r(\c x)=a_3(r)\bigl(\cosh(\c x)-1\bigr)+a_4(r)\bigl(\sinh(\c x)-\c x)+\tfrac{1}{4}x^2+1.\end{equation} 
Moreover, it turns out that $R_* = \frac{r_*}{\gamma}$ where $r_*=\sup\{r  \in (0,+\infty) \,:\,  z_r \in [0,1] \}$. %Therefore, characterizing $R_*$ is equivalent to characterize $r_*$. \bl 
\end{remark}

\begin{proposition} \label{p.R*}

The value of $r_*$ is characterized by the non-linear equation $z_{r_*}''(r_*)=0$, from which it follows the expression:\begin{equation} \label{r*}    r_*=\frac{2(\c+1)\sinh(r_*)}{1+\cosh(r_*)}.\end{equation} 

\end{proposition}
\proof We claim that for any fixed $r>0$, if $z_r \ge 0$ then $z_r \leq 1$.  \bl 
Indeed, if this was not the case, $\hat y=\text{argmax}\{z_r(y)\,:\,y\in(0,r)\}$ would be such that $z_r(\hat y)>1$. The function \begin{equation*}
    \hat{z}(y):=\begin{cases}
      \frac{z_r (y + \hat y )}{ z_r (\hat y)} \qquad & y \in [0, r-\hat y] 
      \\0 &x\in [r-\hat y, r] 
    \end{cases}
\end{equation*} is such that $\hat{z} \in \mathcal{Z}_{r}$ %,\,\hat w(0) = 1,\,\hat w'(0) = \hat w(R)=\hat w'(R)=0$ 
and $\mathcal{J}_r (\hat z) < \mathcal{J}_r (z_r)$, which is a contradiction since $z_r $ is the unique minimizer of $\mathcal{J}_r$ over $\mathcal{Z}_r$. 
Therefore, the admissible solutions $z_r$ of \eqref{EL_r}, i.e.,~those such that $z_r \in[0,1]$, 
must only satisfy the constraint $z_r \geq 0$ (since, as just observed, the condition $z_r \leq 1$ follows \textit{for free} by minimality). As a consequence, $r_*=\sup\{r\in (0,+\infty)\,:\, z_r\geq 0\}$.  

In the following, we will examine the values of $r>0$ for which the admissibility condition $z_r\geq 0$ is met. 
Firstly, let us prove that the condition $z_r''(r)\geq0$ is necessary to guarantee that 
$z_r\geq 0$. 
If $r>0$ is such that $z_r''(r)<0$, then it exists $y<r$ sufficiently close to $r$ such that $z_r(y)<0$. Indeed, by continuity, $z_r''(y)<0$ for $y$ sufficiently close to $r$ and, being $z_r'(r)=0$, 
$z_r$ is increasing near $r$ and thus 
$z_r(y)<z_r(r)=0$. Therefore, for such $r>0$, $z_r$ is not admissible. 

On the other hand, we claim that $z_r''(r)>0$ implies $z_r\geq 0$, and we prove it by contradiction. If this was not the case, $z_r (r') < 0$ where $r':=\argmin\{z_r(y)\,:\,y\in(0,r)\}$. 
We set $r'':=\max\{y\,:\,z_r(y)>0,\,z_r'(y)=0\}$ and observe that $r' < r''< r $; indeed, by the same reasoning as before,  $z_r''(r)>0$ implies
$z_r>0$ in a left neighborhood of $r$. 
Now, $z_r|_{(r',r)}$ is the unique solution of \begin{equation}\begin{cases} 1-2\c^2 z''+2\c^2z^{(4)}=0\qquad\text{on $(r',r)$} \\z(r')=z_r(r')\\z'(r')=0 \\z(r)=0\\z'(r)=0.\end{cases}\end{equation} 
We denote 
$$
	\J_{(r',r)} (z) = \int_{(r',r)} \left(z + \gamma^2 | z'|^2 + \gamma^2 | z''|^2 \right) dy 
$$
and observe that $z_r |_{(r',r)}$ is the unique minimizer of $\J_{(r',r)}$ over the set $\{ z \in H^2(r',r) : z(r')=z_r(r'), \,z'(r')=z(r)=z'(r)=0 \}$. 
We define \begin{align*}
 z_\alpha(y):=\begin{cases}
     z_r(r')+\beta_\alpha(z_r(y)-z_r(r')) \qquad &y \in (r',r'')\\
     \alpha z_r(y)& y \in (r'',r)
 \end{cases}
\end{align*}
for any $\alpha\in(0,1)$ and \aprev{set} $\beta_\alpha=\frac{\alpha z_r(r'')-z_r(r')}{z_r(r'')-z_r(r')}\in (0,1)$ so as to have continuity of $z_\alpha$. We remark that our choice of $r'$ and $r''$ guarantees that $z_\alpha(y)<z_r(y)$. Indeed, for all $x\in (r',r'')$, by definition of $r'$, $z_r(y)> z_r(r')$ and therefore $z_\alpha(y)=
     z_r(r')+\beta_\alpha(z_r(y)-z_r(r'))<z_r(r')+(z_r(y)-z_r(r'))=z_r(y)$. Moreover, $z_r>0$ in a left neighborhood of $r$, $z_r|_{[r'',r]}\geq 0$, ensuring that  $\alpha z_r(y)<z_r(y)$ on $(r'',r)$.
Now,
\begin{align*}
 z'_\alpha(y)=\begin{cases}
     \beta_\alpha z_r'(y) \qquad &y \in (r',r'')\\
     \alpha z_r'(y)& y \in (r'',r)
 \end{cases}
\end{align*}
and, being $r''$ a stationary point, %of maximum, 
$z'_\alpha(r''_-)=z'_\alpha(r''_+)=0$. It follows that $z_\alpha'\in H^1(r',r)$ and its weak derivative is \begin{align*}
 z''_\alpha(y)=\begin{cases}
     \beta_\alpha z_r''(y) \qquad &y \in (r',r'')\\
     \alpha z_r''(y)& y \in (r'',r).
 \end{cases}
\end{align*}
Consequently,  being $z_\alpha \le z_r$ and $\alpha,\,\beta_\alpha\in (0,1)$, we get $\mathcal{J}_{(r',r)}(z_\alpha)<\mathcal{J}_{(r',r)}(z_r|_{(r',r)})$, 
which is in contradiction with the minimality of $z_r$.

Summarizing what has just been proven, 
if  $z_r''(r)>0$
then $z_r$ is admissible, while if $z_r''(r)<0$ then 
$z_r$ is not admissible.
Therefore, in order to characterize the admissible values of $r>0$, one has to study the sign of %$w_R''(R)$ or equivalently that of
\begin{equation}\label{z''}
   z_r''(r)=-\frac{r^2(1+\cosh(r)-4r\sinh(r)+4(\c^2-1)(1-\cosh(r))}{4\c^2(2-2\cosh(r)+r\sinh(r))}.
\end{equation}
Equation \eqref{z''} is obtained by substituting the explicit expressions \eqref{ci} and \eqref{cici} of $a_i(r)$ in $z''_r(r)=a_3(r)\cosh(r)+a_4(r)\sinh(
r)+\frac{1}{2\c^2}$, as follows:
\begin{align*}
    z''_r(r)&=\frac{1}{4\c^2(2-2\cosh(r)+r\sinh(r))}\biggl(-(r^2+4\c^2)\cosh^2(r)-r^2\cosh(r)+ 4\c^2\cosh(r)+2r\sinh(r)\cosh(r)
    \\&\quad+(r^2+4\c^2)\sinh^2(r)+2r\sinh(r)-2r\cosh(r)\sinh(r)+2\bigl(2-2\cosh(r)+r\sinh(r)\bigr)\biggr)\\
    &=\frac{-(r^2+4\c^2)-r^2\cosh(r)+4\c^2\cosh(r)+4r\sinh(r)+4-4\cosh(r)}{4\c^2(2-2\cosh(r)+r\sinh(r))}\\
    &= - \frac{r^2(1+\cosh(r))-4r\sinh(r)+4(\c^2-1)(1-\cosh(r))}{4\c^2(2-2\cosh(r)+r\sinh(r))}.
\end{align*}
For $r>0$, the denominator in \eqref{z''} is always positive: indeed, by direct calculations, as $r\to 0$,  $F(r):=\frac{2(\cosh(r)-1)}{\sinh(r)}$ goes to zero with unit slope, while  $F''(r):=\frac{-2(\cosh(r)-1)^2}{\sinh^3(r)}$ is always negative and hence $r>\frac{2(\cosh(r)-1)}{\sinh(r)}$.
Thereby, to study the sign of $z_r''(r)$, it is sufficient to look for the roots of the numerator (as a quadratic polynomial in $r$), that are the fixed points of
 \begin{equation}
    \label{fixed} G_{1,2}(r):=\frac{2(1\pm \c)\sinh(r)}{1+\cosh(r)}.
\end{equation} 
Observe that it does not exist $r>0$ such that $r=G_2(r)=\frac{2(1-\c)\sinh(r)}{1+\cosh(r)}$: if $\c>1$, then $G_2(r)<0<r$, while if $\c<1$, then $G_2$ is strictly concave (in fact $G_2'(r)=\frac{2(1-\c)}{1+\cosh(r)}$ is decreasing) and, as $r\to 0$, it goes to zero with slope equal to $1-\c$, which is smaller than $1$, the slope of $r$. So, $G_2(r)<r$ for all $r>0$.
On the other hand, the positive fixed point of $G_1(r)=\frac{2(1+\c)\sinh(r)}{1+\cosh(r)}$ is unique. In fact, $G_1$ is strictly concave (as before, %being 
$G_1'(r)=\frac{2(1+\c)}{1+\cosh(r)}$ is decreasing) and as $r\to 0$, it goes to zero with slope equal to $1+\c$, which is greater than $1$, the slope of $r$. 
Therefore, the only admissible values of $r$ are those in between $0$ and the positive fixed point of $G_1$, denoted by $r_*$. 
\qed

\begin{remark} \label{w<1} \normalfont
Since the local profiles $w_R$ are defined in terms of $z_r$ by \eqref{wRz}, % In terms of the local profiles $w_R$,
the natural condition $w_R \le1$ is not actually needed in the optimal profile problem  \eqref{minR} since it follows \textit{for free} by minimality, as highlighted in the proof of Proposition \ref{p.R*}. 
   Moreover, the admissibility condition $ w_R \geq 0$ is equivalent to the condition $w_R''(R)\geq 0$.
\end{remark}

\begin{proposition}\label{wr*}
   The solution of \eqref{e.optprof} is $w_*(x):=  z_{r_*}  (\c x)$  where \begin{equation}
        z_{r_*}  (y)=\begin{cases}-\frac{2\c+1}{2\c^2}(\cosh(y)-1)+\bigl(\frac {r_*}{4\c^2}+\frac{1+\c}{\c r_*}\bigr)(\sinh(y)-y)+\frac{y^2}{4\c^2}+1 \qquad &[0,r_*]
         \\0 &(r_*,+\infty).  
       \end{cases}
   \end{equation}

\end{proposition} 
\proof  
Rearranging \eqref{ci}, we obtain that
\begin{align} \label{ciNEW}
    &a_3(r)=\frac{2\c^2-2\c^2\cosh(r)+r\sinh(r)}{2\c^2(2-2\cosh(r)+r\sinh(r))}-\frac{r^2\bigl(\cosh(r)+1\bigr)}{4\c^2(2-2\cosh(r)+r\sinh(r))}\\ \nonumber
    &a_4(r)=\frac {r}{4\c^2}+\frac{\sinh(r)}{2-2\cosh(r)+r\sinh(r)}.
\end{align} 
Substituting $r=r_* = \frac{2(\c+1)\sinh(r_*)}{1+\cosh(r_*)}$
from \eqref{r*}, we get:
\begin{align*} 
    &a_3(r_*)=\frac{2\c^2-2\c^2\cosh(r_*)+\frac{2(1+ \c)\sinh^2(r_*)}{1+\cosh(r_*)}}{2\c^2\bigl(2-2\cosh(r_*)+\frac{2(1+ \c)\sinh^2(r_*)}{1+\cosh(r_*)}\bigr)}-\frac{\frac{4(1+ \c)^2\sinh^2(r_*)}{1+\cosh(r_*)}}{4\c^2\bigl(2-2\cosh(r_*)+\frac{2(1+ \c)\sinh^2(r_*)}{1+\cosh(r_*)}\bigr)}
    \\&\qquad \quad =\frac{2(\c^2-\c^2\cosh^2(r_*)+(1+ \c)\sinh^2(r_*))}{4\c^2(1-\cosh^2(r_*)+(1+ \c)\sinh^2(r_*))}-\frac{(1+ \c)^2\sinh^2(r_*)}{2\c^2(1-\cosh^2(r_*)+(1+ \c)\sinh^2(r_*))}
    \\&\qquad \quad =\frac{-\c^2+\c+1}{2\c^3}-\frac{(1+ \c)^2}{2\c^3}=-\frac{2\c+1}{2\c^2};\\
    &a_4(r_*)=\frac{r_*}{4\c^2}+\frac{\sinh(r_*)}{2-2\cosh(r_*)+\frac{2(1+ \c)\sinh^2(r_*)}{1+\cosh(r_*)}}
    =\frac{r_*}{4\c^2}+\frac{\sinh(r_*)(1+\cosh(r_*))}{2(1-\cosh^2(r_*)+(1+ \c)\sinh^2(r_*))}
    \\&\qquad \quad =\frac{r_*}{4\c^2}+\frac{1+\cosh(r_*)}{2 \c\sinh(r_*)}=\frac{r_*}{4\c^2}+\frac{1+\c}{\c r_*}.
\end{align*} 
% and 
The thesis immediately follows from \eqref{zr}.  \qed

\begin{lemma}

%For any function $z_r$ that solves \eqref{EL_r}, 
The function $ w_* (x)=z_{r_*} (\c x)$ is such that: \begin{equation} \label{K}
   \K( w_*  )= \frac{1}{2\c} \int_{(0,r_*)} z_{r_*} (s)  \, ds +\c z_{r_*}'''(0). 
\end{equation}
As a consequence, the optimal constant is given by the explicit expression \begin{equation*}%\label{cw}
  c_\v =\frac{2(1+\c)}{r_*}+\frac{(1+2\c)r_*}{2\c^2} -\frac{r_*^3}{24 \c^3 },
\end{equation*}
or equivalently by \eqref{e.cgamma}. 
\end{lemma}
\proof By definition and by a simple change of variable, 
\begin{align*}
    \K ( w_* \bl ) &= \int_{\R_+} \left(w_* + | w_*' |^2 + \frac{1}{\c^2}| w_*''|^2  \right) dx =
     \frac{1}{\c}\int_{(0,r_*) } \left(z_{r_*} + \c^2| z_{r_*}' |^2 + \c^2| z_{r_*}''|^2 \right) dy 
    \\&=   \frac{1}{\c}\biggl(\int_{(0,r_*)} \left((\tfrac12+\tfrac12 )z_{r_*} - \c^2 z_{r_*}''z_{r_*} + \c^2 z_{r_*}^{(4)}z_{r_*} \right) dy +\c^2 z_{r_*}'''(0)\biggr)
    \\&
    =   \frac{1}{\c} \int_{(0,r_*)} \left(\tfrac12z_{r_*} +\tfrac12\bigl(1-2\c^2 z_{r_*}'' + 2\c^2 z_{r_*}^{(4)}\bigr)z_{r_*} \right) dy +\c  z_{r_*}'''(0)
    %\\&
    = 
    \frac{1}{2\c}\int_{(0,r_*)} z_{r_*} \, dy +\c z_{r_*}'''(0), 
\end{align*}
 where the third equality comes from integration by parts. Last equality simply follows from the fact that, since $z_r$ solves \eqref{EL_r}, the terms in parentheses have zero sum.

 \noindent

Now, in order to obtain the optimal constant $c_\v$, which is defined as $c_\v=2\K(w_*)$, it is sufficient to substitute the explicit expression of $z_{r_*}$, given in Proposition \ref{wr*}, into \eqref{K}. In detail, we start by computing: 
\begin{equation*}
z_{r_*}'''(0)=a_3(r_*)\sinh(0)+a_4(r_*)\cosh(0) =a_4(r_*).
\end{equation*}
Subsequently, we calculate:
\begin{align*} 
   \int_{(0,r_*)} z_{r_*}dy=a_3(r_*)\bigl(\sinh(r_*)-r_*\bigr)+a_4(r_*)\left(\cosh(r_*)-1-\frac{r_*^2}{2}\right)+\frac{r_*^3}{12\c^2}+r_*.
\end{align*}
According to \eqref{K}, adding these terms to $2\c^2 z'''_{r_*}(0)=2\c^2 a_4(r_*)$ and, using equation \eqref{r*}, one obtains the following expression for $2\c \K(w_*)$:
\begin{align*} 
   2\c \K(w_*)&= \int_{\R_+} z_{r_*} \, dy +2\c^2 %z_{r_*}'''(0)=
   a_4(r_*)=
   a_3(r_*)\bigl(\sinh(r_*)-r_*\bigr)+a_4(r_*)\left(\cosh(r_*)-1-\frac{r_*^2}{2}+2\c^2\right)+\frac{r_*^3}{12\c^2}+r_*
   \\&=a_3(r_*)\sinh(r_*)+\biggl(\frac {r_*}{4\c^2}+\frac{1+\c}{\c r_*}\biggr)\biggl((\cosh(r_*)+1)+2(\c^2-1)-\frac{r_*^2}{2}\biggr)+\frac{r_*^3}{12\c^2}+r_*(1-a_3(r_*))
   \\&=-\frac{(2\c+1)\sinh(r_*)}{2\c^2}+\biggl(\frac{(\c+1)\sinh(r_*)}{2\c^2(1+\cosh(r_*))}+\frac{1+\cosh(r_*)}{2\c\sinh(r_*)}\biggr)(\cosh(r_*)+1)
   \\&\quad +\biggl(\frac {r_*}{4\c^2}+\frac{1+\c}{\c r_*}\biggr)\biggl(2(\c^2-1)-\frac{r_*^2}{2}\biggr)
   +\frac{r_*^3}{12\c^2}+r_*(1-a_3(r_*))
   \\&=-\frac{(2\c+1)\sinh(r_*)}{2\c^2}+\frac{(\c+1)\sinh(r_*)}{2\c^2}+\frac{(1+\cosh(r_*))^2}{2\c \sinh(r_*)}+\frac{r_*(\c^2-1)}{2\c^2}+\frac{2(\c^2-1)(1+\c)}{\c r_*}\\&\quad-\frac{r_*^3}{8\c^2}-\frac{r_*(1+\c)}{2\c}
   +\frac{r_*^3}{12\c^2}+r_*(1-a_3(r_*))
   \\&=\frac{-\sinh^2(r_*)+1+\cosh^2(r_*)+2\cosh(r_*)}{2\c \sinh(r_*)}+\frac{2(\c^2-1)(1+\c)}{\c r_*} 
   \\&\quad+r_*\biggl(\frac{\c^2-1}{2\c^2}-\frac{1+\c}{2\c}+1+\frac{2\c+1}{2\c^2}\biggr)-\frac{r_*^3}{24\c^2}
   \\&=\frac{1+\cosh(r_*)}{\c\sinh( r_*)}+\frac{2(\c^2-1)(1+\c)}{\c r_*}+r_*\frac{2\c+1}{2\c} -\frac{r_*^3}{24\c^2}
   \\&=\frac{2(1+\c)}{\c r_*}+\frac{2(\c^2-1)(1+\c)}{\c r_*}+r_*\frac{2\c+1}{2\c} -\frac{r_*^3}{24\c^2}
   \\&=\frac{2(1+\c)\c}{ r_*}+r_*\frac{2\c+1}{2\c} -\frac{r_*^3}{24\c^2}.
\end{align*}
Therefore, the optimal constant is given by:
\begin{equation*} 
       c_\v=2\K(w_*)=\frac{2(1+\c)}{ r_*}+r_*\frac{2\c+1}{2\c^2} -\frac{r_*^3}{24\c^3}
\end{equation*}
and the proof is concluded.  \qed

\section{Proof of the $\boldsymbol{\Gamma}$-convergence result \label{convergence}}

For $\eps_n \to 0^+$ we employ the notation $\F_n = \F_{\eps_n}$ and $\psi_n = \psi_{\eps_n}$.

\subsection{Liminf inequality}

%We first prove this lemma on the optimal profile. 
%
%\begin{lemma} \label{l.3.3} Let $S_n \to +\infty$ and $w_n \in H^2(0,S_n ; [ 0,1])$ such that 
%$$
%	w_n (0) \to 1 , \quad w'_n (0) = 0 , \quad w_n ( S_n) \to 0, \quad w'_n (S_n) \to 0 ,
%$$
%then $\liminf_{n \to +\infty} \mathcal{K}_{S_n} ( w_n ) \ge \K (w_*)$.
%\end{lemma}
%
%\proof {\tt extend on $(S_n , S_n +1)$ and control the energy in $(S_n , S_n +1)$ with with the b.c. in $R_n$, then extend in $\R_+$ } \qed

\begin{lemma} \label{l.liminf} Let $(u_n , v_n ) \to (u,v)$ in $L^1 (a,b) \times L^1 (a,b)$ such that $\liminf_{n \to +\infty} \F_n (u_n , v_n ) < +\infty$ then $ u \in SBV^2 (a,b)$ with $\llbracket u \rrbracket >0$, $v=0$ a.e.~in $(a,b)$ and 
\begin{align}
	c_\v \# J_u  & \le \liminf_{n \to +\infty} \int_{(a,b)}  \left(\eps_n^{-1} v_n + \eps_n | v_n'|^2 + \v\,  \eps_n^3 | v_n'' |^2 \right) dx , % \I_n (u_n , z_n ) 
	\label{e.liminfsurf} \\
	\int_{(a,b) \setminus J_u} | u' |^2 \, dx & \le \liminf_{n \to +\infty} \int_{(a,b)} \left(\psi_n (v_n) | u' _n |_+^2 + | u' _n |_-^2 \right) dx . \label{e.liminfbulk}
\end{align}
\end{lemma}

\proof {\bf I.} Clearly, it is enough to consider the case when the liminf in \eqref{e.liminfsurf}-\eqref{e.liminfbulk} are finite. Upon extracting a subsequence (non relabeled) 
it is not restrictive to assume that the liminf in \eqref{e.liminfsurf} and \eqref{e.liminfbulk} is actually a limit. 
%_{n \to +\infty} \F_n (u_n , v_n ) = \lim_{n \to +\infty} \F_n (u_n , v_n )$ 
Thus, $u_n \in \U = H^1(a,b)$  and $v_n \in \V = H^2(a,b ; [0,1])$ for every $n \in \mathbb{N}$. In particular, $v_n \ge 0$ and % there exists $C>0$ such that 
$$
	\int_{(a,b)}  \left(\eps_n^{-1} | v_n | + \eps_n | v_n'|^2 \right) dx \le \int_{(a,b)}  \left(\eps_n^{-1} v_n + \eps_n | v_n'|^2 + \v \, \eps_n^3 | v_n'' |^2 \right) dx .
$$
It follows that
$$
	\mathcal{H}_n (u_n , v_n) = \int_{(a,b)} \psi_n ( v_n ) | u'_n |^2 \, dx +  \int_{(a,b)}  \left(\eps_n^{-1} | v_n | + \eps_n | v_n'|^2 \right) dx 
	\le C \F_n (u_n , v_n ) \le C' . 
$$
The functional $\mathcal{H}_n$ introduced above is a (classical second order) elliptic approximation of  the Mumford-Shah functional; therefore, we can employ the compactness arguments of \cite[Theorem 3.15]{Braides98} from which it follows that $u \in SBV^2(a,b)$, $v=0$    a.e. on $(a,b)$, and that $u_n \weakto u$ in $H^1( (a,b) \setminus J^\delta_u)$ for every $\delta >0$, where $J^\delta_u = \{ x : \mathrm{dist} ( J_u , x) \le \delta \}$.

Let us check that $\jump{u} >0$. Assume by contradiction that $\llbracket u (s) \rrbracket < 0$ for some $s \in J_u$, i.e., $u^+ (s) < u^- (s)$. Let $\eta = \tfrac14 | \jump{u}  |$. For every $\delta >0$ (sufficiently small) there exist $ s - \delta < s' <  s'' < s+\delta $ such that $u_n (s') \to u (s') >  u^- (s) -  \eta$ and $u_n (s'') \to u (s'') <  u^+ (s) + \eta$. If $n \gg 1$  it turns out that 
$$
	\int_{(s'\!,s'')}  | u'_n |_-  \,dx \ge | u_n (s'') - u_n (s') |  =  u_n (s') - u_n (s'') \ge  u^- (s) - u^+ (s) - 2 \eta \ge \tfrac12 | \jump{u(s)}  | .
$$
Hence, by Jensen's inequality and by the boundedness of the energies
$$
	\tfrac14 | \llbracket u (s) \rrbracket |^2 \le \left( \int_{(s'\!,s'')}  | u'_n |_-  \,dx  \right)^2 \le (s'' - s') \int_{(s'\!,s'')}  | u'_n |_-^2  \,dx \le 2 \delta C .
$$
By the arbitrariness of $\delta$ it follows that $ | \llbracket u (s) \rrbracket | = 0$.

{\bf II.} In order to prove \eqref{e.liminfsurf} we ``localize'' around jump points. Let $J_u = \{ s_i : i =1,...,N \}$ and let $0 < \delta \ll 1$ such that the intervals $J_i^\delta = [s_i-\delta , s_i +\delta]$ are disjoint and contained in $(a,b)$.
Clearly 
$$
	 \int_{(a,b)}  \left(\eps^{-1} v_n + \eps | v_n' |^2 + \v\, \eps^3 | v_n'' |^2  \right) dx \geq 
	 %\I_n (u_n , z_n) \ge %\tfrac12 \int_{I^\delta}  (z_n^2 + \eta_n ) | u_n'|^2 \, ds + 
         \sum_{i=1}^N \int_{J^\delta_i}  \left(\eps^{-1} v_n + \eps | v_n' |^2 + \v\, \eps^3 | v_n'' |^2  \right) dx ,
$$
and thus it is enough to prove that for every $i=1,...,N$ it holds 
\begin{equation} \label{e.glinf2}
	 \liminf_{n \to +\infty} \int_{J^\delta_i}   \left(\eps^{-1} v_n + \eps | v_n' |^2 + \v \, \eps^3 | v_n'' |^2  \right) dx \ge c_\v . 
\end{equation}
% by arbitrarity of $\delta$.  
We follow the outline of \cite[Proposition 4.1]{N_CMAME20}. Upon extracting a subsequence (not relabelled) it is not restrictive to assume that $v_n \to 0$ a.e.~in $(-\delta,\delta)$.  %and recall that $v_n \in C^1( - \delta , \delta)$. 
Moreover, for simplicity of notation we assume that the interval $J_i^\delta$ is of the form $[-\delta , \delta]$, i.e.~that $0 \in J_u$.

We claim that there exists $ \hat{x}_n \in (-\delta , \delta) $ with $\hat{x}_n \to 0 $ such that $v_n (\hat{x}_n) \to 1^-$ and $v'_n(\hat{x}_n) =0$. Let $\delta' 
> 0$ (sufficiently small) such that $v_n (-\delta') \to 0$ and  $v_n ( \delta') \to 0$.  Let $x_n \in \mathrm{argmax} \{ v_n (x) : x \in [-\delta' , \delta'] \}$. Let us check that $ \lim_{n \to +\infty} v_n (x_n) = 1$. 
%such that 
Assume by contradiction that there exists a subsequence (non relabeled) such that $\max \{ v_n (x) : x \in [-\delta' , \delta']  \} \le C < 1$. Then $\psi_n ( v_n (x) ) > C' > 0$ for every $x \in [-\delta' , \delta'] $ and thus 
$$
		\F_n ( u_n , v_n) \ge C' \int_{(-\delta' , \delta') } | u'_n|^2 \, dx . % \ge C  \int_{  (-\delta^  {\flat} , \delta^  {\flat})}  | u'_n|^2 \, ds .
$$
Hence, $\{ u_n \}$ is bounded in $H^1( -\delta' , \delta' )$ and its limit $u$ belongs to $H^1 ( -\delta' , \delta' )$, which contradicts $0 \in J_u$. As $v_n ( -\delta') \to 0$ and $v_n (\delta') \to 0$ the maximizer $x_n$ actually belongs to the open interval $( -\delta' , \delta')$ and thus $v'_n (x_n) =0$. The claim follows by the arbitrariness of $\delta'$.

For simplicity, let us assume that $\hat{x}_n =0$ and define the rescaled functions $w_n (s) =  v_n ( \eps_n  s) $ for $ | s | \le R_n = \delta / \eps_n \bl$. Then, by a change of variable, % $s= \eps_n r + s_n^\flat$ we have 
\begin{align*}
	\int_{J_i^\delta} \left(\eps_n^{-1} v_n  + \eps_n | v_n'|^2 +  \v \, \eps_n^3 | v_n'' |^2 \right) dx 
	=
	\int_{(-R_n, R_n)} \left(w_n  + | w_n'|^2 +  \v \, | w_n'' |^2 \right) dx . 
\end{align*}
Note that $ w_n ( 0 ) \to 1^-$ and $w'_n (0) = 0$. Clearly, for $n \gg 1$ and every $M \in \mathbb{N}$ we have 
$$
	\int_{(0,M)} \left(w_n  + | w_n'|^2 +  \v \, | w_n'' |^2 \right) dx \le C 
$$
and thus $w_n$ is a bounded sequence in $H^2 (0,M)$. By a diagonal argument it follows that there exists a subsequence (non relabelled) and a function $w \in H^2_{\text{\sl loc}} ( 0, + \infty)$ such that $w_n \weakto w $ in $H^2(0,M)$  for every $M \in \mathbb{N}$. As a consequence
%$$
%	\liminf_{n \to +\infty} \int_{(0, R_n)} \left(w_n  + | w_n'|^2 +   \v \, | w_n'' |^2 \right) dx 
%	\ge
%	\liminf_{n \to +\infty} \int_{(0, M)} \left(w_n  + | w_n'|^2 + \v \,  | w_n'' |^2 \right) dx \ge
%	\int_{(0, M)} \left(w  + | w '|^2 +  \v \, | w'' |^2 \right) dx .
%$$
\begin{align*}
	\liminf_{n \to +\infty} \int_{(0, R_n)} \left(w_n  + | w_n'|^2 +   \v \, | w_n'' |^2 \right) dx 
	&\ge
	\liminf_{n \to +\infty} \int_{(0, M)} \left(w_n  + | w_n'|^2 + \v \,  | w_n'' |^2 \right) dx\\
    &\ge \int_{(0, M)} \left(w  + | w '|^2 +  \v \, | w'' |^2 \right) dx . 
\end{align*}
Taking the supremum with respect to $M$ yields
$$
	\liminf_{n \to +\infty} \int_{(0, R_n)} \left(w_n  + | w_n'|^2 +\v \,  | w_n'' |^2 \right) dx 
	\ge
	\int_{(0, +\infty)} \left(w  + | w '|^2 +\v \,   | w'' |^2 \right) dx = \mathcal{K} (w) \geq \mathcal{K} (w_*) . 
$$
Arguing in the same way in the interval $(-\infty, 0)$ we finally get \eqref{e.glinf2}, since $c_\v = 2 \mathcal{K} (w_*)$.

{\bf III.} It remains to prove \eqref{e.liminfbulk}. First, note that $\psi_n (v_n) - \eta_n = (v_n -1)^2  \to 1$ strongly in $L^2(a,b)$ (for instance by dominated convergence) and remember that $u_n \weakto u$ in $H^1( (a,b) \setminus J^\delta_u)$ for every $\delta >0$. As a consequence $(\psi_n ( v_n ) - \eta_n)^{1/2}  u'_n \weakto u'$ in $L^2 ( (a,b) \setminus J^\delta_u)$ for every $\delta >0$ (by weak-strong convergence). Therefore, 
$$
	\liminf_{n \to +\infty}  \int_{(a,b) \setminus J_u^\delta}  \left(\psi_n ( v_n ) | u' _n |_+^2 + | u'_n|_-^2 \right) dx  \ge
	\liminf_{n \to +\infty}  \int_{(a,b) \setminus J_u^\delta}  ( \psi_n ( v_n )  - \eta_n )| u' _n |^2  \, dx  \ge 
	 \int_{(a,b) \setminus J_u^\delta} | u' |^2  \, dx   .
$$
Taking the supremum with respect to $\delta>0$ yields \eqref{e.liminfbulk}.  \qed 

From Lemma \ref{l.liminf} and by a standard supremum of measures argument  \cite[Proposition 1.16]{Braides98} the following $\Gamma$-liminf estimate holds. 

\begin{lemma} \label{l.liminf} Let $(u_n , v_n ) \to (u,v)$ in $L^1 (a,b) \times L^1 (a,b)$ such that $\liminf_{n \to +\infty} \F_n (u_n , v_n ) < +\infty$ then $ u \in SBV^2 (a,b)$ with $\jump{u} >0$, $v=0$ a.e.~in $(a,b)$, and 
$$
	\int_{(a,b) \setminus J_u} \kappa | u' |^2 \, dx  + G_c  \# J_u  \le \liminf_{n \to +\infty} \F_n ( u_n , v_n ) . 
$$
\end{lemma}

\subsection{Limsup inequality \label{Gsup}}

\begin{proposition} \label{p.Glsup} Let $u \in SBV^2(a,b)$ with $\jump{u} > 0$. 
There exist $u_n \in \U = H^1 (a,b)$ and $v_n \in \V = H^2 (a,b ; [0,1])$ such that $u_n \to u$ in $L^1 (a,b)$, $v_n \to 0$ in $L^1 (a,b)$ and 
\begin{align} 
	 \int_{(a,b)} \left(\psi_n (v_n)  | u'_n |_+^2 + |u'_n|_-^2  \right) dx  &  \ \to \  \int_{(a,b) \setminus J_u} | u' |^2 \, dx , \label{e.18} \\
	  \int_{(a,b)}  \left(\eps_n^{-1} v_n^2 + \eps_n | v'_n |^2 + \v \, \eps^3_n | v''_n |^2  \right) dx 
	 & \ \to \ c_\v \# J_u .  \label{e.19}
\end{align}
\end{proposition} 

\proof It is not restrictive to consider the case of a single crack, i.e., ~$J_u = \{ s \}$. Let $\eta_n \ll \delta_n \ll \eps_n$. As $\jump{u}>0$ for $n \gg 1$ we have $u (s+\delta_n) > u (s - \delta_n)$. Let us define the functions
\begin{gather*}
	u_n (x) = 
	\begin{cases}
		u (x) & \text{for } x < s - \delta_n \\
		t u( s +\delta_n) + (1-t) u ( s -\delta_n) & \text{for }  s - \delta_n \le x \le s + \delta_n \text{ and } t = (x - s + \delta_n) / 2 \delta_n  \\
		u(x) & \text{for }  x > s + \delta_n , 
	\end{cases}
	\\
	v_n (x) = 
	\begin{cases}
		w_* ( (s-\delta_n - x )/\eps_n) & \text{for } x < s - \delta_n \\
		1 \bl & \text{for } s- \delta_n \le x \le s+ \delta_n  \\
		w_*((x-s-\delta_n)/\eps_n) & \text{for }  x > s + \delta_n . 
	\end{cases}
\end{gather*}
Note that $u_n \in H^1(a,b)$ and that $u'_n >0$ in $(s-\delta_n , s+\delta_n)$. Note also that $v_n \in H^2 (a,b)$ since $v_n(0)=1$ and $ v_n' (0)=0$ and that $\mathrm{supp} ( v_n) \subset (s -\delta_n -\eps_n  R_* , s + \delta_n + \eps_n R_* )$. 
It is easy to check that $u_n \to u$ in $L^1(a,b)$ and $v_n \to 0$ in $L^1(a,b)$. 
Let us write
\begin{align*}
	 \int_{(a,b)} \left(\psi_n (v_n)  | u'_n |_+^2 + |u'_n|_-^2  \right) dx 
	 & = 
	 \int_{(a,s-\delta_n)} \left(\psi_n (v_n)  | u' |_+^2 + |u' |_-^2  \right) dx +
	 \int_{(s-\delta_n,s+\delta_n)} \eta_n  | u'_n |_+^2 \, dx \, + \\
	& \quad  \int_{(s+\delta_n,b)} \left(\psi_n (v_n)  | u' |_+^2 + |u' |_-^2  \right) dx . 
\end{align*}
As $v_n \to 0$ we have $\psi_n (v_n) \to 1$ a.e.~in $(a,b)$. By dominated convergence it follows that 
$$
	\int_{(a,s-\delta_n)} \left(\psi_n (v_n)  | u' |_+^2 + |u' |_-^2  \right) dx \ \to \ \int_{(a,s)}  \left(| u' |_+^2 + |u' |_-^2  \right) dx = \int_{(a,s)}  | u' |^2 \, dx . 
$$
The same reasoning can be applied to the interval $(s+\delta_n,b)$. Moreover, 
$$
	\int_{(s-\delta_n,s+\delta_n)} \eta_n  | u'_n |_+^2 \, dx = \eta_n ( u ( s + \delta_n ) - u( s - \delta_n ) )^2 / 2 \, \delta_n \to 0 ,
$$
since $\eta_n = o (\delta_n)$ and $u ( s + \delta_n ) - u( s - \delta_n ) \to \llbracket u (s) \rrbracket$. Thus \eqref{e.18} is proved.

Let us write 
\begin{align*}
	\int_{(a,b)}  \left(\eps_n^{-1} v_n^2 + \eps_n | v'_n |^2 + \v\, \eps^3_n | v''_n |^2  \right) dx  
	& = \int_{(a, s-\delta_n)}  \left(\eps_n^{-1} v_n^2 + \eps_n | v'_n |^2 + \v\, \eps^3_n | v''_n |^2  \right) dx  \ + \\ 
	& \quad \int_{(s-\delta_n,s+\delta_n)}  \eps_n^{-1} \, dx \ + \\
	& \quad \int_{(s+\delta_n,b)}  \left(\eps_n^{-1} v_n^2 + \eps_n | v'_n |^2 + \v\, \eps^3_n | v''_n |^2  \right) dx  . 
\end{align*}
As $\delta_n = o (\eps_n)$ we have 
$$
	\int_{(s-\delta_n,s+\delta_n)}  \eps_n^{-1} \,  dx = 2 \, \delta_n \eps_n^{-1} \to 0 . 
$$
Moreover, by a simple change of variable, 
$$
	\int_{(s+\delta_n,b)}  \left(\eps_n^{-1} v_n^2 + \eps_n | v'_n |^2 + \v \, \eps^3_n | v''_n |^2  \right) dx = \int_0^{R_*} \left(w_* + | w_*' |^2 + \v \, | w_*''|^2 \right) dx = \mathcal{K}(w_*) . 
$$
We can proceed in the same way in the interval $(a , s- \delta_n)$.  The proof is concluded. \qed

% !Tex root=AT1IV.tex

\section{Numerical tests}
\label{sec:numerical test}
In this section, we provide an extensive comparison between $AT_1$, second- and fourth-order models. We can view the proposed fourth-order model (see §\ref{optimal}) as an extension of the second-order phase-field functional,  weighted by the coefficient $\v$, which we initially set equal to 1 (a study of the dependence on $\rho$ is given in \S \ref{sec: sensitivity study}).% coherently with the discussion carried out in the previous sections. 
Consequently, we set $c_\rho = 4.4485$ and $R_* = 3.83$ and we analyze the numerical performance of our $AT_1$ fourth-order model, considering several benchmarks. First, a pure tensile test is carried out to assess the elastic limit and its relative error with {respect} to the theoretical expected value. In this case, we compare the results only for second- and fourth-order $AT_1$ functionals (since $AT_2$ does not define an elastic limit).

\begin{figure}[h!!]
	\begin{center}
		\begin{subfigure}[b]{0.75\textwidth}
			\includegraphics[width=10cm]{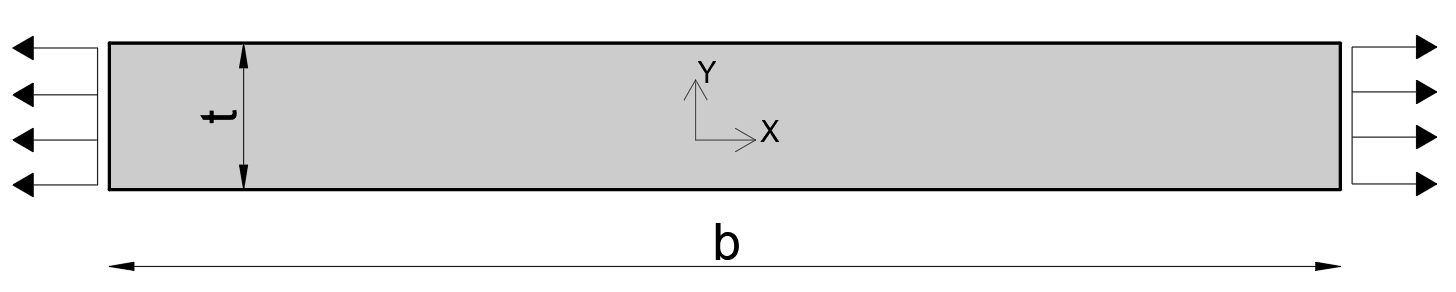}
			\caption{Pure traction test.}\label{fig:geometry_bc_a}
			\vspace{0.25cm}
		\end{subfigure}
		\begin{subfigure}[b]{0.3\textwidth}
			\includegraphics[width=4.35cm]{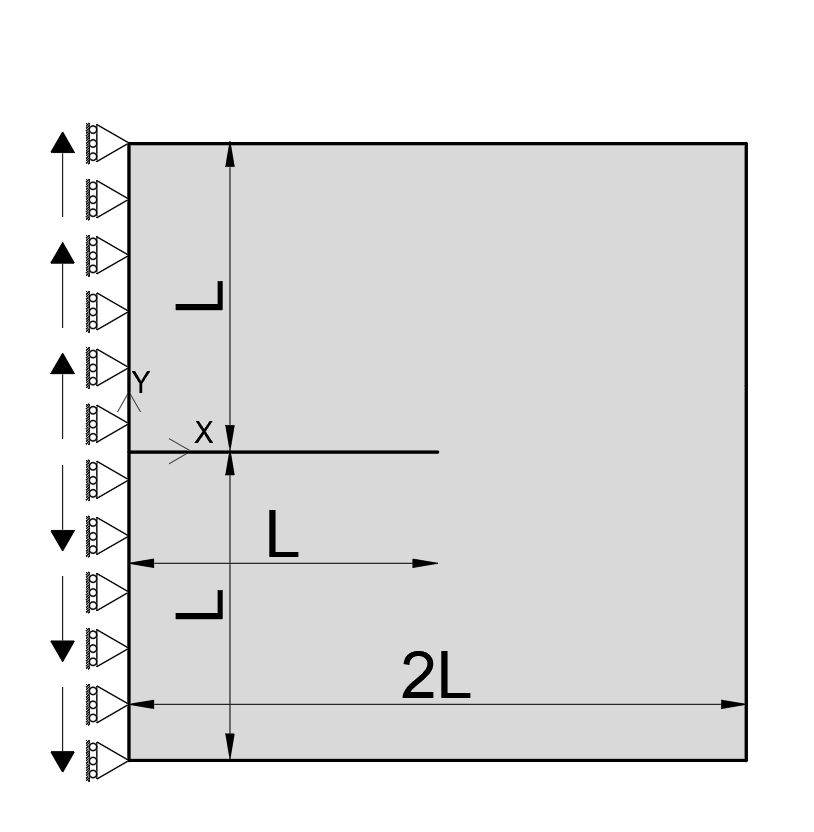}
			\caption{DCB.}\label{fig:geometry_bc_b}
			\vspace{3cm}
		\end{subfigure}
		%	\hfill
		\begin{subfigure}[b]{0.3\textwidth}
			\includegraphics[width=4.35cm]{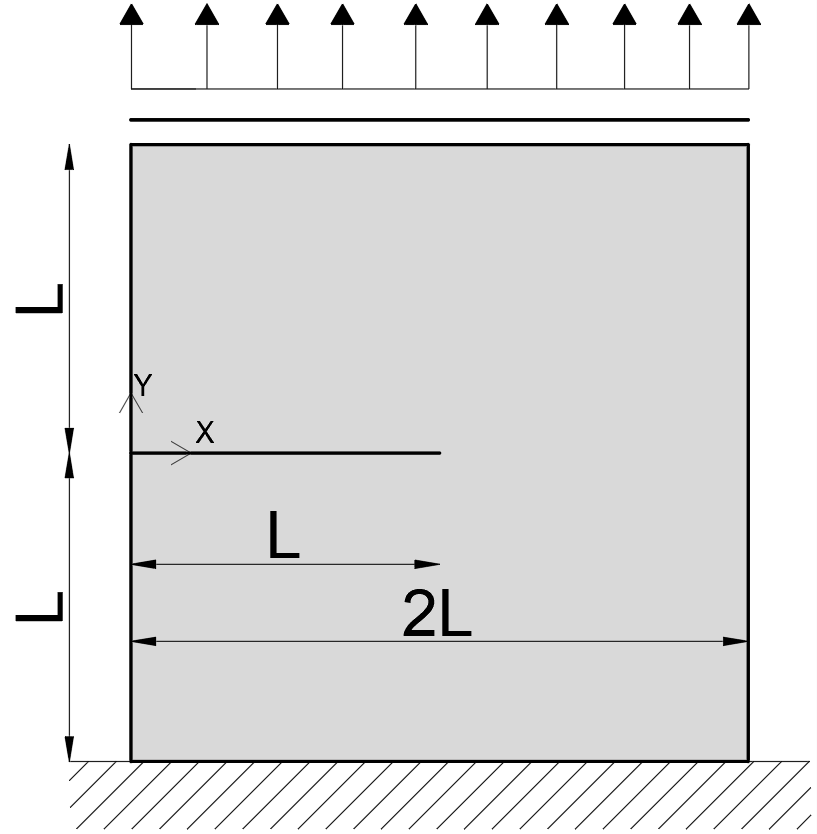}
			\caption{SEN tensile test.}\label{fig:geometry_bc_c}
			\vspace{3cm}
		\end{subfigure}
		%   \hfill
		\begin{subfigure}[b]{0.3\textwidth}
			\includegraphics[width=4.35cm]{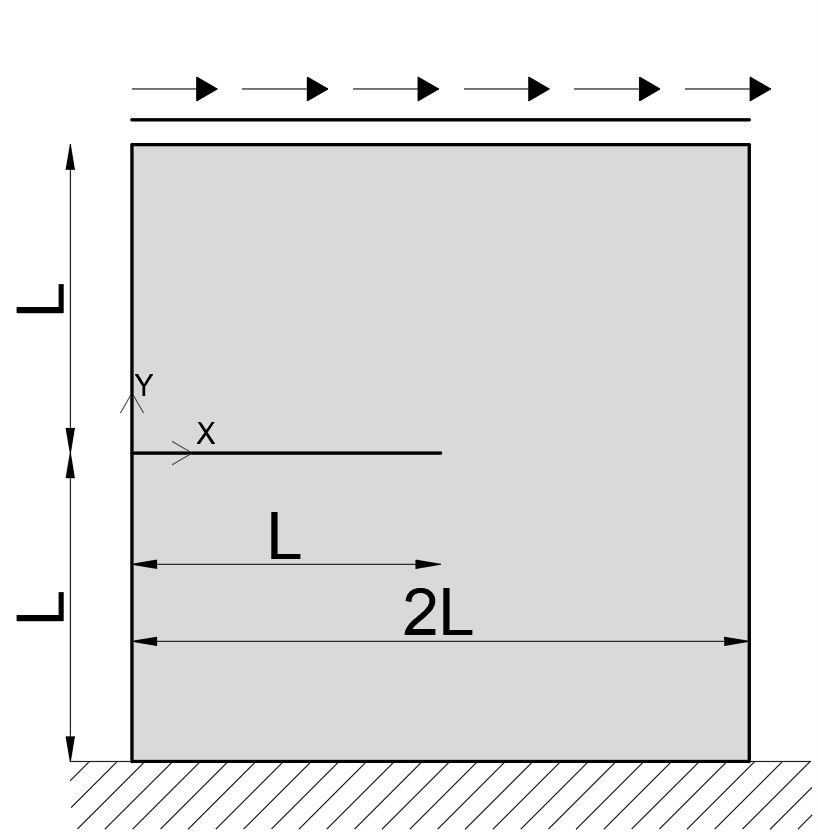}
			\caption{SEN shear test.}\label{fig:geometry_bc_d}
			\vspace{3cm}
		\end{subfigure}
	\end{center}
	\vspace{-80pt}
	\caption{Geometry and boundary conditions for the considered numerical tests.}
	\label{geometry and bc}
\end{figure}

After that, three benchmarks {well established in the phase-field literature for fracture} are considered, namely the DCB (\cite{BORDEN2014100,  GOSWAMI2020112808}), the SEN tensile test \cite{GERASIMOV2019990,GOSWAMI2020112808}, and the SEN shear test \cite{MieheWelschingerHofacker10}. To  make the comparison more homogeneous, %avoid affecting the accuracy results comparing these benchmarks, all test features were homogenized and 
we set all the above examples with the same conditions, avoiding differences in geometry, material, and imposition of the initial pre-field. Therefore, in all these tests (see Figures \ref{fig:geometry_bc_b} - \ref{fig:geometry_bc_d}) the uncracked domain is a square with a side of $1$ mm, a fracture toughness $G_c = 2.7\cdot 10^{-3}$ kN/mm, an internal length $\eps = 0.01$ mm, a Young's modulus $E = 210$ kN/mm$^2$, and a Poisson's ratio $\nu = 0.3$. In all these benchmarks,  we compare second- and fourth-order $AT_1$, based on a new definition of the mesh size in dependence of the $R_*$ value {(see §\ref{s.explicit})}, in terms of accuracy of the dissipated energy.
%as well as classical AT2 models in terms of accuracy  of the dissipated energy. 
More precisely, for each considered mesh, we compute the dissipated energy $\mathcal{D}_{\eps,h}$ at the end of the analysis and the corresponding (sharp) crack length $\ell$. Thus, the effective toughness $G_c^{\text{eff}}$ is evaluated as:
	\begin{equation}
		\label{eq:effective_toughness}
		G_c^{\text{eff}} =  \frac{\mathcal{D}_{\eps,h}}{\ell}\,.
	\end{equation}
	Then, given the input toughness $G_c$, the relative error \cite{BORDEN2014100} is computed as:
	\begin{equation}
		\label{relative error}
		\text{R.error} =  \left| \frac{G_c^{\text{eff}} - G_c}{G_c} \right|  = \left| \frac{\Delta G_c^{\text{eff}}}{G_c}  \right| \,.
\end{equation}
Note that, by $\Gamma$-convergence, $\mathcal{D}_{\eps,h} \to G_c \ell$, therefore $G_c^{\text{eff}} \to G_c$ and $\text{R.error} \to 0$ (as $h \ll \eps \ll 1$). 

All the simulations are performed with quadratic $C^1$-continuous B-Splines {(for a comprehensive discussion on the definition of these isogeometric functions, readers may refer to \cite{pieg1996nurbs,HUGHES20054135,CottrellCMAME2007} and references therein)} to directly compare all the models as in \cite{greco2024higher} independently on the order of the functional. 

%\NOTE{AP: Plotterei il dominio della figura 3(b) di dimensioni uguali alla 3(c) e 3(d), se possibile. Ad esempio mettendo una striscia bianca a layer zero per la figura 3(b) di dimensione uguale all'area tratteggiata delle figure 3(b) e 3(c) in modo tale da avere subfigures uguali. Dettagli minori: volendo, userei lo stesso trick per figura 3(d) perche' piu' piccina (e.g., striscia sopra le frecce orizzontali). In figura 3(c) al ``100'' (quello che leggiamo) la linea orizzontale che mimica i carrelli orizzontali si confonde con quella delle frecce verticali, quindi separerei almeno con offset pari a quanto utilizzato per le frecce orizzontali in figura 3(d).}
%\NOTE{mi sembra che qui ci starebbe bene dire quali sono le funzioni utilizzate, se ricordo bene, si usano polinomi dello stesso ordine indipendentemente dal grado del funzionale}

\subsection{Pure traction test\label{s.puretrac}} 

%\NOTE{AP: Cambiato altezza da $h$ a $t$ perche' uguale a mesh size altrimenti (aggiornare Figure 3(a). Tolto [] da unita' nel testo, manterrei solo nei grafici.}
A pure traction test is considered to measure the difference in terms of the elastic limit value for the second- and fourth-order $AT_1$ models with {respect} to a reference analytical solution. In order to simulate a pure 1D traction test, a specimen with \ap$t = 1$ mm \bl and $b = 20$ mm is analyzed. Imposed horizontal displacements, mimicking a traction condition, are applied at $x=\pm b/2$ (see Figure \ref{fig:geometry_bc_a}). 
%whereas vertical displacements are neglected by setting simple supported boundary conditions at $y=\pm \ap t \bl/2$. \NOTE{AP: Controllare->ci sono i carrelli in questo test? Perche' in figura 3(a) non ci sono.} 
{For this benchmark, we consider} a Young's modulus $E = 100$ kN/mm$^2$, a Poisson's ratio $\nu = 0$, an internal length $\eps= 0.125$ mm, and a material toughness $G_c= 0.01$ kN/mm. {Thus, the theoretical  value $\sigma_{\text{th}}$ of the elastic limit is computed as a function of the regularization parameter $c_\v$, the material toughness \textit{G}$_c$, and the shear modulus $\mu=E/(2(1+\nu))$, and reads:
\begin{equation}
	\label{eq:stress1D_analytical_solution}
	\sigma_{\text{th}} =  \sqrt{\frac{2G_c\mu}{c_\rho \eps}}\,.
\end{equation}
}
The quasi-static displacement-controlled loading history consists of $2000$ steps {with uniform increments, starting from} a minimum horizontal  (boundary) displacement $u_{\text{min}}= 0$ mm {up to} a maximum horizontal displacement $u_{\text{max}} = 20\cdot10^{-2}$ mm. 

%\NOTE{non mi torna $2000 \Delta u = u_{\text{max}}$\\ Risposta di Luigi: c'era un errore nel valore massimo dello spostamento.} 
%\NOTE{MN. Farei un double check dei risultati in Figura 4. Al quarto ordine la zona di transizione dovrebbe essere circa il doppio rispetto al secondo ordine, perch\`e si passa da $R_*=2$ a $R_* \sim 4$ ma qui mi sembra circa il triplo.}

\begin{figure}[!h]
	\begin{center}
		\begin{subfigure}[b]{1\textwidth}
			\includegraphics[width=14.1cm]{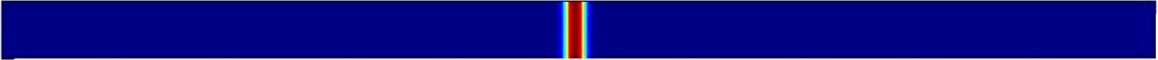}
			\caption{2$^{\text{nd}}$ order crack pattern.} 
		\end{subfigure} 			%	\hfill
		\begin{subfigure}[b]{1\textwidth}
			\includegraphics[width=0.9cm,,angle=-90]{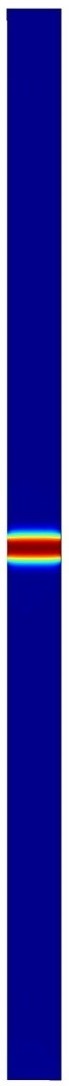}
			\caption{4$^{\text{th}}$ order crack pattern.}
		\end{subfigure}
	\end{center}
	\vspace{-10pt}
	\caption{Last step crack patterns for the $AT_1$ models.}
	\label{crack patter pure traction}
\end{figure}

All analyses are performed considering a uniform mesh with size $h = 0.0625$ mm and,
%size h = 0.0625 [mm] in y-direction for both models, instead in the x-direction h = 0.0625 [mm] for the second-order AT1 and h = $0.0625\cdot$R [mm] for the fourth-order AT1. 
in Figure \ref{crack patter pure traction}, we report the crack patterns corresponding to the last loading steps for the second- and fourth-order $AT_1$ models\footnote{{Due to the homogenous response in terms of stress in the body for this test, a crack could nucleate anywhere in the domain. Therefore, a small ($\sim10^{-4}$)
		%\NOTE{Inserirei l'order of magnitude della perturbazione $\sim10^{-xx}$ qui nella nota.}) 
perturbation is applied at the center of the domain to allow the crack to nucleate in the middle of the specimen.}}.

The elastic limits, obtained by computing the {numerical} stress at the step before crack onset (denoted by $\sigma_{c}$) are compared in Table \ref{Tab: Accuracy pure traction } to $\sigma_{\text{th}}$ (see Equation \eqref{eq:stress1D_analytical_solution}){, showcasing that the proposed $AT_1$ high-order model is able to correctly capture the stress peak, providing a relative error in line with the second-order $AT_1$ model.}
	
\begin{table}[!h] \centering
	\caption{Pure traction test: relative percentage error.}\label{Tab: Accuracy pure traction } 
	\vspace{-5pt}
	\renewcommand\arraystretch{1.5}
	\renewcommand\tabcolsep{4pt}
	\begin{tabular}{c c c c}
		\hline
		\centering
		{$AT_1$-model} & {$\sigma_{\text{th}}$} & {$\sigma_{c}$} & {R.error} \\
				\hline
		{[-]} & {[kN/mm$^2$]} & {[kN/mm$^2$]} & {[\%]}\\
		\hline
		$2^{\text{nd}}$ order & 1.73205 & 1.73200 & 0.00293\\
		$4^{\text{th}}$ order & 1.34103 & 1.34101 & 0.00206\\
		\hline
	\end{tabular}
\end{table}

\subsection{Mesh-size considerations and {its} appropriate selection}\label{sec: mesh size study}

As far as the choice of the mesh size, we argue in the following way. From the theoretical point of view (see Theorem \ref{t.gammah}), the $\Gamma$-limit is recovered for $h \ll \eps$, i.e.,~in the limit when the mesh size is much smaller then the internal length. In practice, this condition means that the mesh size should be small enough to ``resolve'' the phase-field profile in its transition region, which is of order $R_* \eps$. This observation entails an important consideration in the choice of the mesh size, crucial for the phase-field approach. From Proposition \ref{p.eq} and \ref{p.eq-app}, we know that $R_* \sim 3.8$ for the fourth order and $R_*=2$ for the second order. Therefore, we choose 
\begin{equation}
	h = \frac{R_*\,\eps}{n} ,\, \text{ for $n = 4, 8,  16$, }
	\label{mesh definition}
\end{equation}
as reference mesh sizes for the both the $AT_1$ functionals. Note that, in this way, we provide also a set of results for roughly the same mesh (see for instance Figure \ref{DCB accuracy vari modelli}) since $h = 0.005$ mm corresponds to $R_* \, \eps / 4$ for the second order and to $R_* \eps / 8$ for the fourth order.  Actually, to better check the error trend of our fourth-order functional, we provide also the numerical results for $h =  R_* \eps/ 2$, clearly with $R_* = 3.83$.  For $AT_2$ the above argument, based on $R_*$, does not apply, since the theoretical results would give $R_* = + \infty$, therefore, for a direct comparison we choose the mesh sizes corresponding to the second-order $AT_1$. 

Finally, in Figure \ref{fig:c}  it can be seen how the 1D phase-field profile is influenced by progressively decreasing the coefficient $\rho$ in front of the high-order term of the surface energy. %(namely $c={1,1/2...}$) within the interval $I = \pm R$ where the profile is positive defined. 
As $\rho$ decreases, the fourth-order optimal profile tends to the second-order optimal one, providing a ``less'' regular function with smaller crack width. As a consequence, in the parametric study of \S \ref{sec: sensitivity study} we compare the numerical results obtained both  by fixing the mesh size and linking it to $\rho$, through the value $R_*(\rho)$. 

\subsection{DCB test}

The DCB test has been employed in \cite{BORDEN2014100, GOSWAMI2020112808, greco2024higher} to compare low- and high-order $AT_2$ formulations. In our case, we consider a squared sample of side {$2L=1$ mm} with an initial phase-field crack located at $y = 0$ mm and {$x \in [0, \,L]$ mm}, modeled by the {interpolated phase-field variable (IPF)} technique  \cite{greco2024higher}. A quasi-static displacement-controlled loading history is applied {at the left vertical edge (see Figure \ref{fig:geometry_bc_b})} consisting of 37 loading steps, with a minimum vertical boundary displacement $u_{\text{min}} = 0$ mm, a uniform increment $\Delta u = 5\cdot10^{-4}$ mm, and a maximum vertical boundary displacement $u_{\text{max}} = 18\cdot10^{-3}$ mm, {tuned} to completely break the specimen. {Additionally, we completely restrain the horizontal displacement at $x=0$ mm.}

%\NOTE{non mi torna $37 \Delta u = u_{\text{max}}$ \\ Risposta di Luigi: errore di scrittura dei valori: avevo invertito la storia di carico del pure traction con quella del DCB.} 
{We provide an extensive comparison of the proposed fourth-order $AT_1$ model with respect to the second-order one, as well as second- and fourth-order $AT_2$ models, such that for the second-order $AT_1$ the considered meshes are $h = R_* \eps/4$, $h = R_* \eps/8$. Instead, for the fourth-order $AT_1$ we consider $h = R_* \eps/4$, $h = R_* \eps/8$, and $h = R_* \eps/16$, whereas for the classical second- and fourth-order $AT_2$ formulations we examine $h = \eps/4$, $h = \eps/8$. All $AT_2$ and the $AT_1$ second-order tests are performed with locally refined mesh, obtained using a non-uniform knot-vector, while all fourth-order $AT_1$ tests are performed with uniform meshes. Despite this difference in the mesh definition, the $AT_1$ fourth-order model showcases computational advantages even performing analyses with uniform mesh (see \S \ref{sec: degree saving} for an in-depth discussion in terms of trade-off between accuracy and computational cost for all considered models on several benchmarks).}
%  All tests are performed with uniform mesh and, due to the ``larger'' width profile given by the IV-order model, the number of elements saved for each test is above 50\%.
\begin{figure}[h!!]
	\begin{center}
		\includegraphics[width=10cm]{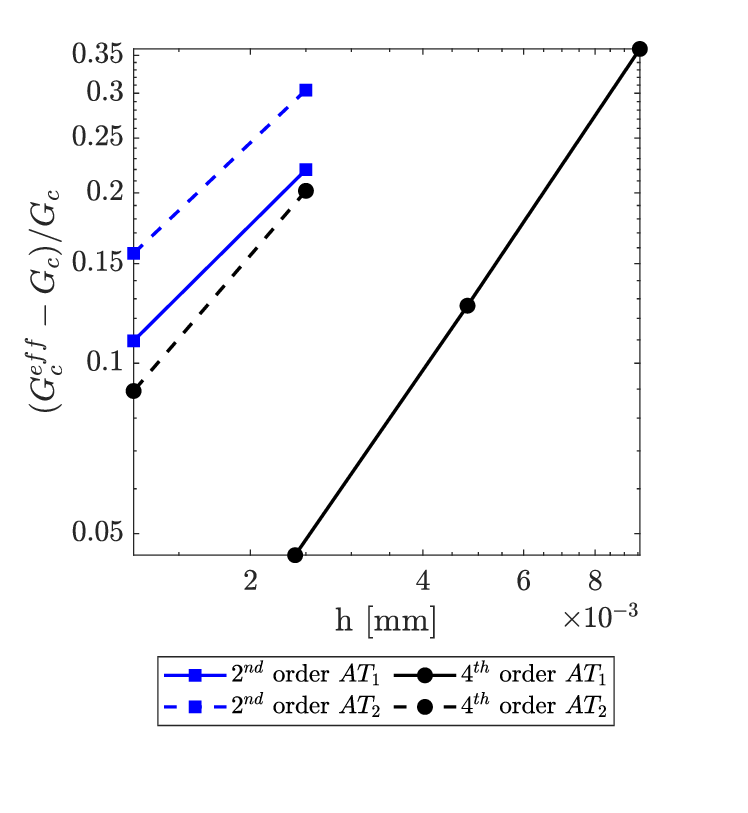}
	\end{center}
	\vspace{-20pt}
	\caption{DCB test: accuracy analysis considering all $AT$-models represented in double logarithmic scale. }
	\label{DCB accuracy vari modelli}
%	\NOTE{Per quanto ne so, se si mette il $\log$ nella label (in ascissa o in ordinata) poi i valori sugli assi devono essere i $\log$ delle variabili, in questo caso diventerebbero entrambi dei valori negativi. In questo grafico (e in quelli seguenti)  invece mi sembra che i valori sugli assi non siano i $\log$ ma i valori ``reali''. Insomma, i valori e la label devono essere consistenti, o entrambi $\log$ o entrambi non-$\log$. Forse suggerirei di non mettere il $\log$ nella label e lasciare i valori, scrivendo nel caption che il grafico e' in scale logaritmica.}
\end{figure}
%\NOTE{AP: Luigi, cambieresti la notazione della legenda nella Figure 5 per uniformare (->$AT_1$, $AT_2$)? Grazie!}
%\NOTE{AP: Figure 5 e' non log, ma i trattini piu' fitti nello spacing mi confondono, controlleresti? Penso che eventualmente questo commento di applichi a tutte le figure inerenti.}

{In Table \ref{Tab: errors DCB AT}, it can be observed that for the fourth-order $AT_1$ we have a smaller error values, allowing for a significant reduction in terms of degrees of freedom, and a higher convergence rate, in comparison with the other $AT_1$ model, computed as:
\begin{equation}
	\text{CR} = \frac{\log{\left((\Delta G_c^{\text{eff}})_{k+1}\right)} - \log{\left((\Delta G_c^{\text{eff}})_{k}\right)} }{ \log\left(h_{k+1}\right) - \log\left(h_{k}\right)}\,,
	\label{CR}
\end{equation}
where $k$ represent the loading steps, $\Delta G_c^{\text{eff}}$ is the toughness variation between two steps and $h$ is the mesh size.	
This result also holds in the comparison with the $AT_2$ models as we can observe in Figure \ref{DCB accuracy vari modelli}. Furthermore, we remark that, for this test, all $AT_1$ models are more accurate than the $AT_2$ ones and the proposed $AT_1$ fourth-order functional is able to provide more accurate results than all investigated classical literature models. More specifically, it can be seen observed for comparable meshes (namely in the case of $h=R_*\eps/16$) that the fourth-order $AT_1$ model shows a relative error on the toughness about ten times lower than the other models.}
%\NOTE{AP: Ho unificato perche' c'erano ripetizioni, e' corretta l'info sulla mesh sopra?}

%\begin{table}[!h]  \centering
		%\caption{DCB test: relative percent error \\ on the toughness for mesh size $h = %0.0025$ [mm]}\label{Tab: errors DCB AT} 
		%\renewcommand\arraystretch{1.5}
		%\renewcommand\tabcolsep{8pt}
		%\begin{tabular}{ c c c c c }
			%\hline {Quantities} & \multicolumn{2}{c}{{2$^{\text{nd}}$ order}} & %\multicolumn{2}{c}{4$^{\text{th}}$ order}  \\[-4pt]
%			{-} & {AT2} & {AT1} & {AT2} & {AT1} \\
%			\hline 
%			R.error [\%]  & 30.96 &  21.99 & 20.17 & \bb $\sim 4.5$ \bl  \\ 
%			\hline
%			CR [-] & 0.959 & 1.005 & 1.175 & 1.507 \\
%			\hline
			%  \botrule 
%	\end{tabular} 
%\end{table}
%
\begin{table}[!htbp]  \centering
	\caption{DCB test: $AT_1$ model relative percent error on the toughness.}\label{Tab: errors DCB AT} 
	\renewcommand\arraystretch{1.5}
	\renewcommand\tabcolsep{15pt}
	\begin{tabular}{ c c c }
		\hline
		 {\textit{h}} & {{2$^{\text{nd}}$ order}} & {4$^{\text{th}}$ order}  \\[-4pt]
		{[mm]} & {[\%]} & {[\%]} \\
		\hline 
		$R_*\eps/8$  & 21.99 &  12.63  \\ 
		$R_*\eps/16$ & 10.95 & 4.50     \\ 
		\hline
		CR [-] & 1.01 & 1.51 \\
		\hline
		%  \botrule 
	\end{tabular} 
\end{table}
%\NOTE{AP: Cifre significative in Table 4 (es. 4.5). Dovrebbero essere 4 ovunque.}
%
%\begin{table}[!h] \centering
%	\caption{DCB test: relative percent error on the toughness for different AT-models. \bb Bold: results for comparable mesh size. \bl}\label{Tab: errors DCB AT} 
%	\vspace{-5pt}
%	\renewcommand\arraystretch{1.5}
%	\renewcommand\tabcolsep{8pt}
%	\begin{tabular}{ c c c c }
%	\hline {$h$} & \multicolumn{2}{c}{{2$^{\text{nd}}$ order}} & {4$^{\text{th}}$ order} \\
%	\hline {[mm]} & {AT1} & {AT2} & {AT2} \\
%		\hline 
%		\boldmath{$\eps/4 \sim 2.1 \cdot 10^{-3}$} & {\bf 21.99} & {\bf 30.96} &  {\bf 20.17}  \\ 
%		$\eps/8$ & 10.95 & 15.63 & 8.93  \\ 
%		\hline
%		$ CR $  & 1.005 & 0.959 & 1.175 \\
%		\hline
%		%  \botrule 
%	\end{tabular}
%	\hspace{2cm}
%	\begin{tabular}{ c c }
%	\hline
%	  {$h$} & {4$^{\text{th}}$ order} \\
%	\hline 
%	 {[mm]} & {AT1}\\
%		\hline 
%		$R_*\eps/4 $ & 35.93   \\ 
%		 $R_* \eps/8 $ & 12.63   \\ 
%		  \boldmath{$2.1 \cdot 10^{-3}$} & \boldmath{$\sim 4$}   \\ 
%		\hline
%		 $ CR $ & 1.507 \\
%		\hline
%		%  \botrule 
%	\end{tabular}
%	\qquad
%\end{table}

\begin{figure}[!htbp] 
	\begin{center}
		\begin{subfigure}[b]{0.3\textwidth}
			\includegraphics[width=4cm]{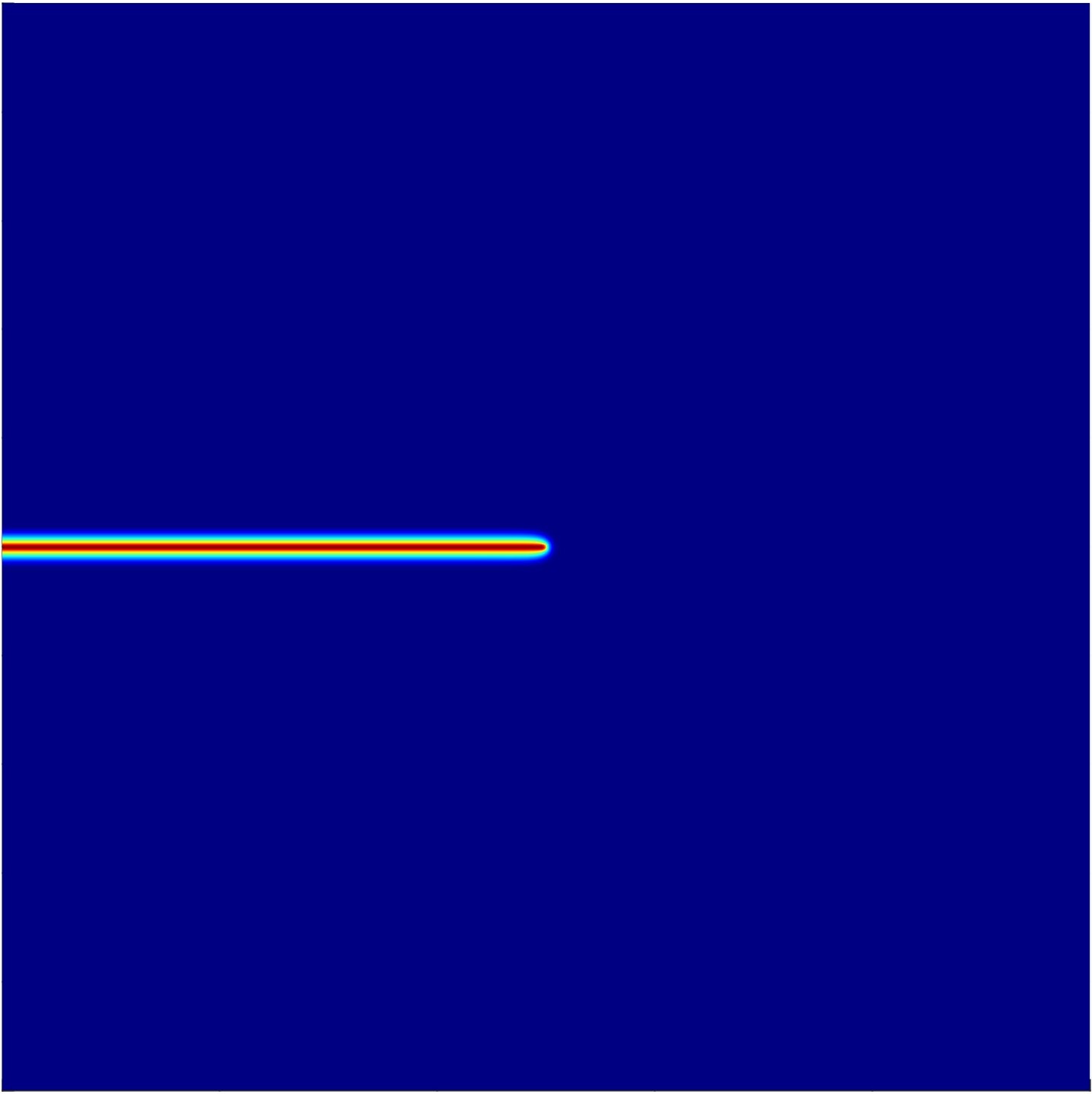}
			\caption{Step 1:  \\initial condition.}
		\end{subfigure}
		%	\hfill
		\begin{subfigure}[b]{0.3\textwidth}
			\includegraphics[width=4cm]{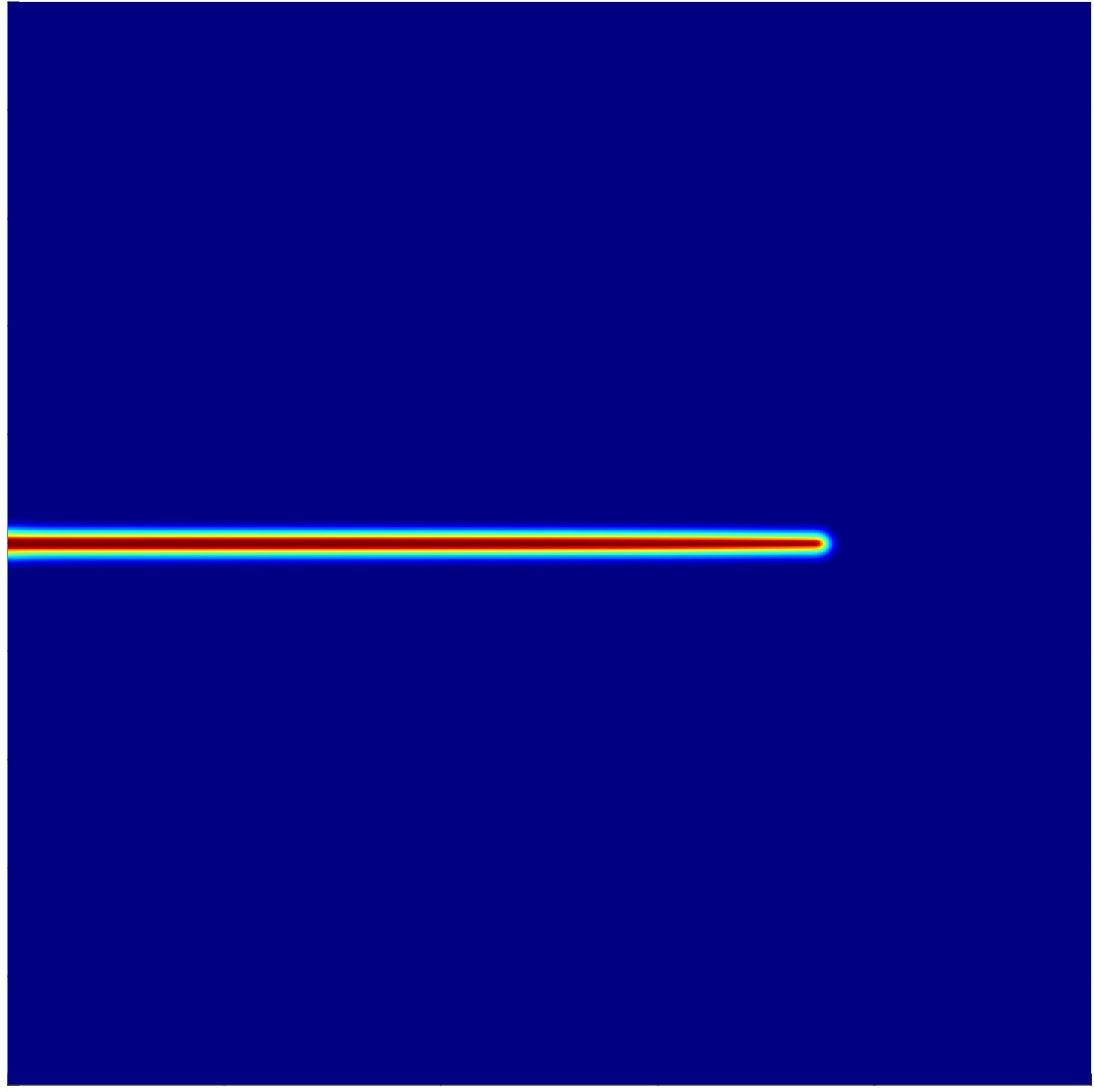}
			\caption{Step 18: fracture\\propagation in the domain. }
		\end{subfigure}
		%   \hfill
		\begin{subfigure}[b]{0.3\textwidth}
			\includegraphics[width=4.55cm]{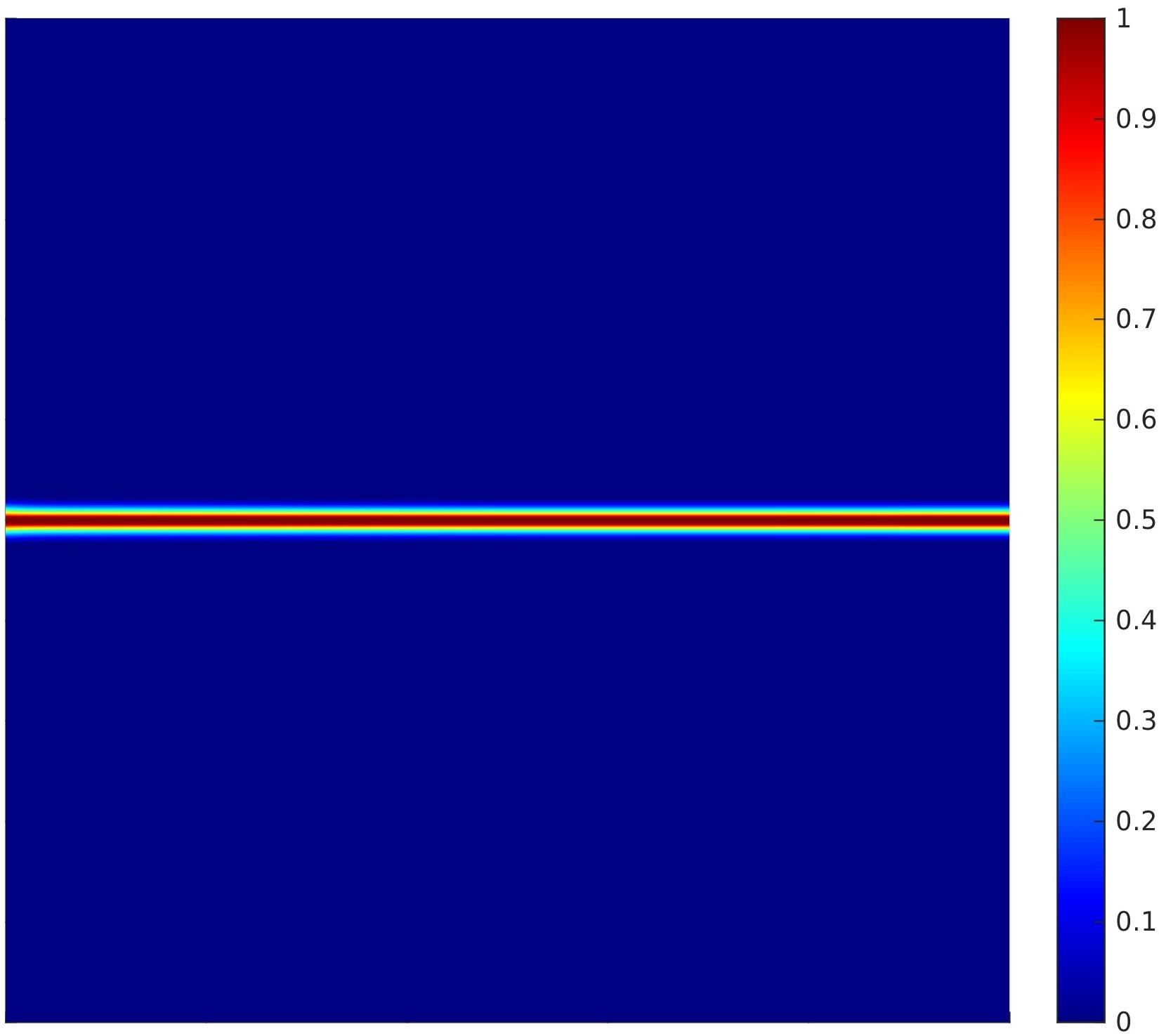}
			\caption{Step 37: total breakage \\of the specimen.}
		\end{subfigure}
		\begin{subfigure}[b]{0.3\textwidth}
			\includegraphics[width=4cm]{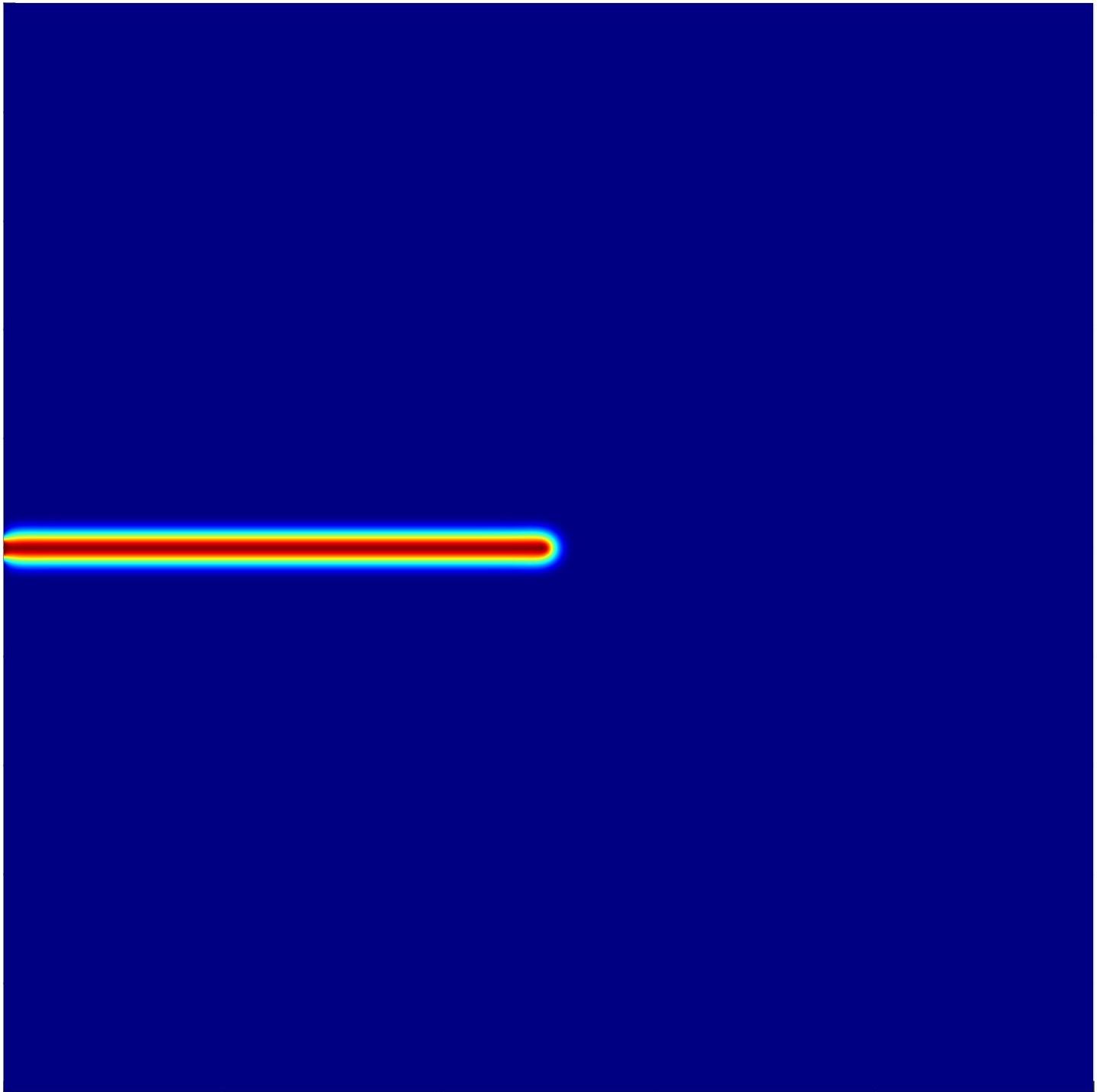}
			\caption{Step 1:  \\initial condition.}
		\end{subfigure}
		%	\hfill
		\begin{subfigure}[b]{0.3\textwidth}
			\includegraphics[width=4cm]{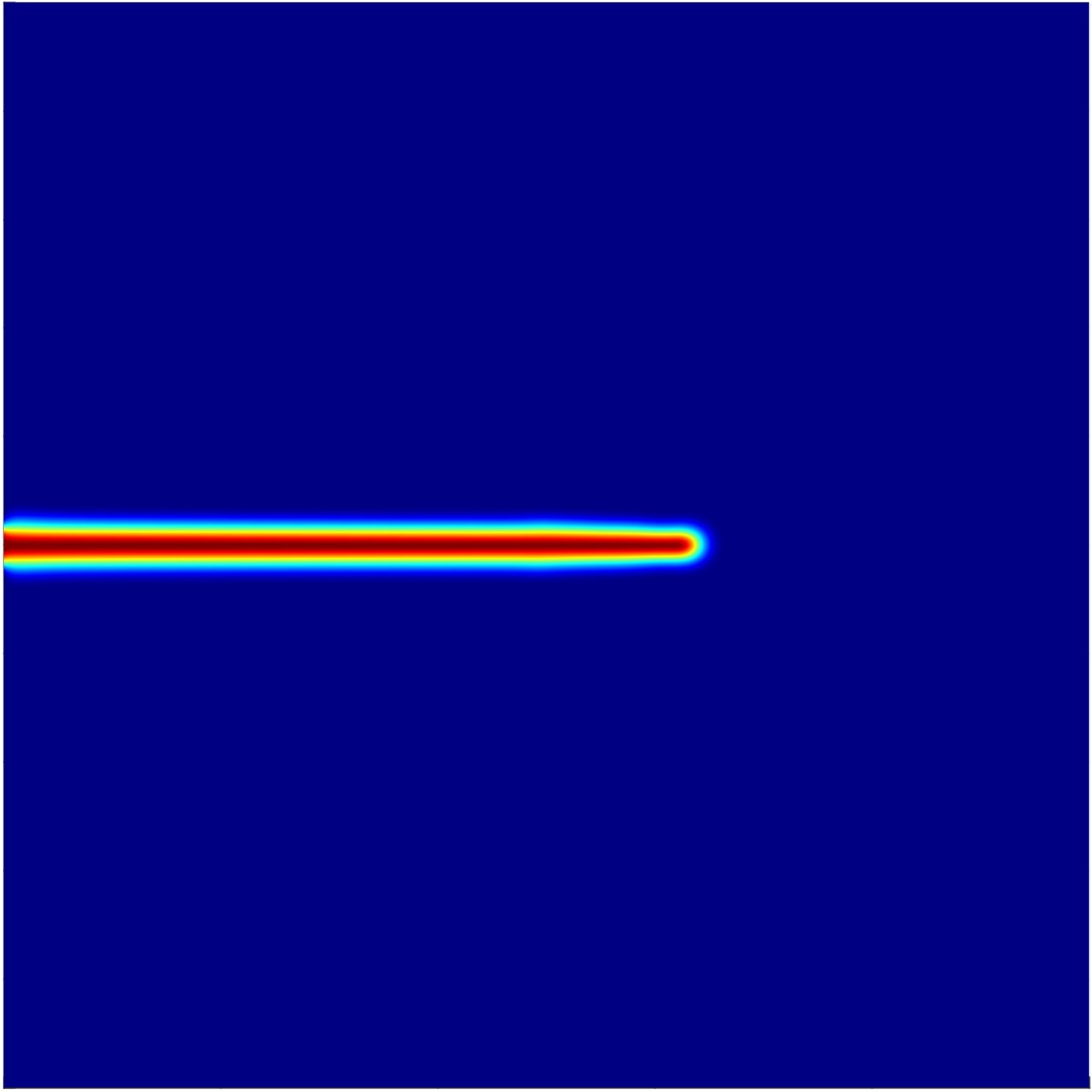}
			\caption{Step 18: fracture \\propagation in the domain. }
		\end{subfigure}
		%   \hfill
		\begin{subfigure}[b]{0.3\textwidth}
			\includegraphics[width=4.6cm]{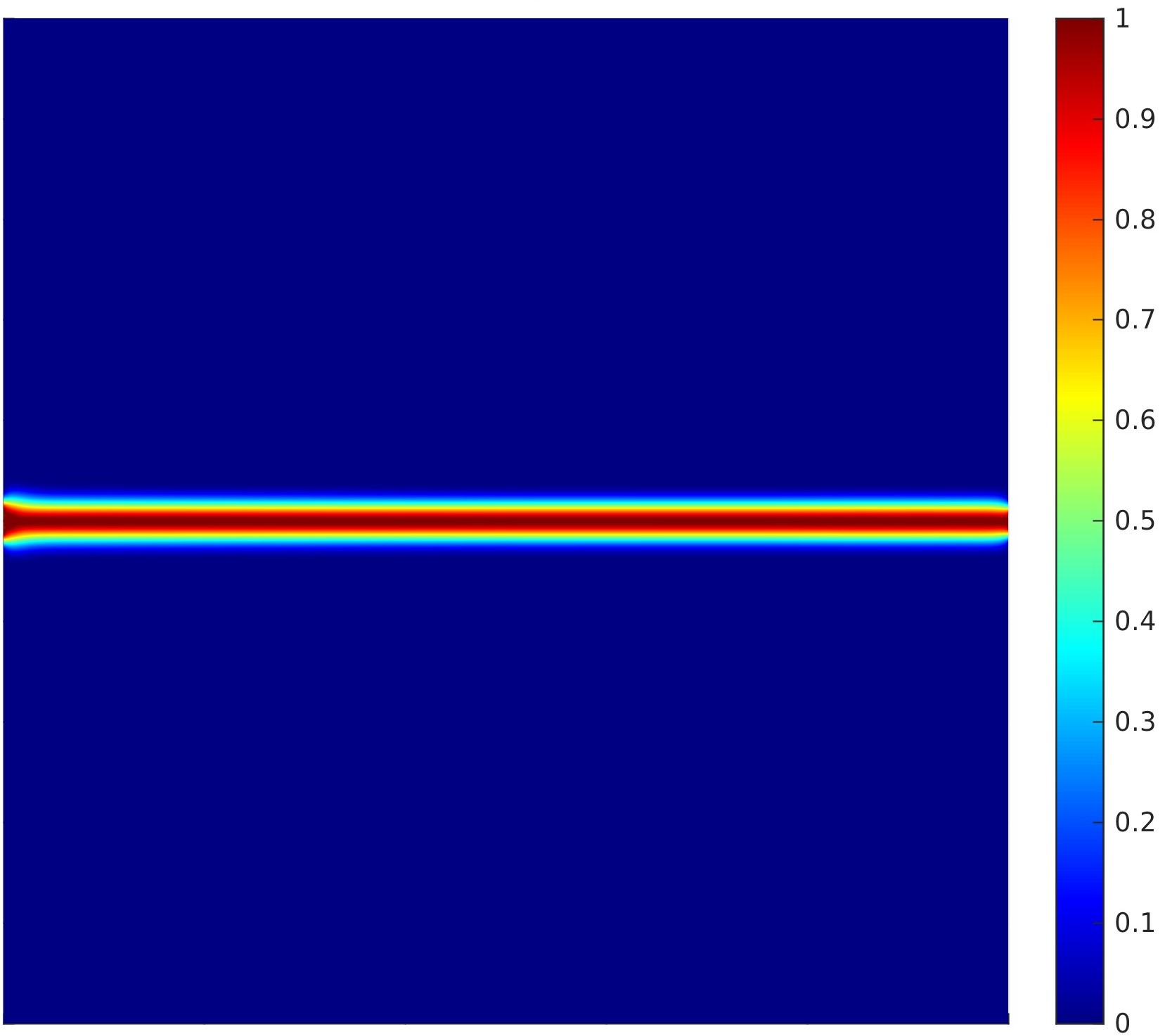}
			\caption{Step 37: total breakage \\of the specimen.}
		\end{subfigure}
	\end{center}
	\vspace{-15pt}
	\caption{DCB test: crack patterns of three different loading steps for the $AT_1$ second-order (top row) and $AT_1$ fourth-order model (bottom row) considering $ h= 0.0025 $ mm.}
	\label{crack patter DCB novel AT}
\end{figure}

In Figure \ref{crack patter DCB novel AT}, we compare the phase-field contour plots for the $AT_1$ second-order and the novel $AT_1$ fourth-order formulations. It can be observed that the proposed model provides a thicker profile than the one obtained for the second-order formulation due to the higher regularity of the model. Nevertheless, the novel model allows for comparable values in terms of relative error on the toughness (see Table \ref{Tab: errors DCB AT}) showing a significant reduction of the relative error. 

 \subsection{SEN tensile test}
 \label{senT}
We consider now the SEN tensile test \cite{MieheWelschingerHofacker10}, another mode I fracture test leading to complete breakage of the specimen. Geometry and boundary
conditions of this benchmark are summarized in Figure \refeq{fig:geometry_bc_c}. Also in this case, the initial pre-crack is located at $y = 0$ mm and {$x \in [0, \,L]=[0, \,0.5]$ mm} and is modeled with the technique of the IPF. A quasi-static displacement-controlled loading history is applied {at the top horizontal edge}, comprising 21 loading steps, a minimum boundary displacement $u_{min} = 0$ mm, a uniform increment $\Delta u = 3\cdot10^{-4}$ mm, and a maximum boundary displacement $u_{max} = 6\cdot10^{-3}$ mm. {Instead, we constrain the horizontal displacement at $y=L$. For this benchmark, we consider mesh sizes $h = R_* \eps/4$ and $h = R_* \eps/8$ for the $AT_1$ models, whereas for the $AT_2$ functionals we examine $h = \eps/2$ and $h = \eps/4$.}

\begin{figure}[h!!]
	\begin{center}
		\includegraphics[width=12cm]{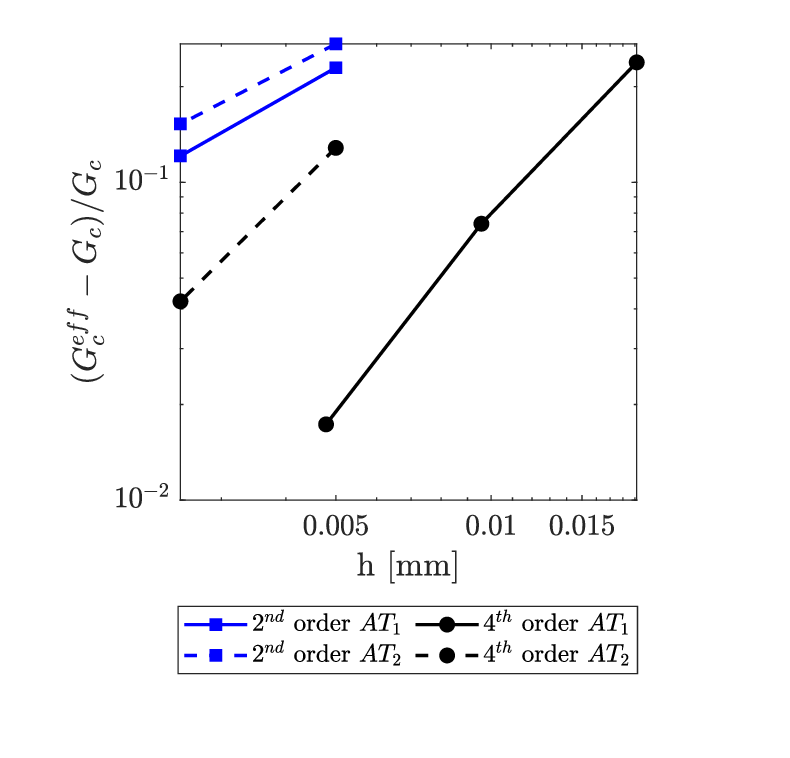}
	\end{center}
	\vspace{-20pt}
	\caption{SEN tensile test: accuracy analysis considering all the $AT$-models.}
	\label{tensile accuracy vari modelli}
\end{figure}
%\NOTE{in figura 7 e 9 metterei nella legenda prima riga AT1 e seconda riga AT2}
%\NOTE{AP: Luigi, mi cambieresti la notazione della legenda nella Figure 7 per uniformare (->$AT_1$, $AT_2$)? Thx.}

\begin{figure}[h!!]
	\begin{center}
		\begin{subfigure}[b]{0.3\textwidth}
			\includegraphics[width=4cm]{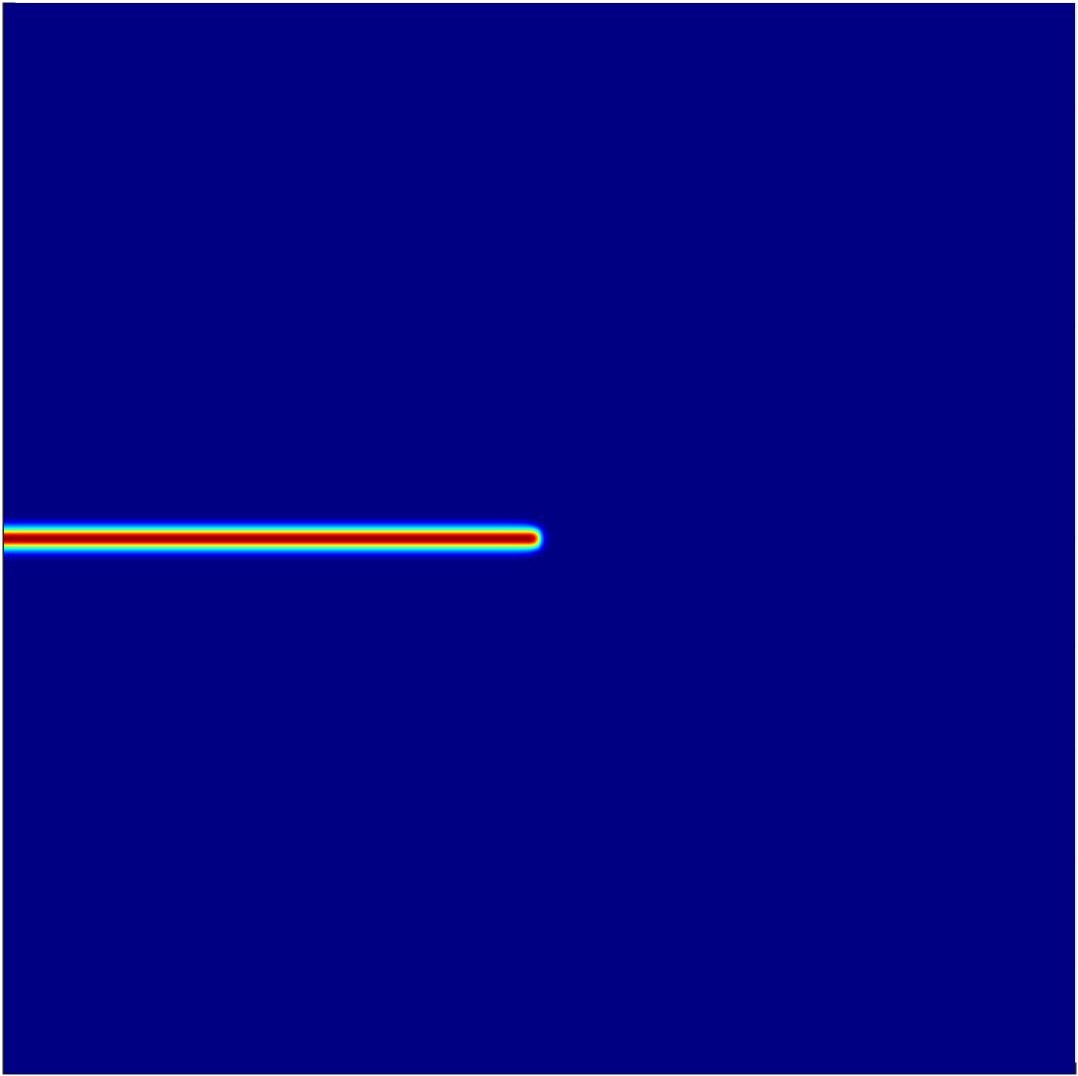}
			\caption{Step 1:  \\initial condition.}
		\end{subfigure}
		%	\hfill
		\begin{subfigure}[b]{0.3\textwidth}
			\includegraphics[width=4cm]{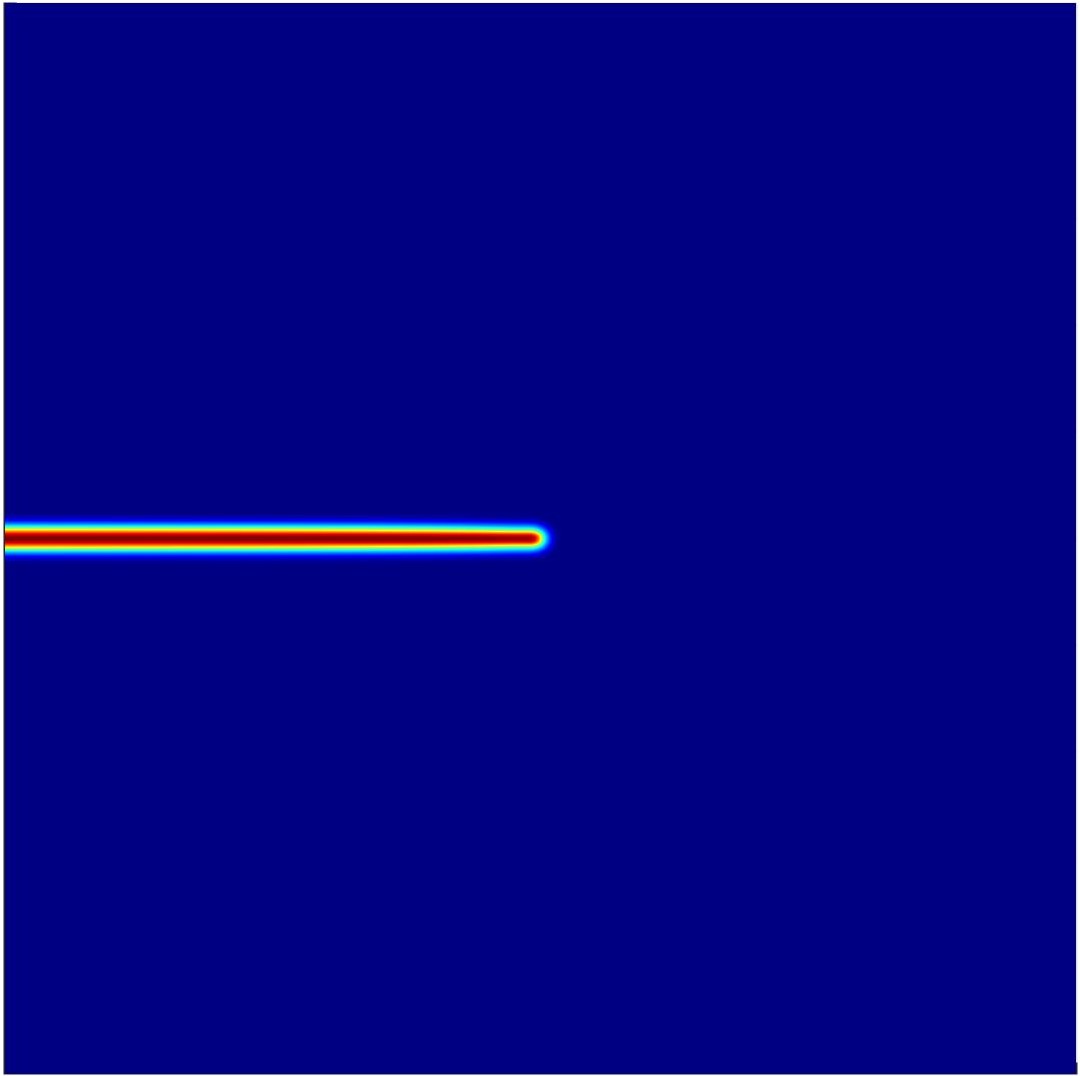}
			\caption{Step 17: fracture \\propagation in the domain. }
		\end{subfigure}
		%   \hfill
		\begin{subfigure}[b]{0.3\textwidth}
			\includegraphics[width=4.7cm]{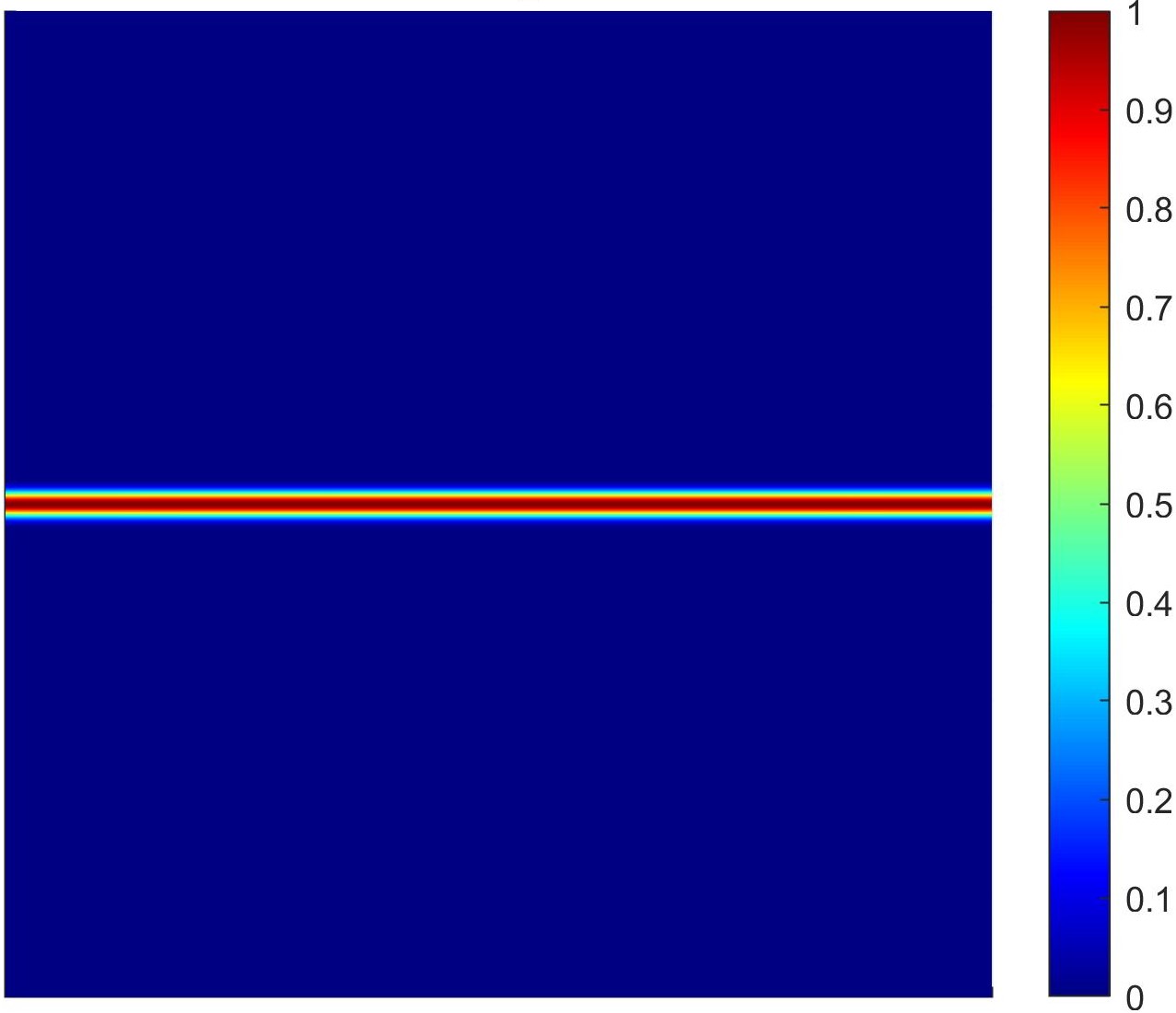}
			\caption{Step 21: total breakage \\of the specimen.}
		\end{subfigure}
		\begin{subfigure}[b]{0.3\textwidth}
			\includegraphics[width=4cm]{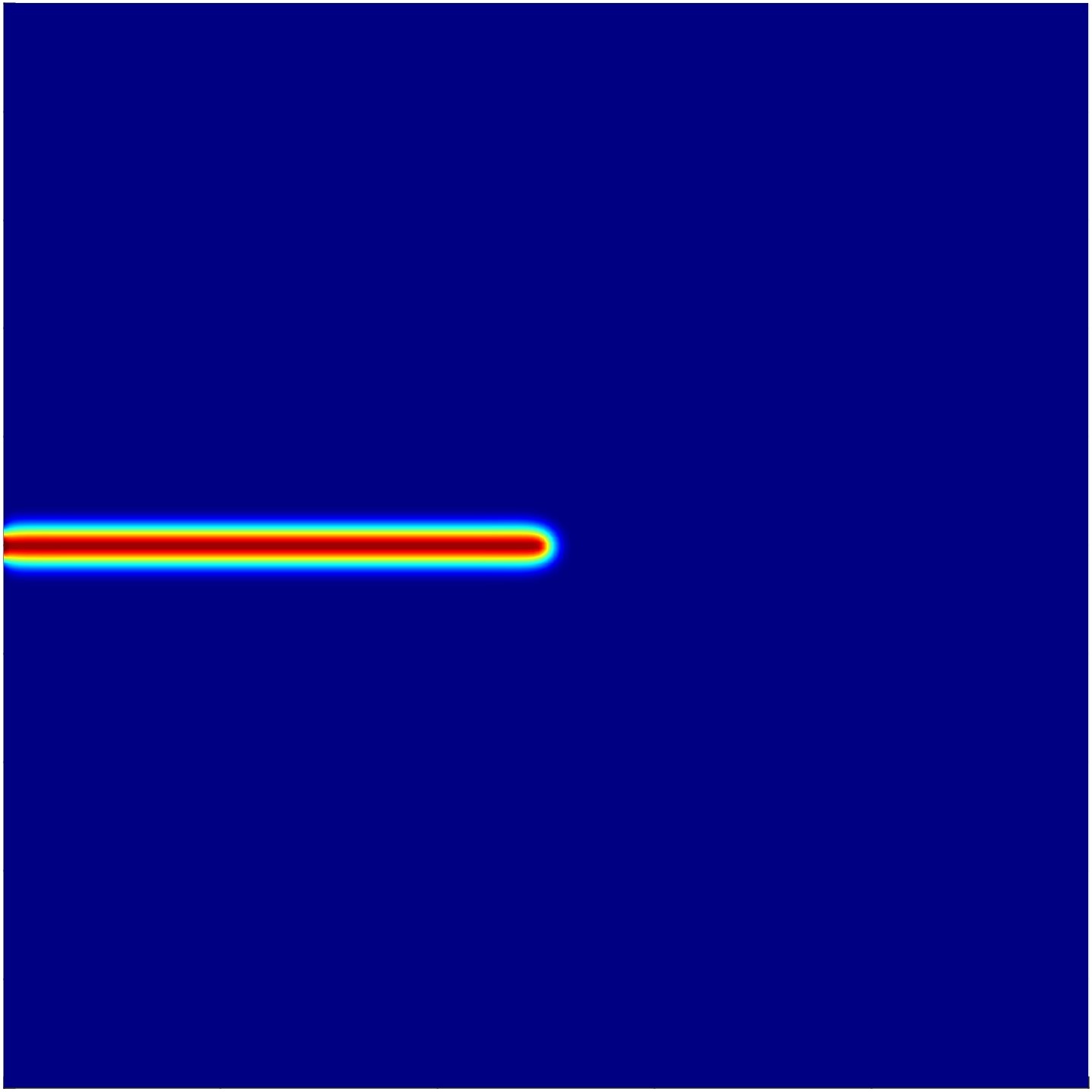}
			\caption{Step 1:  \\initial condition.}
		\end{subfigure}
		%	\hfill
		\begin{subfigure}[b]{0.3\textwidth}
			\includegraphics[width=4cm]{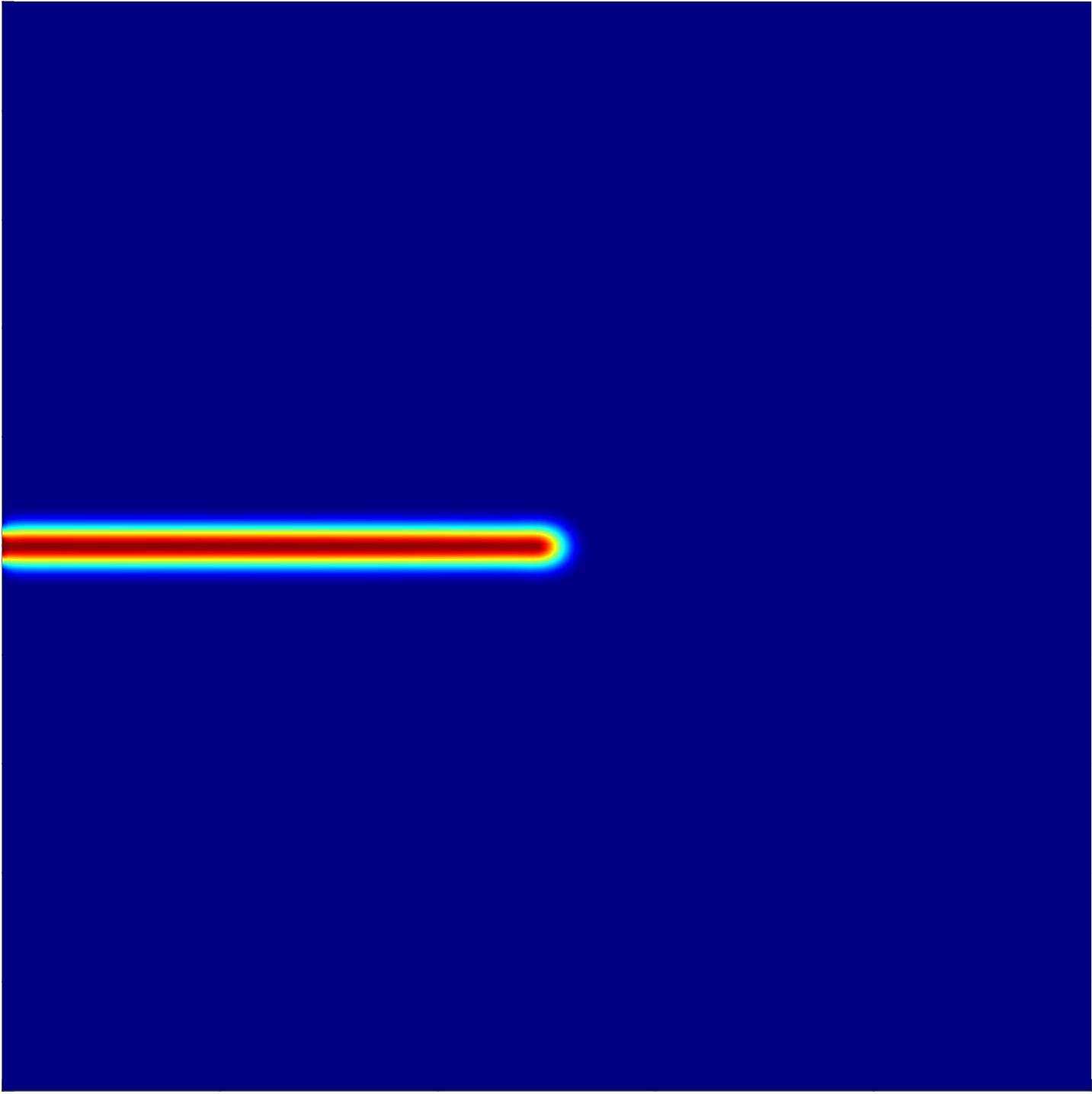}
			\caption{Step 16: fracture \\propagation in the domain. }
		\end{subfigure}
		%   \hfill
		\begin{subfigure}[b]{0.3\textwidth}
			\includegraphics[width=4.6cm]{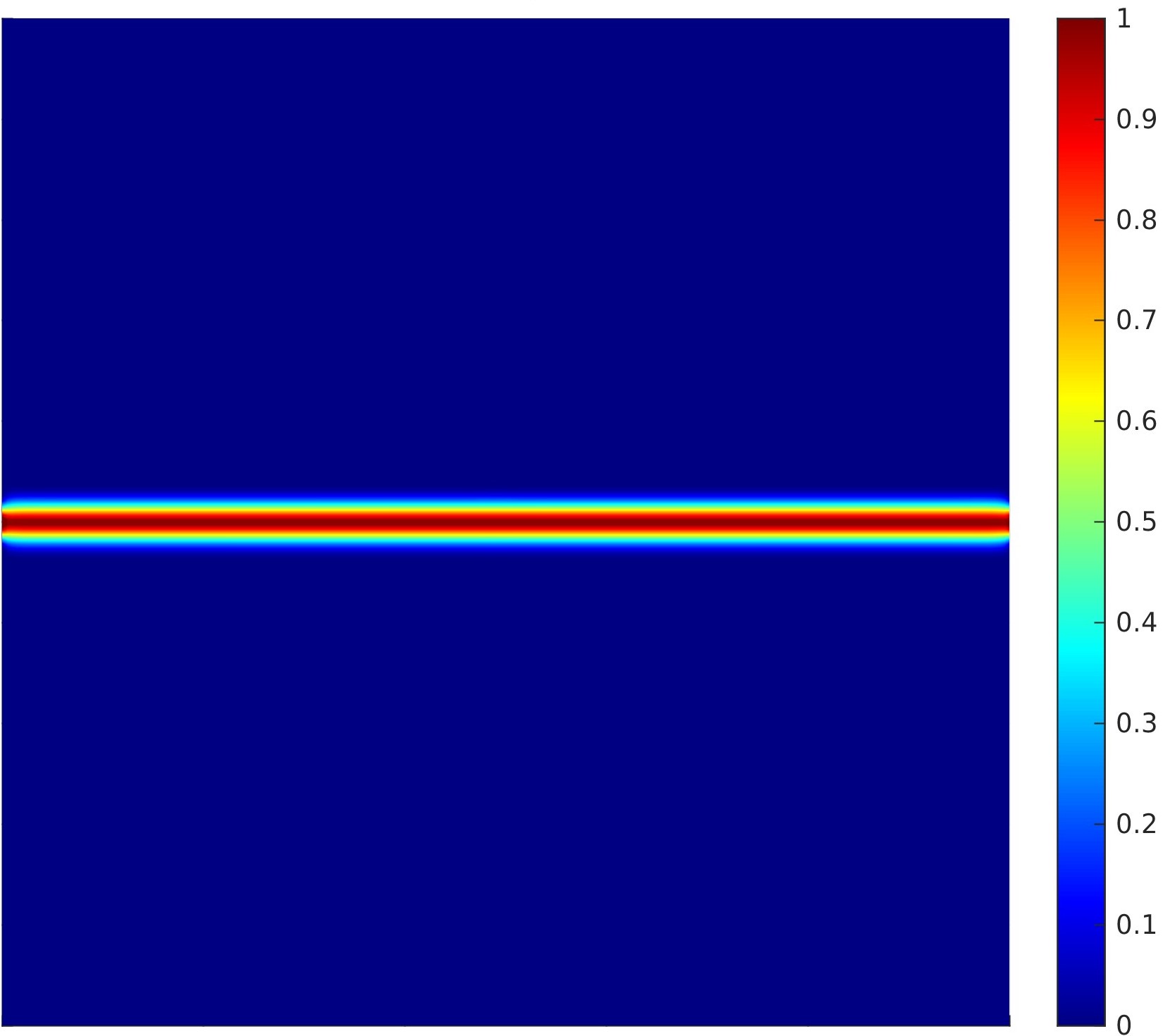}
			\caption{Step 17: total breakage \\of the specimen.}
		\end{subfigure}
	\end{center}
	\vspace{-15pt}
	\caption{SEN tensile test: crack patterns of three different loading steps for the $AT_1$ second-order (top row) and $AT_1$ fourth-order model (bottom row) considering $ h= 0.005 $ mm.} 
	\label{crack patter Tensile novel AT}
\end{figure}

{Let us comment on the accuracy trend reported in Figure \ref{tensile accuracy vari modelli}: it can be observed that the proposed $AT_1$ model has a slope comparable with the $AT_2$ fourth-order model, however, for fixed mesh size, it is more accurate than the analyzed second-order models; resulting in practice in a significant reduction of the computational cost (see §\ref{sec: degree saving} for more details). Namely, it can be seen in Figure \ref{tensile accuracy vari modelli} that, for the fourth-order $AT_1$ model, if we  refine the mesh from $h = R_* \eps/4$  to $h = R_* \eps/8$, the relative error on the toughness decreases from 7.41\% to 1.73\%. Also for this example, for comparable meshes, the $AT_1$ fourth-order model shows an error approximately ten times smaller than the error obtained with the other models and we observe that the $AT_1$ models are more accurate than the $AT_2$ functional counterparts. As expected, also for the SEN tensile test, it can be seen in Figure \ref{crack patter Tensile novel AT} that the proposed $AT_1$ model provides a smoother crack pattern than the $AT_1$ second-order model, but entails much lower values of relative error on the toughness for a fixed mesh size in all cases (see Table \ref{Tab: errors tensile AT})}. 

\begin{table}[!h]  \centering
	\caption{SEN tensile test: $AT_1$ model relative percent error on the toughness represented in double logarithmic scale.}\label{Tab: errors tensile AT}   
	\renewcommand\arraystretch{1.5}
	\renewcommand\tabcolsep{15pt}
	\begin{tabular}{ c c c }
		\hline
		{\textit{h}} & {{2$^{\text{nd}}$ order}} & {4$^{\text{th}}$ order}  \\[-4pt]
		{[mm]} & {[\%]} & {[\%]} \\
		\hline 
		$R_*\eps/4$  & 22.94 &  7.41  \\ 
		$R_*\eps/8$ & 12.11  & 1.70     \\ 
		\hline
		CR [-] & 0.92 & 1.69 \\
		\hline
		%  \botrule 
	\end{tabular} 
\end{table}
%\NOTE{AP: Cifre significative in Table 4, intendo e.g., 1.7. Dovrebbero essere 4 ovunque.}

%\begin{table}[!h]  \centering
%	\caption{SEN tensile: relative percent error \\ on the toughness for mesh size $h = 0.005$ [mm]}\label{Tab: errors tensile AT}  
%	\renewcommand\arraystretch{1.5}
%	\renewcommand\tabcolsep{8pt}
%	\begin{tabular}{ c c c c c }
%		\hline {Quantities} & \multicolumn{2}{c}{{2$^{\text{nd}}$ order}} & \multicolumn{2}{c}{4$^{\text{th}}$ order}  \\[-4pt]
%		{-} & {AT2} & {AT1} & {AT2} & {AT1} \\
%		\hline 
%		R.error [\%]  & 27.26 &  22.94 & 12.84 & \mg $\sim 1.7 $ \bl  \\ 
%		\hline
%		CR [-] & 0.837 & 0.922 & 1.604 & 1.686 \\
%		\hline
		%  \botrule 
%	\end{tabular} 
%\end{table}
\bl

%\NOTE{AP: Forse ha senso mettere lo step prima della rottura in Figure. 8(e) per dire che il IV anticipa e in Figure 8(f) lo step attuale 17, che ne pensate? E mantenere le captions come per il II At1.}
%\NOTE{questo dato con $R_* \eps/8$ andrebbe inserito in tabella}

\subsection{SEN shear test}
\label{senS}
{To provide a more challenging
crack pattern evolution,} we consider the SEN shear test having material properties and geometry as in \cite{GERASIMOV2019990, GOSWAMI2020112808, greco2024higher}. In our case, the SEN domain features an internal length $\eps = 0.01$ mm and an initial pre-crack located at $y = 0$ and $x \in [0, \,L]=[0, \,0.5]$ mm (see Figure \ref{fig:geometry_bc_d}) modeled with the IPF technique. The specimen is loaded with a quasi-static displacement-controlled history consisting of 20 loading steps, with a minimum displacement $u_{\text{min}} = 6\cdot10^{-3}$ mm, a uniform increment $\Delta u = 3\cdot10^{-4}$ mm, and a maximum displacement $u_{\text{max}} = 12\cdot10^{-3}$ mm. {Also in this case, for the $AT_1$ models the considered mesh sizes comprise $h = R_* \eps/4$ and $h = R_* \eps/8$, whereas for the $AT_2$ formulations we examine $h = \eps/2$ and $h = \eps/4$.}

\begin{figure}[h!!]
	\begin{center}
		\includegraphics[width=13cm]{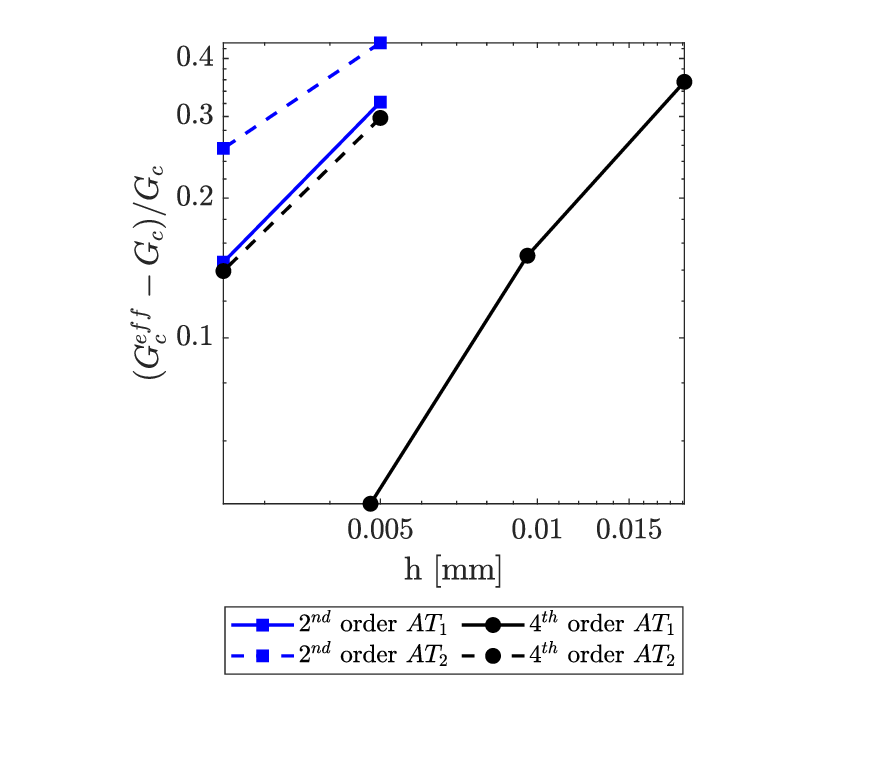}
	\end{center}
	\vspace{-20pt}
	\caption{SEN shear test: accuracy analysis considering all the $AT$-models represented in double logarithmic scale.}
	\label{SEN shear accuracy vari modelli}
\end{figure}
%\NOTE{AP: Luigi, mi cambieresti la notazione della legenda nella Figure 9 per uniformare (->$AT_1$, $AT_2$)? Thx.}

We study the accuracy trend reported in Figure \ref{SEN shear accuracy vari modelli} and, again, we can observe that the novel $AT_1$ model {showcases a higher slope when compared to the $AT_2$ fourth-order model} and comparable accuracy with second-order models allowing for a significant reduction of the computational effort (about 90\% of the degrees of freedom when {using a mesh size {\textit{h} = $0.0096$ mm}).\bl
{ Additionally, we highlight that further refinement of the mesh preserves the convergence trend of the error, that decreases from 15.30\% to 4.39\%, confirming that also for this test in the case of comparable meshes the proposed $AT_1$ fourth-order model exhibits an error ten times lower than the ones provided by other models.} In Figure \ref{crack patter SEN shear novel AT}, we compare the crack patterns obtained with both the proposed $AT_1$ fourth-order and $AT_1$ second-order model, remarking that our fourth-order $AT_1$ model provides a crack pattern which exhibits a kink due to the utilized coarse mesh. However, this tests confirms once again comparable values of relative error on the toughness (see Table \ref{Tab: errors shear AT}) when comparing our $AT_1$ fourth-order formulation against the $AT_1$ low-order counterpart.
\begin{table}[!h]  \centering
	\caption{SEN shear test: AT1 model relative percent error on the toughness.}\label{Tab: errors shear AT}    
	\renewcommand\arraystretch{1.5}
	\renewcommand\tabcolsep{15pt}
	\begin{tabular}{ c c c }
		\hline
		{\textit{h}} & {{2$^{\text{nd}}$ order}} & {4$^{\text{th}}$ order}  \\[-4pt]
		{[mm]} & {[\%]} & {[\%]} \\
		\hline 
		$R_*\eps/4$  & 32.22 &  15.30  \\ 
		$R_*\eps/8$ & 14.56  & 4.39     \\ 
		\hline
		CR [-] & 1.15 & 1.83 \\
		\hline
		%  \botrule 
	\end{tabular} 
\end{table}
\begin{figure}[h!!]
	\begin{center}
		\begin{subfigure}[b]{0.3\textwidth}
			\includegraphics[width=4cm]{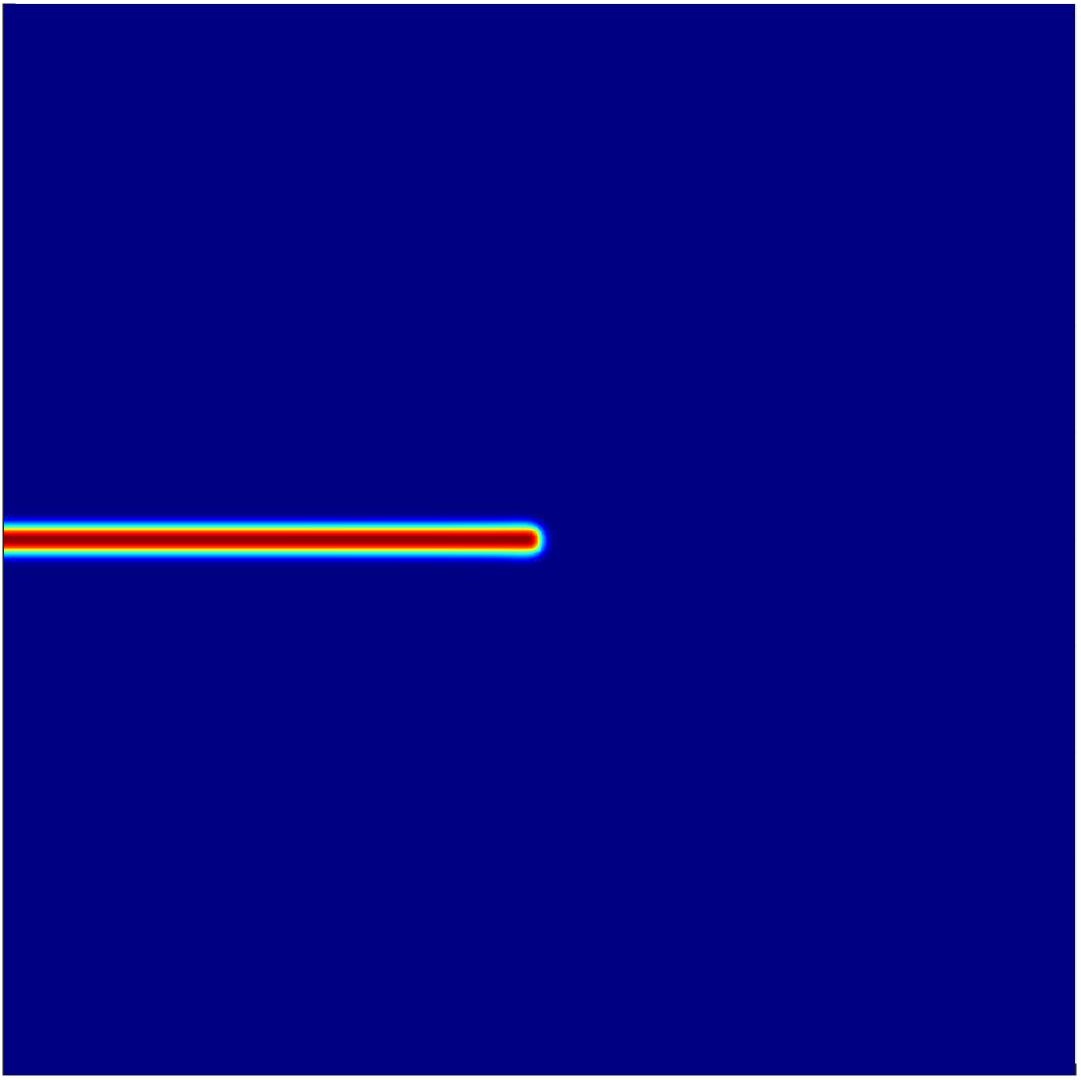}
			\caption{Step 1:  \\initial condition.}
		\end{subfigure}
		%	\hfill
		\begin{subfigure}[b]{0.3\textwidth}
			\includegraphics[width=4cm]{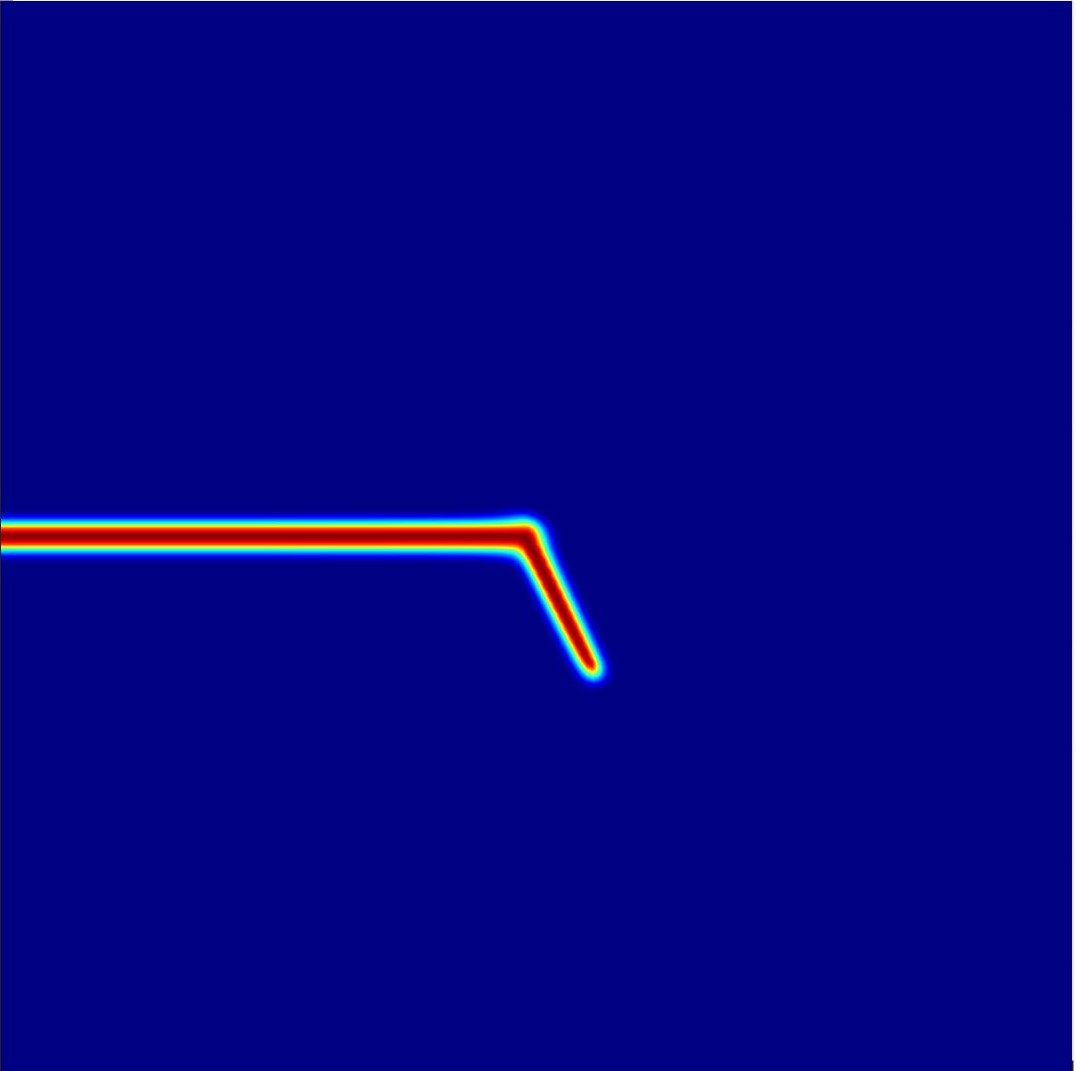}
			\caption{Step 17: fracture \\propagation in the domain. }
		\end{subfigure}
		%   \hfill
		\begin{subfigure}[b]{0.3\textwidth}
			\includegraphics[width=4.7cm]{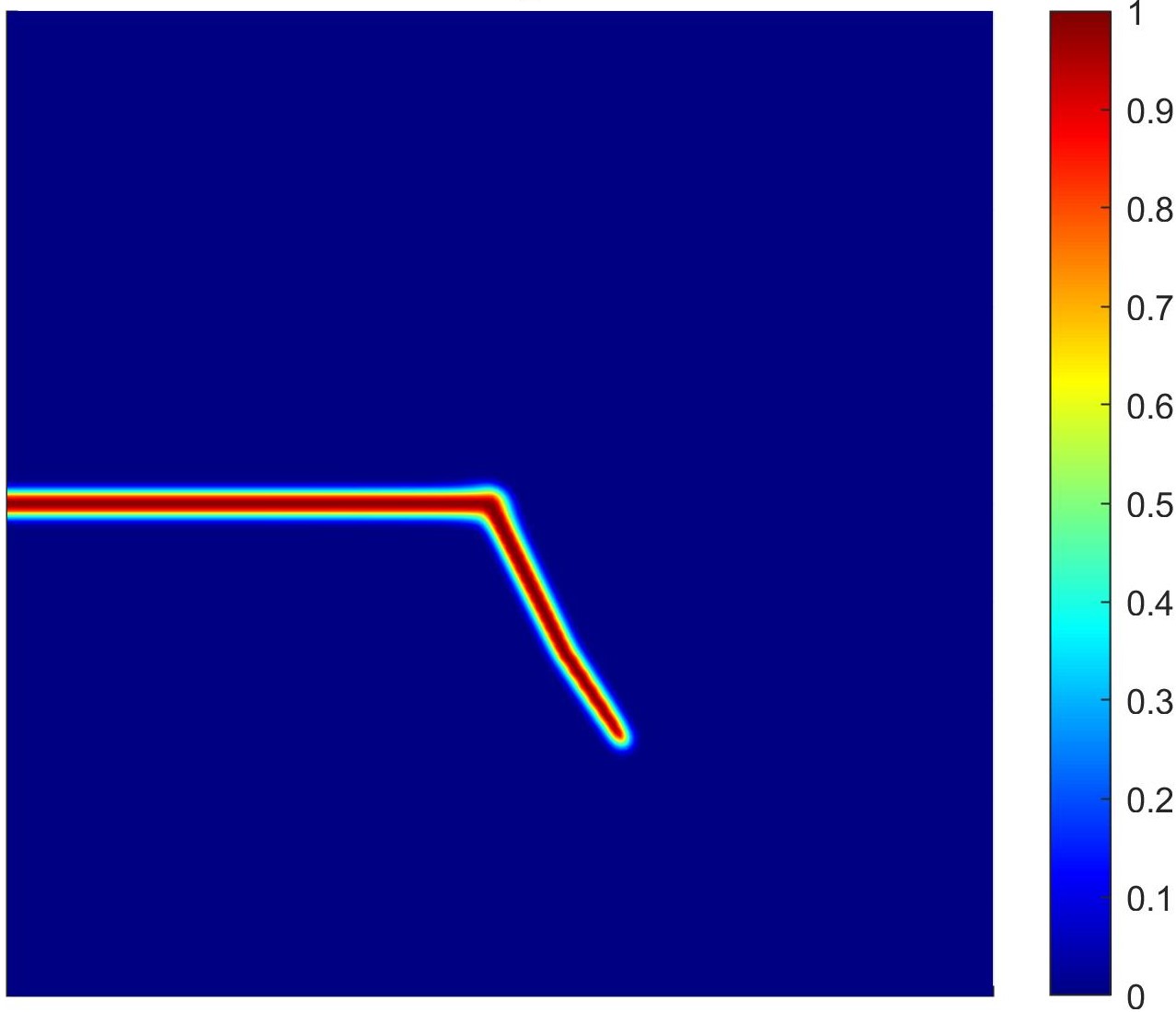}
			\caption{Step 21: \\last loading step.}
		\end{subfigure}
		\begin{subfigure}[b]{0.3\textwidth}
			\includegraphics[width=4cm]{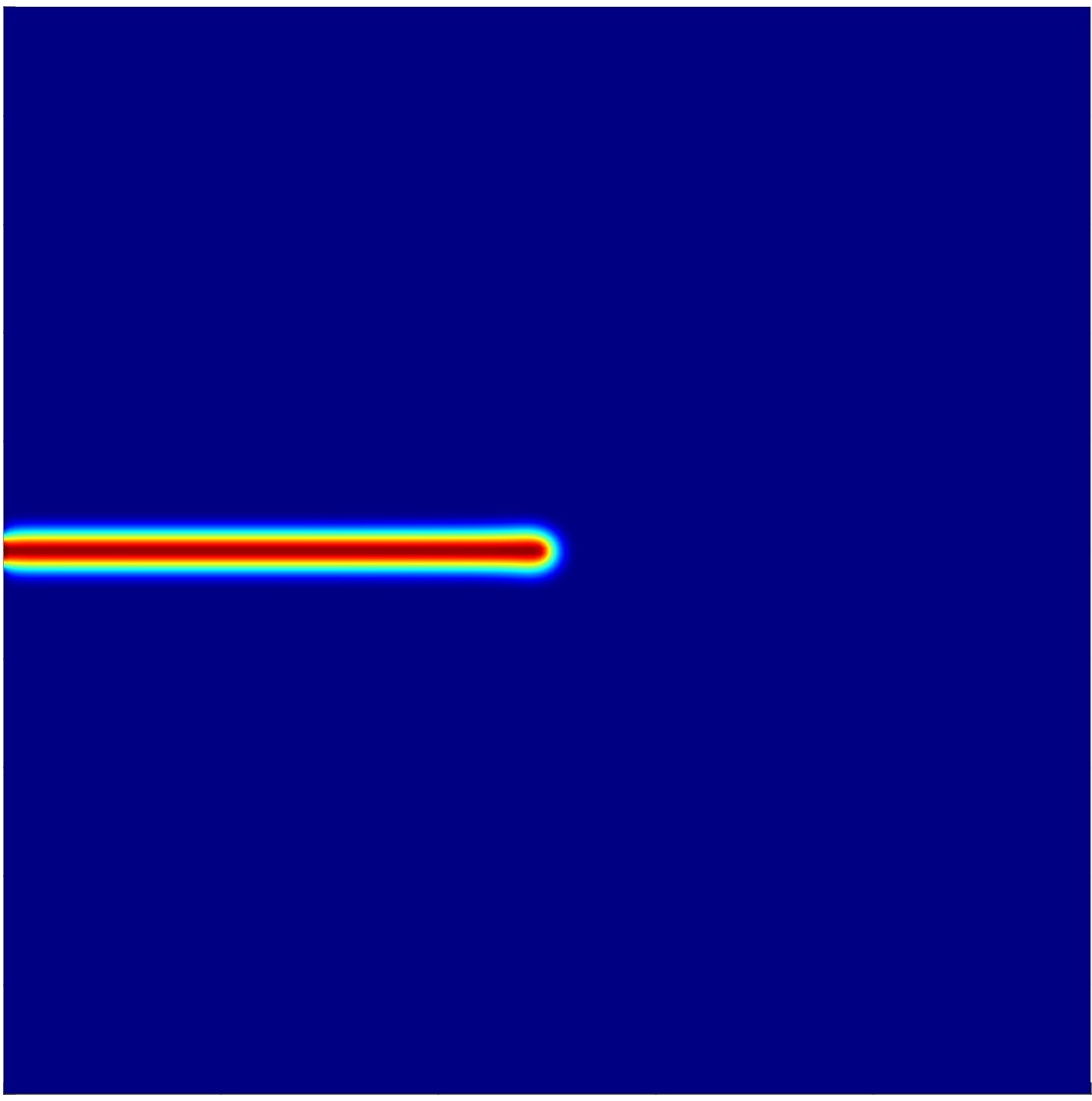}
			\caption{Step 1:  \\initial condition.}
		\end{subfigure}
		%	\hfill
		\begin{subfigure}[b]{0.3\textwidth}
			\includegraphics[width=4cm]{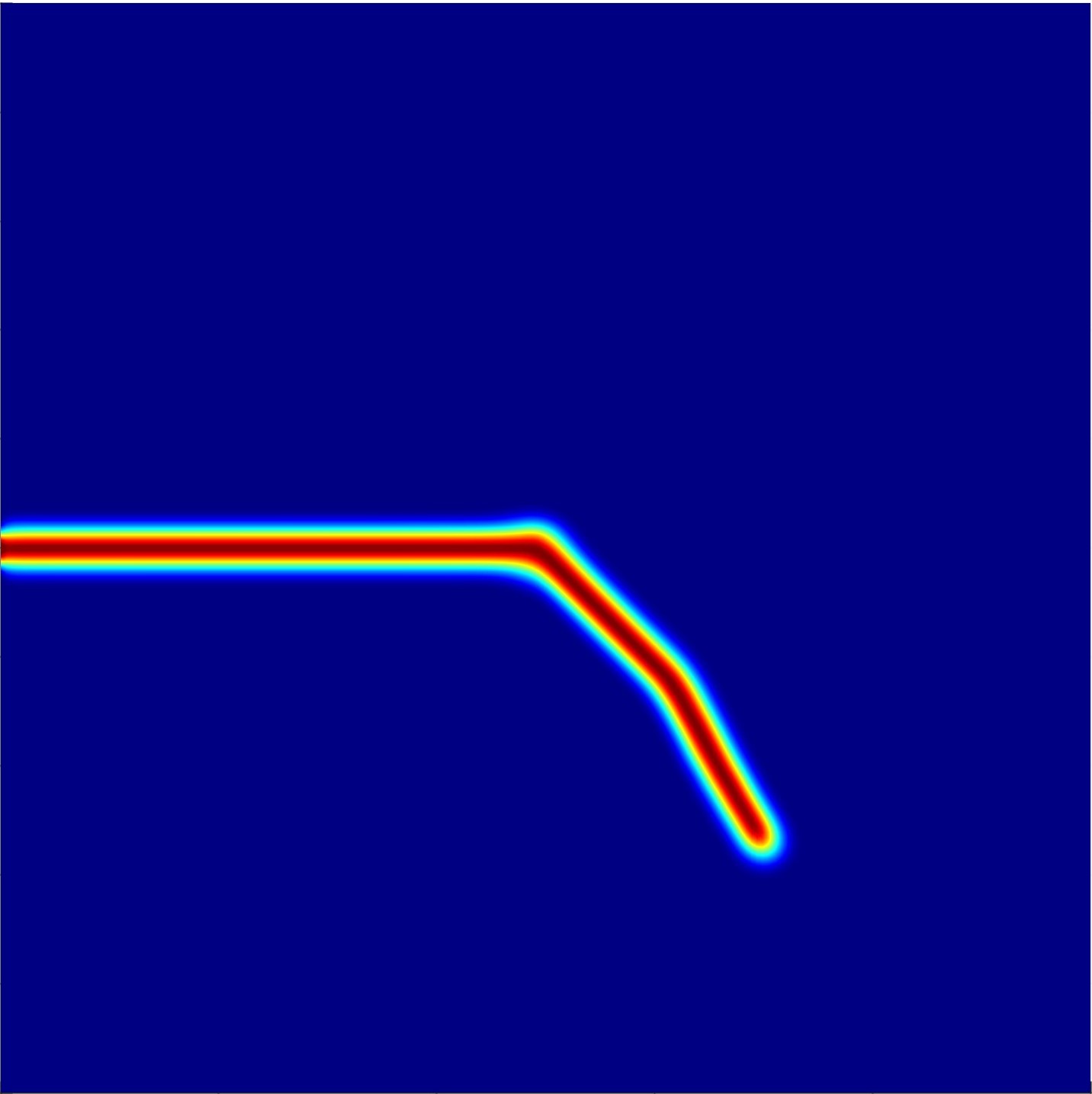}
			\caption{Step 17: fracture \\propagation in the domain. }
		\end{subfigure}
		%   \hfill
		\begin{subfigure}[b]{0.3\textwidth}
			\includegraphics[width=4.7cm]{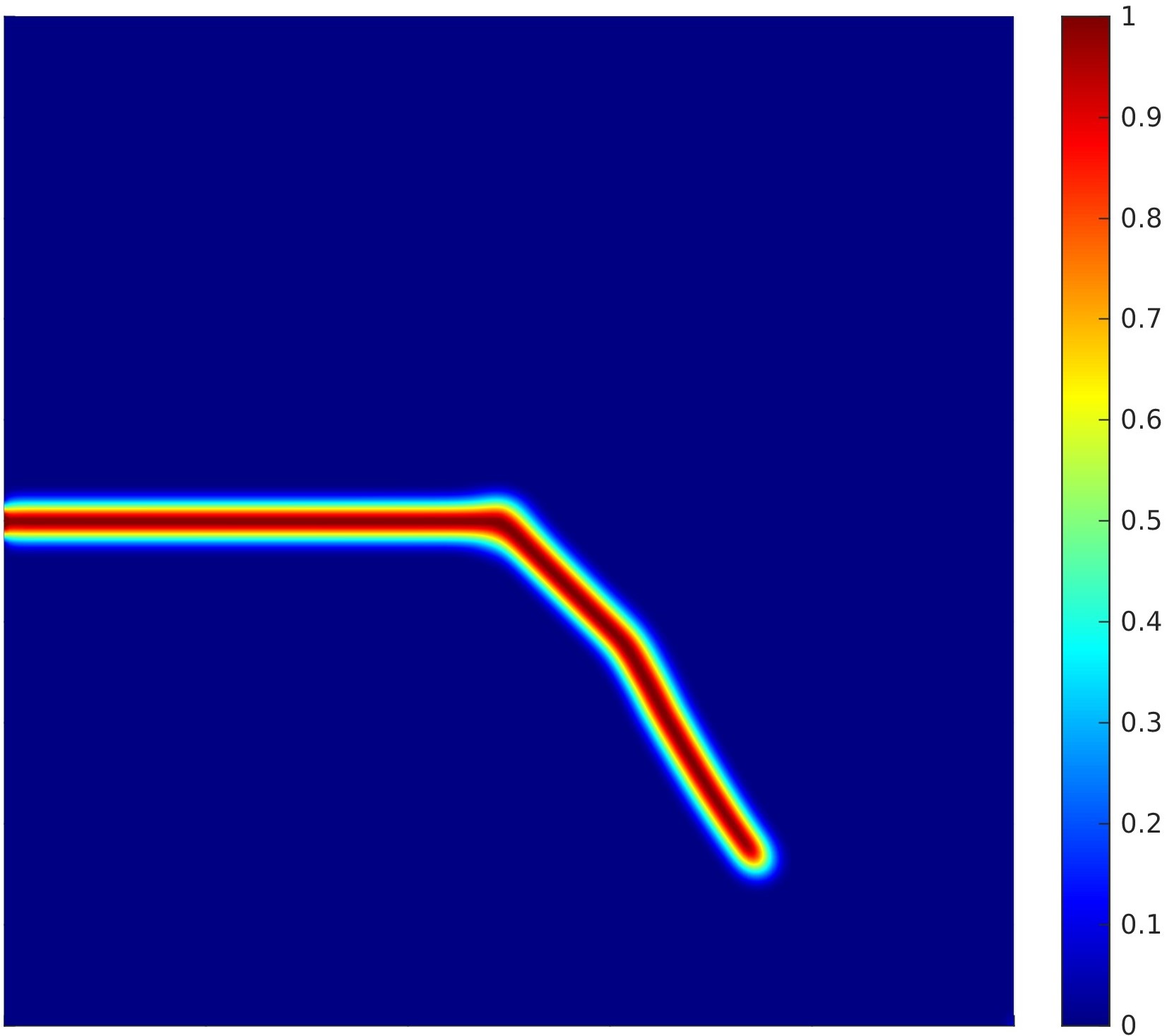}
			\caption{Step 21: \\last loading step.}
		\end{subfigure}
	\end{center}
	\vspace{-15pt}
	\caption{SEN shear test: crack patterns of three different loading steps for the $AT_1$ second-order (top row) and  $AT_1$ fourth-order model (bottom row) considering $ h= 0.005 $ mm.}
	\label{crack patter SEN shear novel AT}
\end{figure}
\section{Parametric study {on the coefficient $\rho$}}
\label{sec: sensitivity study}
%In section \refeq{sec:numerical test}, we propose for the new AT1 model to choose the mesh size according to $R$ parameter, which is a function of the coefficients that characterize the energy dissipation functional. Thus, 
In this section we study how the coefficient $\v$, that weights the high-order term of the {phase-field} energy, affects the accuracy of the effective toughness as well as the elastic limit. {To this end, we} consider the following values: $\v = 2^i$ for $i=-4,...,4$.

\subsection{Elastic limit comparison}
Considering a variable coefficient associated to the laplacian term influences the value of the regularization parameter $c_\v$ from which depends the theoretical value of elastic limit. As in \S \ref{s.puretrac}, we consider a pure traction test (see Figure \ref{fig:geometry_bc_a}) for the proposed $AT_1$ model. The results are summarized in Table \ref{Tab: errors pure traction c vary}{, highlighting a minimum value in terms of relative error for $\v= 1$,} that appears to be a good choice to perform all analyses. Furthermore, it can be observed that the value of $R_*$ is increasing with respect to $\v$ and, consequently, the mesh size $h \sim R_* \eps$ is increasing with respect to $\v$. %following Eq. \eqref{Analytical solution}, the theoretical elastic limit decreases. 
Thus, other possible choices of $\v$ balance the trade off between accuracy and computation time.

%\medskip
%\NOTE{con un errore sotto l'$1\%$ (se capisco bene) mi sembra che vadano tutti ugualmente bene} 

\begin{table}[h!] \centering
	\caption{{Pure traction test:} relative percentage error on the elastic limit for {a fixed mesh corresponding to the reference case $\v=1$ (i.e., $R_* = 3.83$).}}\label{Tab: errors pure traction c vary} 
%	\vspace{-5pt}
	\renewcommand\arraystretch{1.25}
	\renewcommand\tabcolsep{8pt}
\begin{tabular}{c c c c c c}
	\hline
	\centering
	{{$\v$}} & {{{$R_*$}}} & {c$_\v$} & {$\sigma_{\text{th}}$} & {$\sigma_{c}$} & {{R.err}} \\
	\hline
	{[-]} & {[-]} & \textbf{[-]} & {[kN/mm$^2$]} & {[kN/mm$^2$]} & \textbf{[\%]}\\
	\hline
	16 & 7.1041 & 7.7811 & 1.0140 & 1.0199 & 0.5871  \\ 
	8 & 6.0364 & 6.6769 & 1.0946 & 1.0973 & 0.2440  \\ 
	4 & 5.1514 & 5.7717 & 1.1773 & 1.1784 & 0.0929 \\ 
	2 & 4.4230 & 5.0369 & 1.2603 & 1.2605 & 0.0197  \\ 
	1 & 3.8300 & 4.4485 & 1.3410 & 1.3410 & 0.0026    \\ 
	1/2& 3.3554 & 3.9852 & 1.4168 & 1.4167 & 0.0115  \\ 
	1/4 & 2.9847 & 3.6281 & 1.4849 & 1.4848 & 0.0101 \\ 
	1/8 & 2.7045 & 3.3593 & 1.5432 & 1.5430 & 0.0123\\ 
	1/16 & 2.4998 & 3.1615 & 1.5907 & 1.5900 & 0.0468 \\ 
	\hline
\end{tabular}
\end{table}
%\NOTE{AP: Table 6: uniformare le cifre significative a categoria. $R$ e' $R_*$?}
\subsection{Toughness accuracy}
To evaluate the toughness accuracy as a function of $\v$, the SEN tensile and the SEN shear tests are reconsidered, providing respectively a mode I and mode II fracture benchmark. {All analyses feature mesh size $h(R_*) = R_* \eps/4$, where $R_* = 3.83$, namely the value in the reference case $\v = 1$.}
%\NOTE{AP: Ho capito correttamente quello che intendi rispetto alla mesh Luigi nelle modifiche fatte?}
\begin{figure}[h!!]
	\begin{center}
		\begin{subfigure}[b]{0.18\textwidth}
		\includegraphics[width=2.5cm]{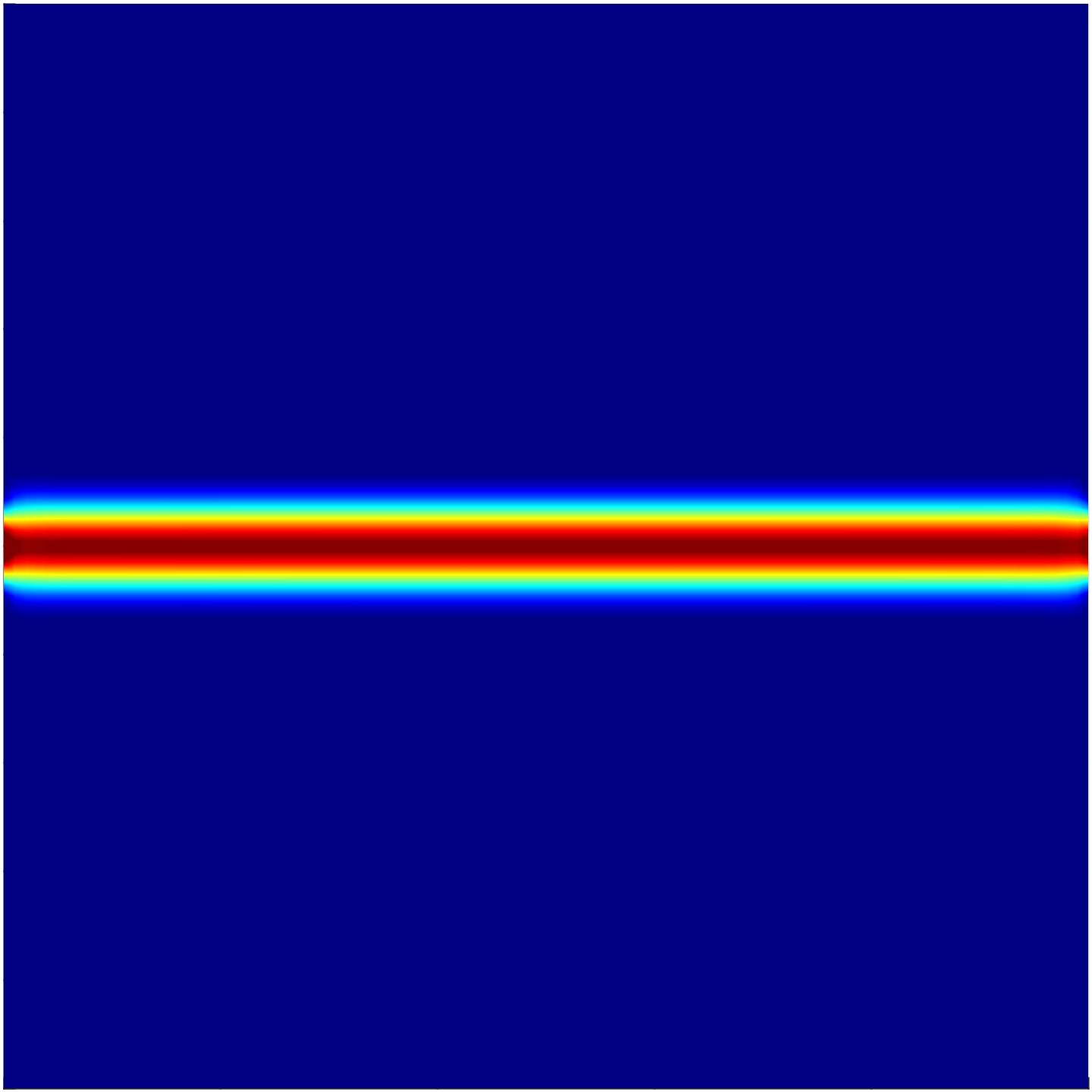}
		\caption{$\v$= 16.}
	\end{subfigure}
	\begin{subfigure}[b]{0.18\textwidth}
		\includegraphics[width=2.5cm]{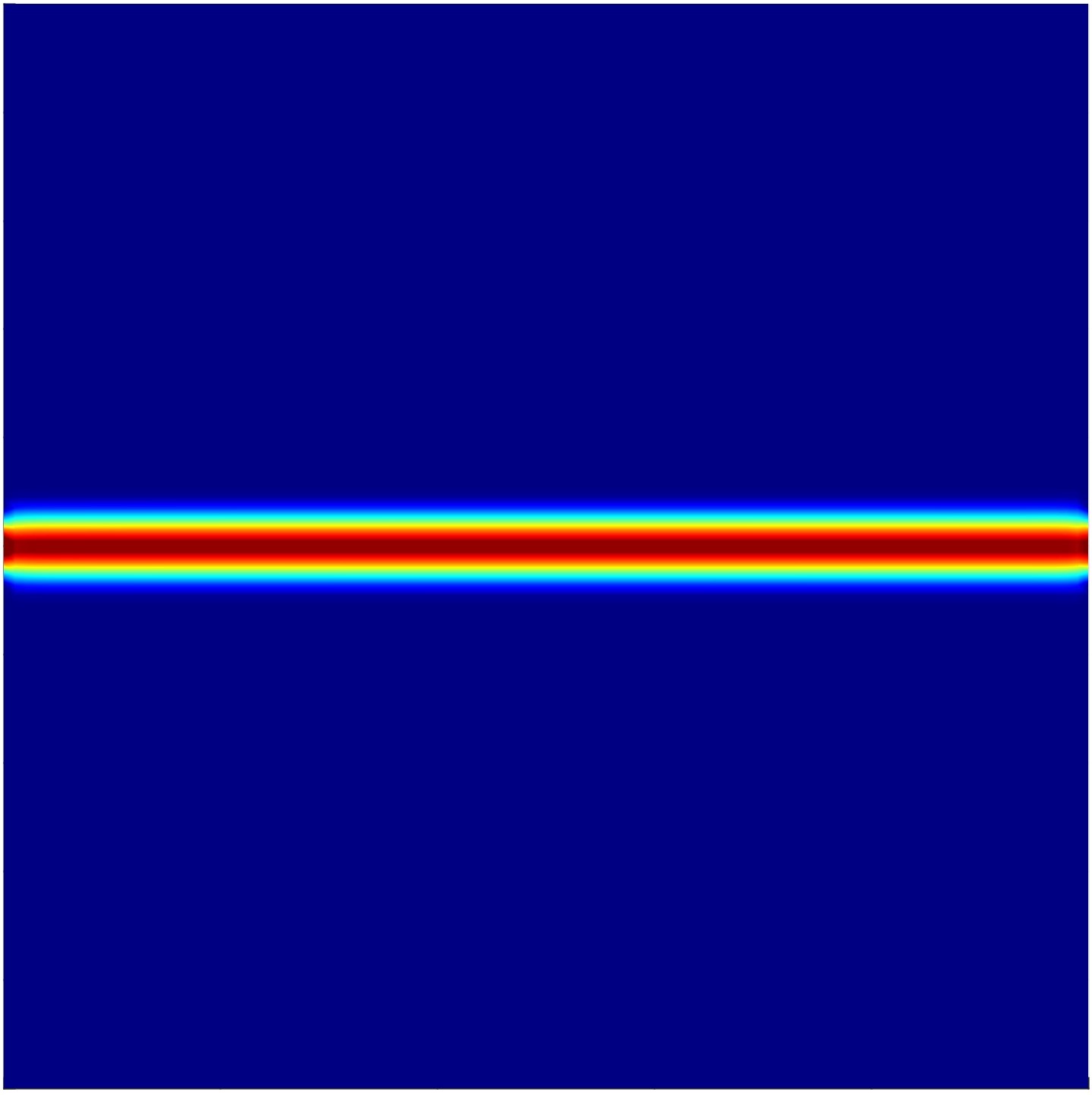}
		\caption{$\v$= 4.}
	\end{subfigure}
	\begin{subfigure}[b]{0.18\textwidth}
		\includegraphics[width=2.5cm]{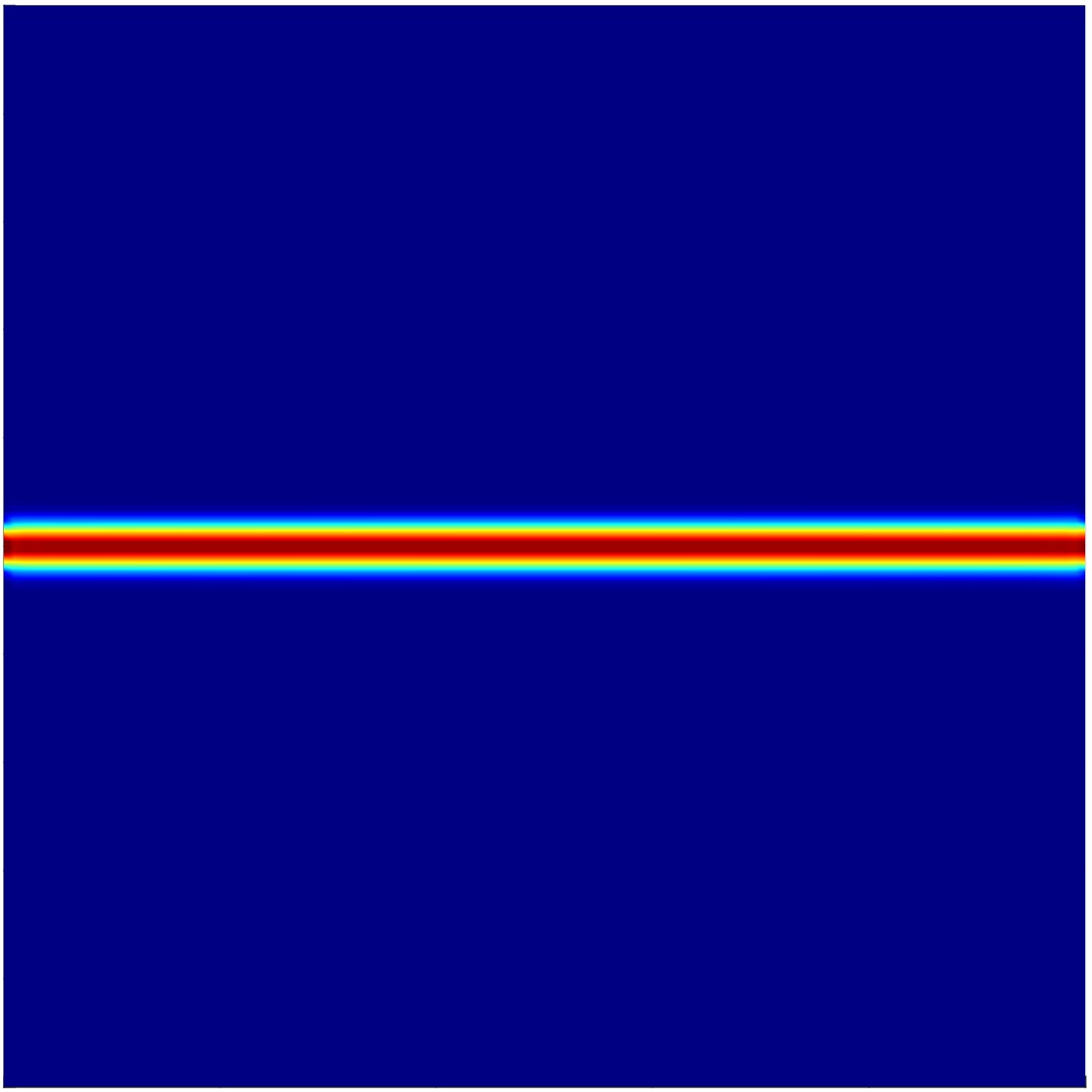}
		\caption{$\v$= 1.}
	\end{subfigure}
	\begin{subfigure}[b]{0.18\textwidth}
		\includegraphics[width=2.5cm]{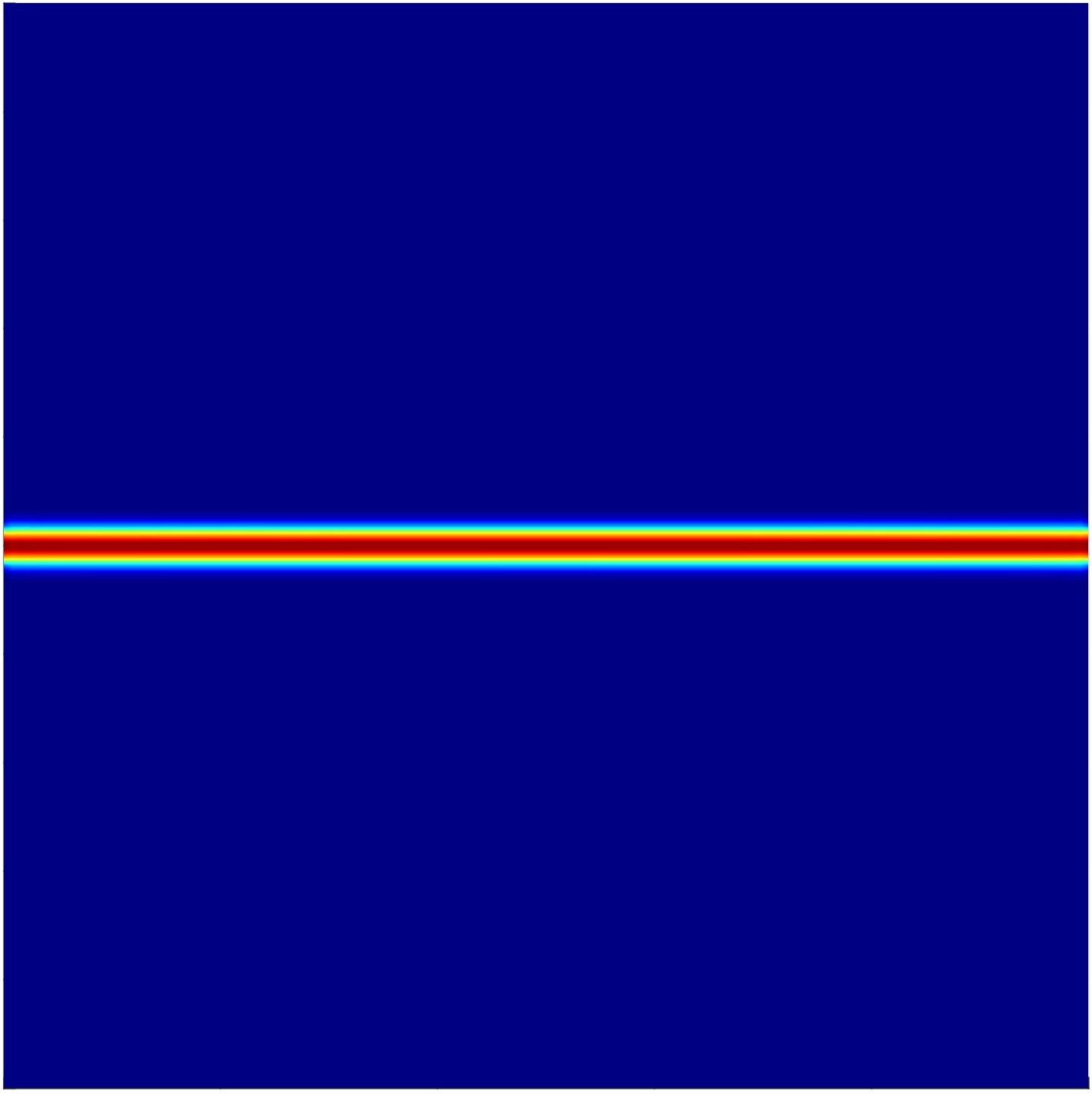}
		\caption{$\v$= 1/4.}
	\end{subfigure}
	\begin{subfigure}[b]{0.18\textwidth}
		\includegraphics[width=2.5cm]{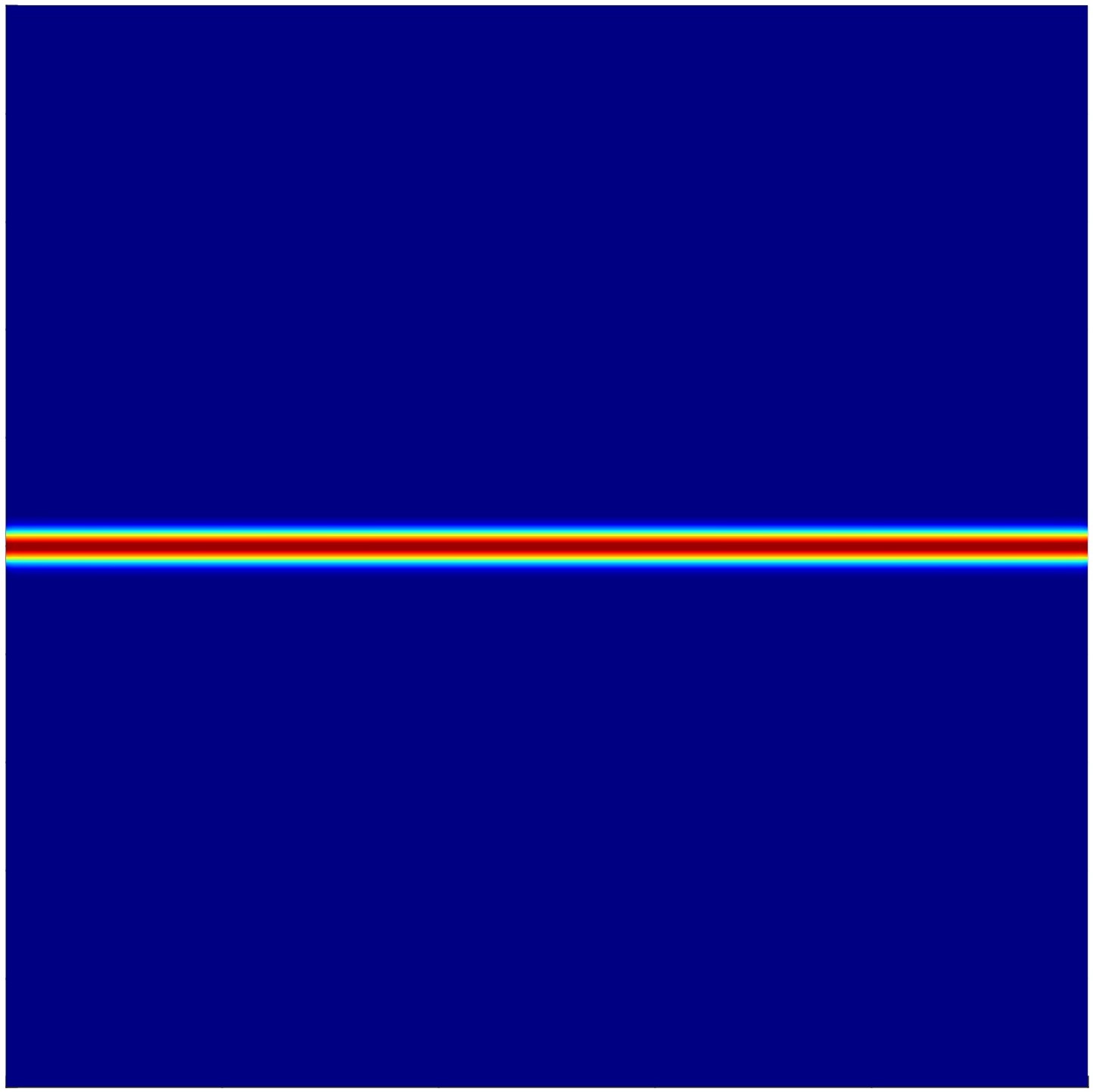}
		\caption{$\v$= 1/16.}
	\end{subfigure} 
	\end{center}
\vspace{-15pt}
	\caption{SEN tensile test: influence of the variation of the $\v$ coefficient on the final crack pattern for a fixed mesh size.}
\label{Crack coeff vary SEN T}
\end{figure}
{In Figure \ref{Crack coeff vary SEN T}, it can be observed that for the SEN traction test as the $\v$ coefficient decreases, the fracture width is reduced. The same behavior can be noticed also for the SEN shear benchmark in Figure \ref{Crack coeff vary SEN S}. Additionally, we notice in this latter case that the crack length corresponding to the last loading step reduces as $\v$ decreases.}

\begin{figure}[h!!]
	\begin{center}
		\begin{subfigure}[b]{0.18\textwidth}
			\includegraphics[width=2.5cm]{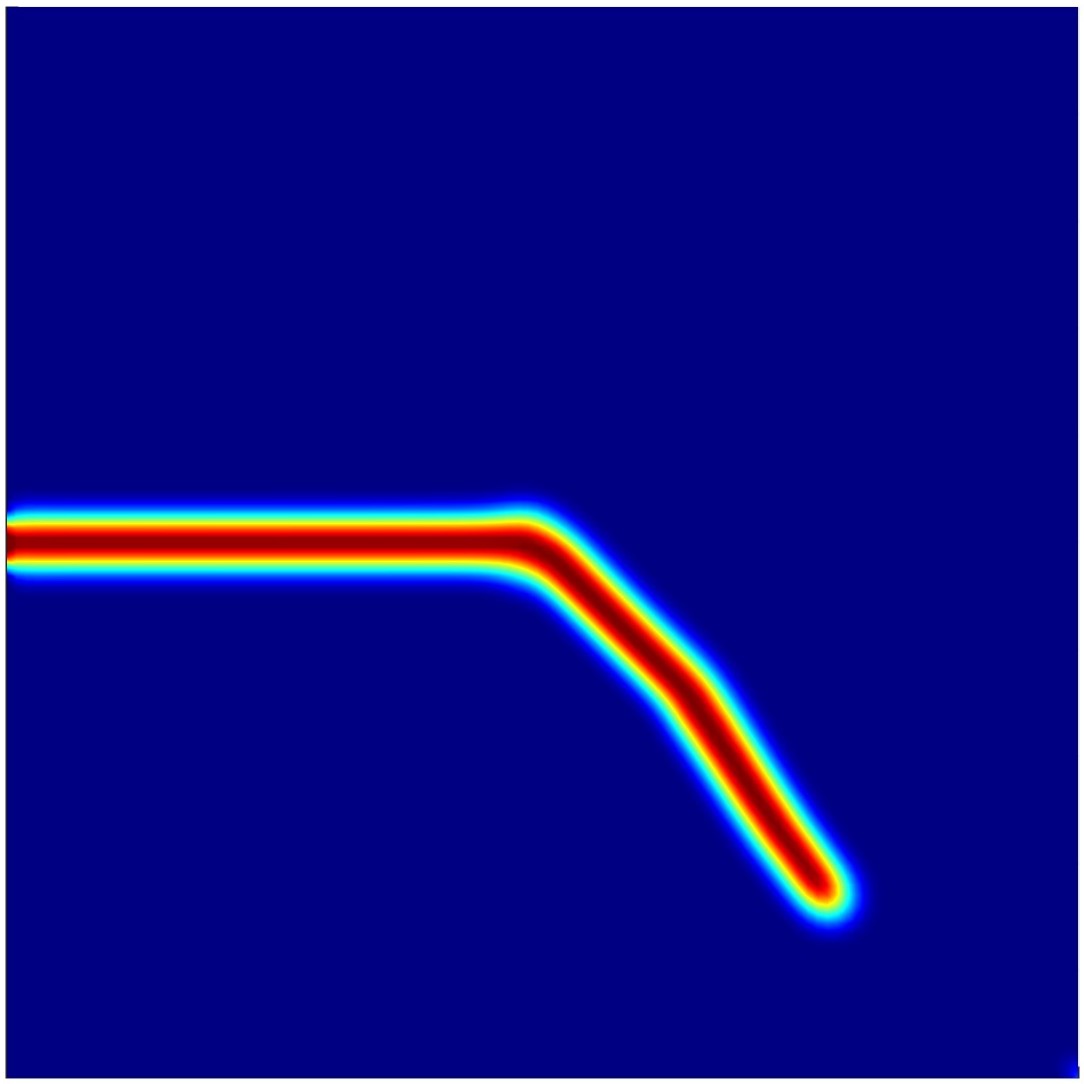}
			\caption{$\v$ = 4.}
		\end{subfigure}
		\begin{subfigure}[b]{0.18\textwidth}
			\includegraphics[width=2.5cm]{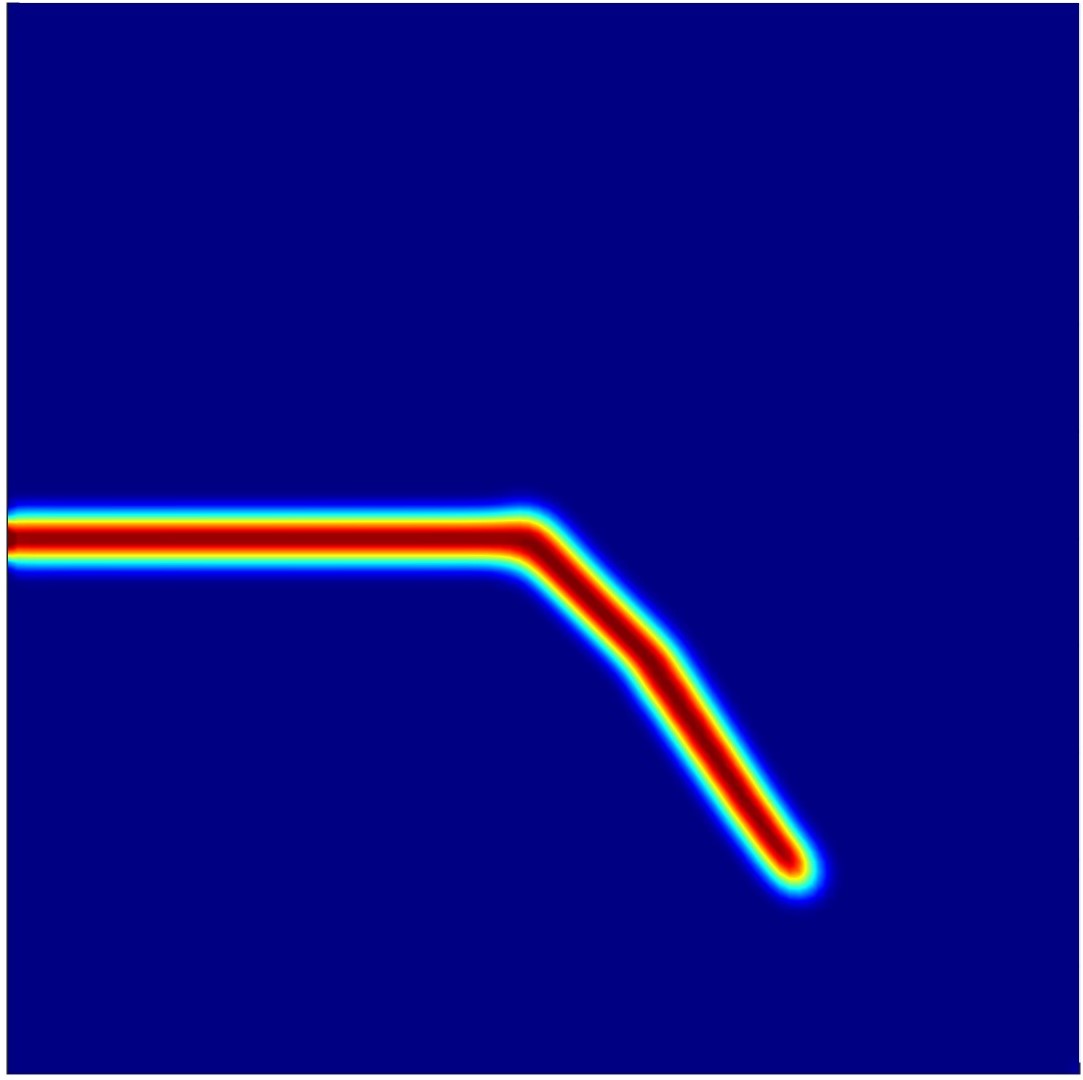}
			\caption{$\v$= 2.}
		\end{subfigure}
		\begin{subfigure}[b]{0.18\textwidth}
			\includegraphics[width=2.5cm]{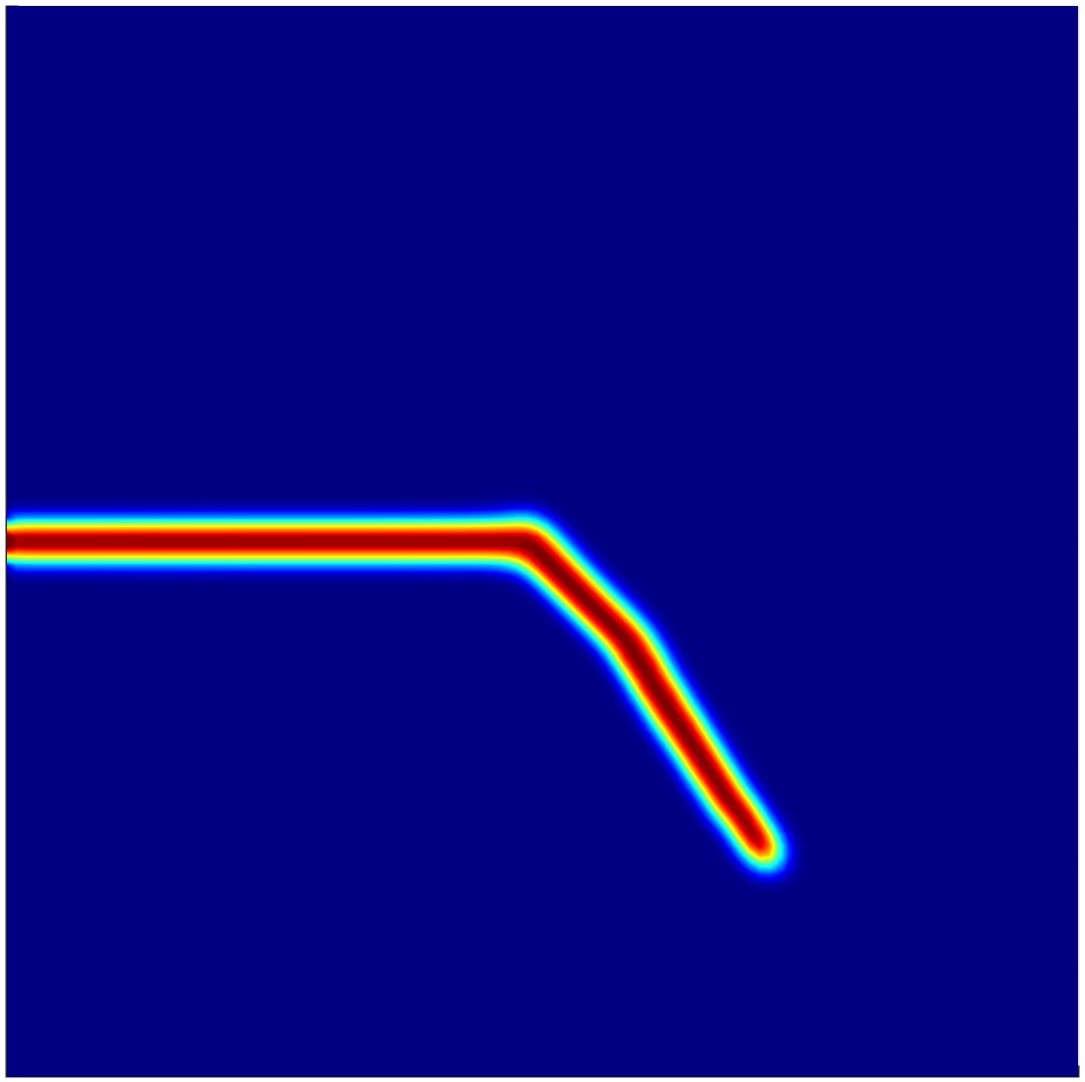}
			\caption{$\v$= 1.}
		\end{subfigure}
		\begin{subfigure}[b]{0.18\textwidth}
			\includegraphics[width=2.5cm]{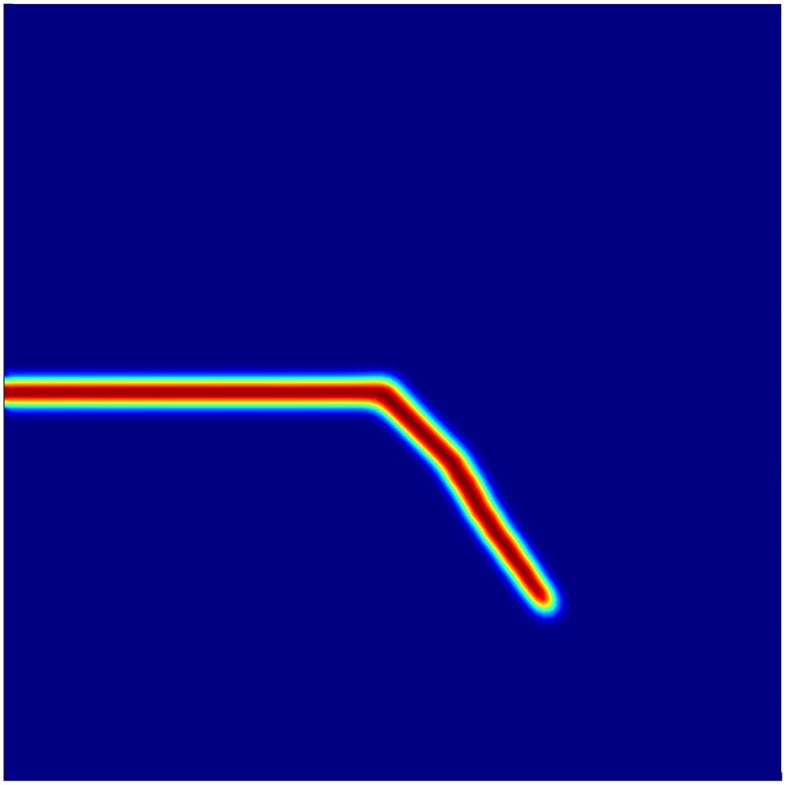}
			\caption{ $\v$= 1/2.}
		\end{subfigure}
		\begin{subfigure}[b]{0.18\textwidth}
			\includegraphics[width=2.5cm]{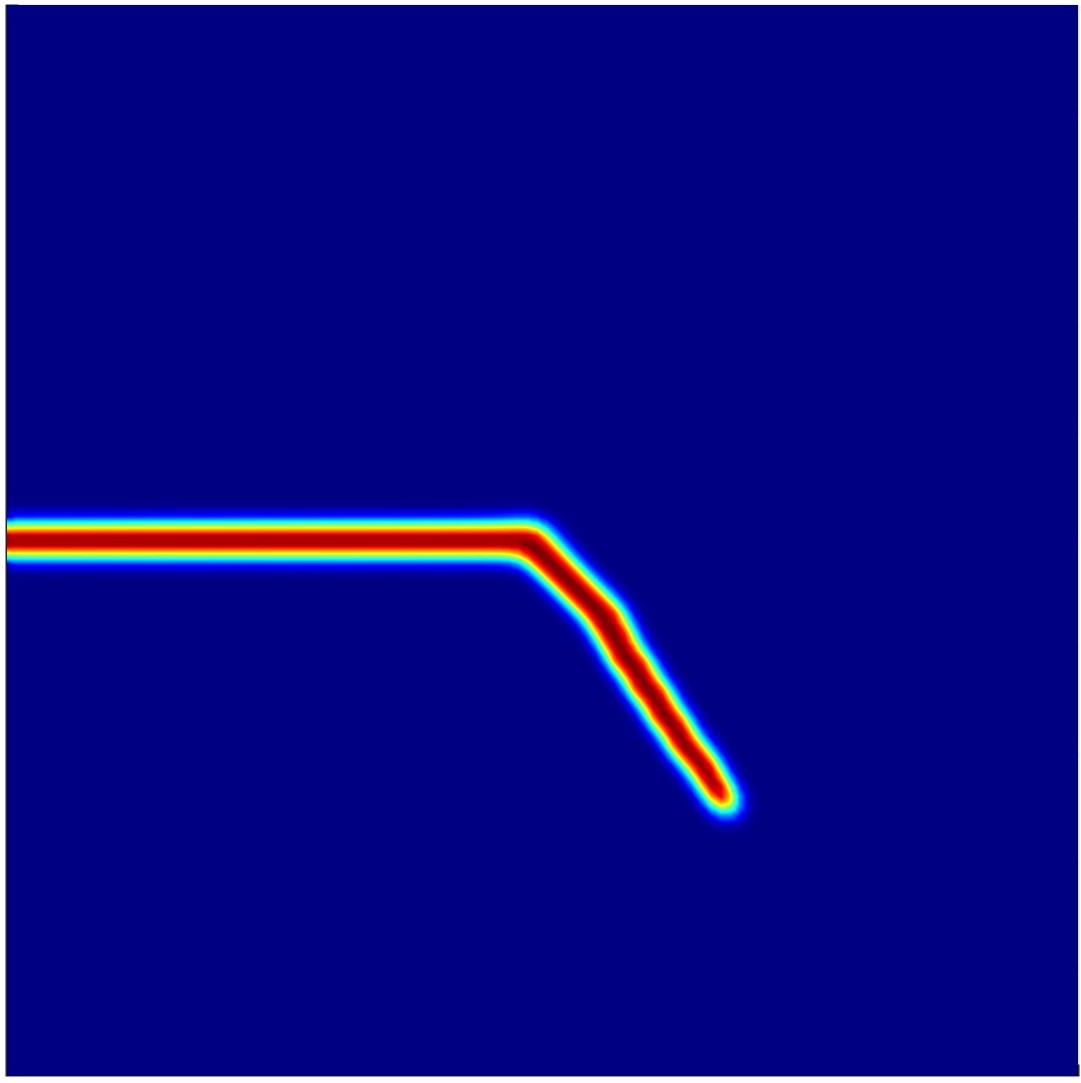}
			\caption{$\v = 1/4$.}
		\end{subfigure} 
	\end{center}
	\vspace{-15pt}
	\caption{SEN shear test: influence of the variation of the $\v$ coefficient on the final crack pattern for a fixed mesh size.}
	\label{Crack coeff vary SEN S}
\end{figure}

%The same variation on $c$ for the SEN shear test instead leads to a fracture that propagates more in the domain (Figure \ref{Crack coeff vary SEN S}). %
\begin{figure}[h!!]
	\begin{center}
		\begin{subfigure}[b]{0.5\textwidth}
			\includegraphics[width=9.5cm]{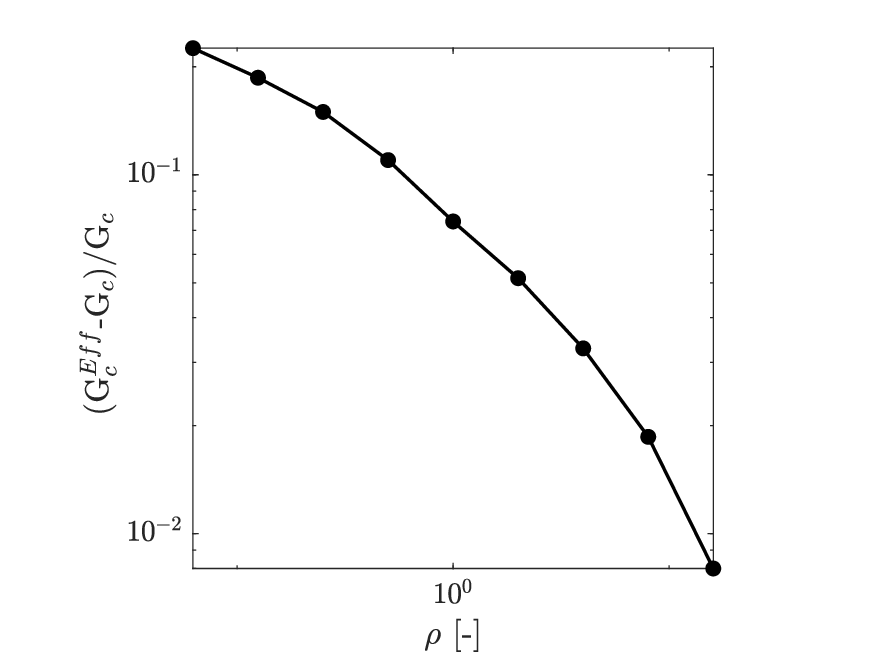}
			\caption{SEN tensile test.}\label{errore c vary SEN_a}
		\end{subfigure}
		\begin{subfigure}[b]{0.49\textwidth}
			\includegraphics[width=9.5cm]{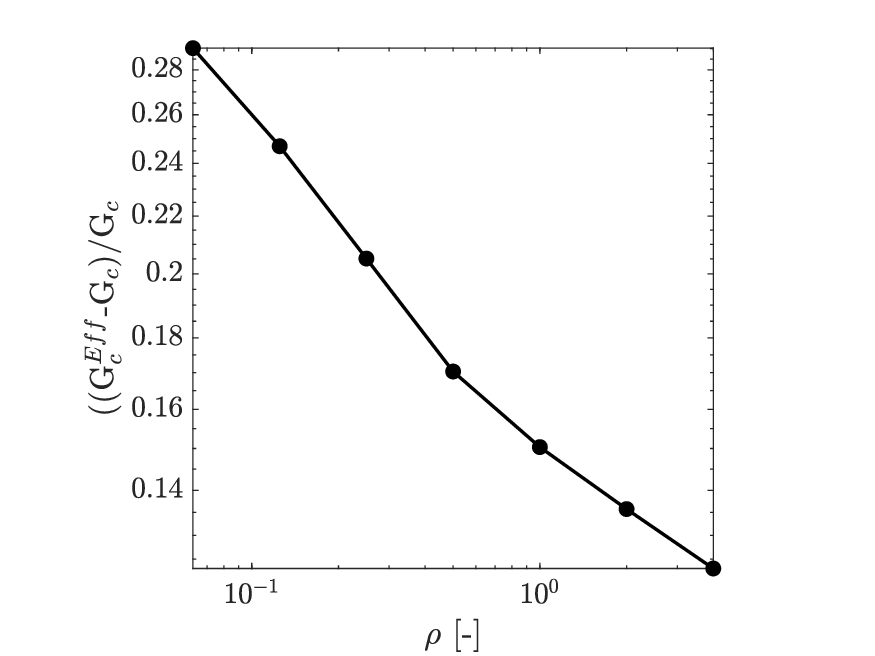}
			\caption{SEN shear test.}\label{errore c vary SEN_b}
		\end{subfigure}
%		\vspace{-5pt}
	\end{center}
	\caption{Trend of the relative error \bl on the toughness for the SEN traction (left) and SEN shear (right) tests as a function of the $\v$ coefficient for a fixed mesh level.}
	\label{errore c vary SEN}
\end{figure}
%\NOTE{AP: Figura 13: $\rho$ [-], space is missing.}
{In Figure \ref{errore c vary SEN_a}, we highlight that for the SEN tensile test the relative error on the toughness is reducing while the $\v$ laplacian coefficient  increases. This is due to the fact that by increasing the coefficient related to the higher-order term, the values of $R_*$ increase and, consequently, the profile exhibits a wider shape. This observation underlines the influence of the parameter $R_*$ on the mesh size that will be the object of an in-depth study in the next section. Also in the case of the SEN shear test (see Figure \ref{errore c vary SEN_b}) the relative error on the toughness decreases as the coefficient $\v$ increases. However, for this latter test a significantly high value of the $\v$ coefficient provides an inaccurate crack patter that is far from the reference one in the literature. Therefore, we omit this result in Figure \ref{errore c vary SEN_b} and prefer to represent the error trend only in the interval $\v \in [1/16 - 4]$.}

%allows for lower error values with very large meshes. On the other hand, large values of $\v$ result into a crack pattern with a thicker fracture that moves away from the idea of a sharp fracture (see e.g., Figure \ref{Crack coeff vary SEN T}a). Instead the SEN shear test exhibits an opposite trend with respect to the SEN tensile one (Figure \ref{errore c vary SEN}b), as we have a reduction of the error once the Laplacian coefficient $\v$ decreases. This is a consequence of how the error is calculated, i.e., by relating the energy dissipated with the length of the fracture. As it can be noticed in  Figure \ref{Crack coeff vary SEN S}, the fracture has a greater propagation once the $c$ coefficient decreases and consequently we have a lower error.

\section{Trade-off between accuracy and {computational} cost \label{trade}}
{As the proposed $AT_1$ high-order model allows to perform simulation with higher accuracy, we focus in this section on assessing the computational cost in terms of necessary number of control points for a fixed and comparable level of accuracy with respect to the identified reference models in the literature. More specifically, in §\ref{sec: degree saving} we set an acceptable range of error and study how much the mesh can be coarsened, whereas in §\ref{influence of  $R_*$} we study the $R_*$ value influence on the mesh size and the accuracy.}

\subsection{$AT_1$ - $AT_2$ comparisons for different mesh sizes}
\label{sec: degree saving}
{While performing numerical simulations, it is of primal importance to balance the trade-off between providing reliable results and keeping the computational effort of the analyses to an acceptable level.
Thus, we herein focus again on two of the benchmarks discussed in §\ref{sec:numerical test}, namely the SEN shear and traction tests, and we assess the computational gain in terms of number of utilized control points for a fixed level of error on the fracture toughness.
First, we focus on the SEN shear benchmark, which is modeled using a uniform quadratic $C^1$ IgA mesh and we observe in Table \ref{Tab: cp saved shear} that, for a fixed toughness error of about 15\%, the $AT_1$ high-order model needs 11,664 control points, whereas all other $AT$-models\footnote{{We remark that the $AT_2$ second-order model has been excluded from this discussion because it is out of scale. In fact, as Figure \ref{SEN shear accuracy vari modelli} clearly highlights, it would need even more control points than the considered $AT_1$ second-order and  $AT_2$ fourth-order models.}} utilize 161,604 degrees of freedom. Consequently, we have in this case a saving in terms of control points of $\sim$93\%. Then, we analyze the SEN tensile test, which considers a non-uniform knot vector in the case of all $AT_2$ and the $AT_1$ low-order models, whereas, even considering a uniform mesh for the $AT_1$ fourth order model, we obtain a saving of $\sim$43\% in terms of control points (see Table \ref{Tab: cp saved tensile}) to reach an error level of $\sim$10\%.}

\begin{table}[!h] \centering
	\caption{{SEN shear test: number of control points for a fixed value of error on the toughness.}}
	\label{Tab: cp saved shear} 
	\vspace{-5pt}
	\renewcommand\arraystretch{1.5}
	\renewcommand\tabcolsep{8pt}
	\begin{tabular}{ c c c c c }
		\hline {Error} & \multicolumn{2}{c}{{2$^{\text{nd}}$ order}} & \multicolumn{2}{c}{4$^{\text{th}}$ order}\\
		\hline {[\%]} &  {$AT_1$} & {$AT_2$} & {$AT_1$} & {$AT_2$} \\
		\hline 
		$\sim$ 15  & 161,604  & - & 11,664 & 161,604\\ 
		\hline
		%  \botrule 
	\end{tabular}
\end{table}

\begin{table}[!h] \centering
\caption{{SEN tensile test: number of control points for a fixed value of error on the toughness.}}
	\label{Tab: cp saved tensile} 
	\vspace{-5pt}
	\renewcommand\arraystretch{1.5}
	\renewcommand\tabcolsep{8pt}
	\begin{tabular}{ c c c c c }
		\hline {Error} & \multicolumn{2}{c}{{2$^{\text{nd}}$ order}} & \multicolumn{2}{c}{4$^{\text{th}}$ order}\\
		\hline {[\%]} &  {$AT_1$} & {$AT_2$} & {$AT_1$} & {$AT_2$} \\
		\hline 
		$\sim$ 10  & 20,100  & 20,100 & 11,664 & 20,100\\ 
		\hline
		%  \botrule 
	\end{tabular}
\end{table}

%\textcolor{red}{Due to the computational limit linked to our homecode, it was necessary for the DCB and SEN tensile test perform the analyses with the AT2 models and the second-order AT1 by locally refining the fracture propagation area, which for these benchmarks is known} \bb a priori, \bl\textcolor{red}{using an anisotropic knot-vector refinement guarantees the minimum number of elements in the internal length to accurately solve the phase-field problem. \mg In Table \refeq{Tab: cp saved tensile} and Table \refeq{Tab: cp saved dcb} is reported the savings in terms of control points and the relative error of the toughness for these tests provides from the use of the novel fourth-order AT1 model with uniform mesh. \bl}

\medskip

\subsection{Influence of the $R_*$ parameter}
\label{influence of  $R_*$}
{In §\ref{sec: sensitivity study}, we investigate how the different values of $R_*$ change the elastic limit (see Table \ref{Tab: errors pure traction c vary}) as well as the crack pattern width. All these considerations are made for a fixed mesh size (referring to the mesh size of the case considering $\v$ = 1, $R_* = 3.83$). However, the mesh-size definition in Equation \eqref{mesh definition} allows to considerably reduce the number of elements of a simulations for higher values of the $\v$ coefficient. Namely, the user has the choice to balance the trade-off between having a very accurate simulation (considering a mesh size equal to the one utilized for the $AT_2$-model analyses) or less accurate results (providing an acceptable level of error) but faster simulations by taking into account the $R_*$ parameter in the mesh size definition. To investigate this latter possibility, we performed the SEN tensile and SEN shear tests changing the mesh size as in Equation \eqref{mesh definition} for our fourth-order formulation and we observed that, for the SEN tensile test (see Figure \ref{Crack coeff vary SEN T mv}), there is not a significant variation in the crack pattern. However, in the SEN shear test case (see Figure \ref{Crack coeff vary SEN S mv}), choosing a mesh size based on the $R_*$ value as in Equation \eqref{mesh definition} leads to a comparable fracture length with a smoother crack pattern. It can be observed also that reducing the coefficient $\v$ the curvature of the fracture makes a difference, reaching the one obtained for the second-order models and the fourth-order $AT_2$, where the coefficient $\v  =  1/16$.}
%\NOTE{AP: Cambiato $c=1/16$ in $\v$.}

\begin{figure}[h!!]
	\begin{center}
		\begin{subfigure}[b]{0.18\textwidth}
			\includegraphics[width=2.5cm]{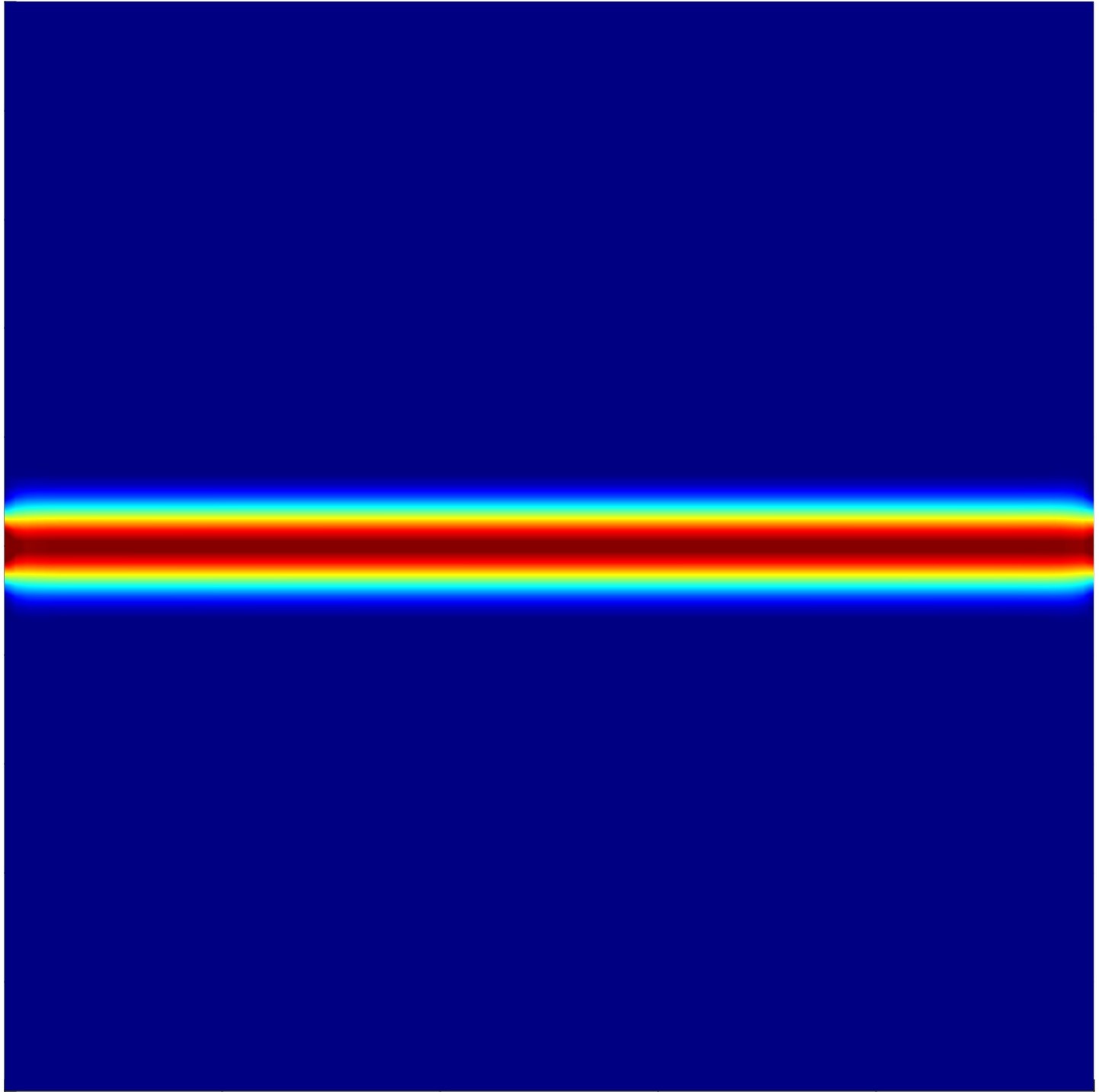}
			\caption{ $\v$= 16.}
		\end{subfigure}
		\begin{subfigure}[b]{0.18\textwidth}
			\includegraphics[width=2.5cm]{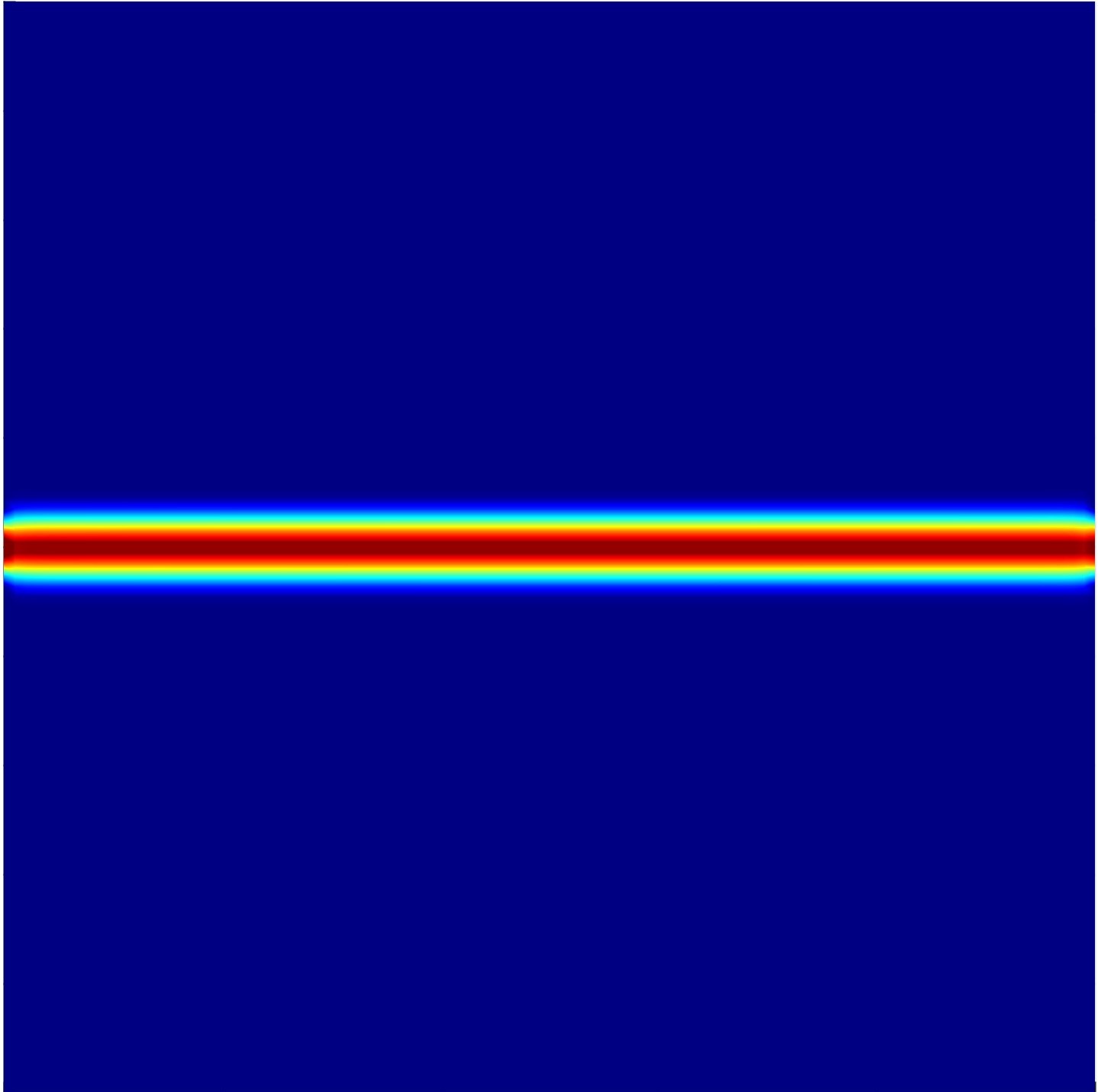}
			\caption{$\v$= 4.}
		\end{subfigure}
		\begin{subfigure}[b]{0.18\textwidth}
			\includegraphics[width=2.5cm]{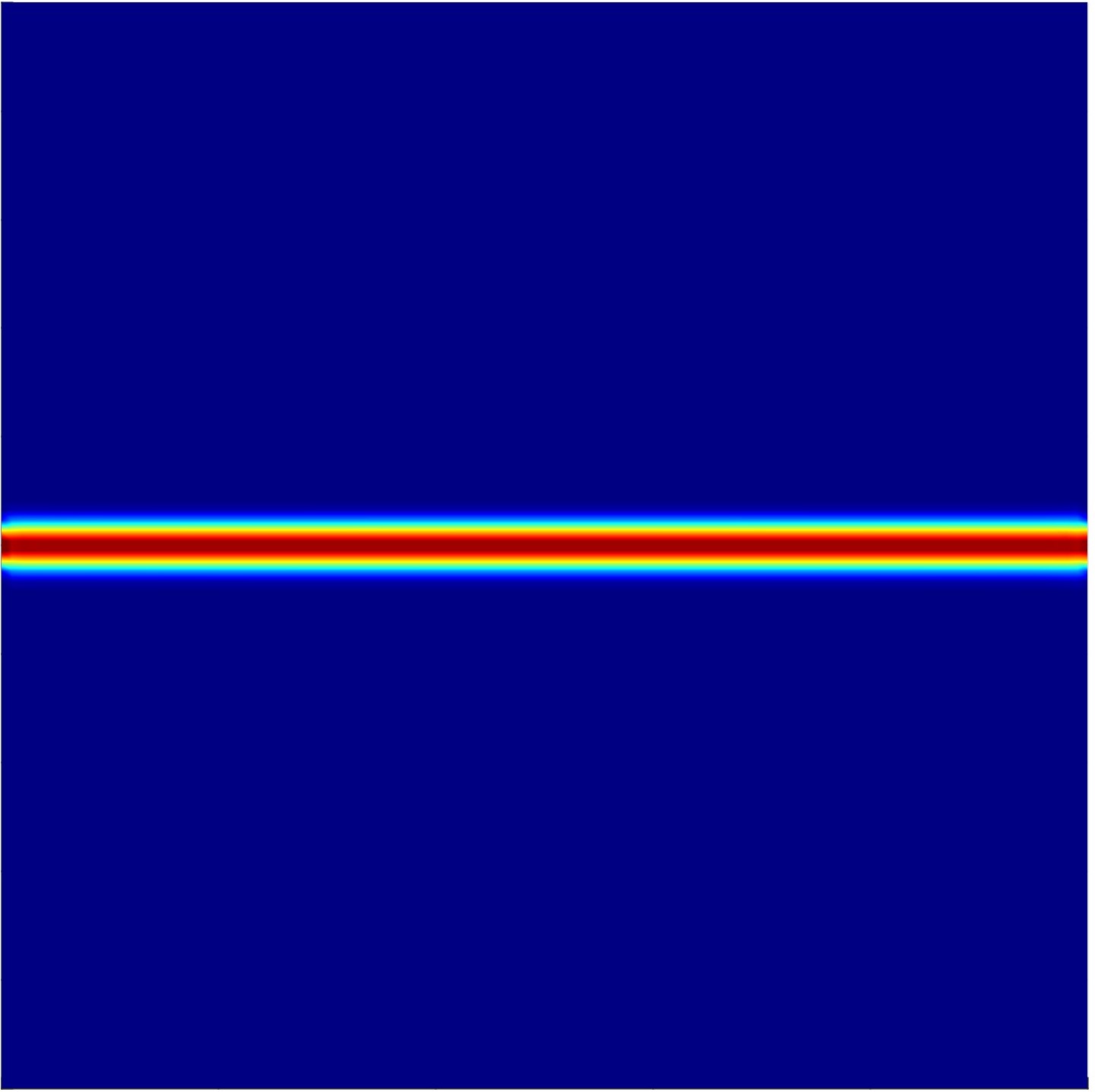}
			\caption{$\v$= 1.}
		\end{subfigure}
		\begin{subfigure}[b]{0.18\textwidth}
			\includegraphics[width=2.5cm]{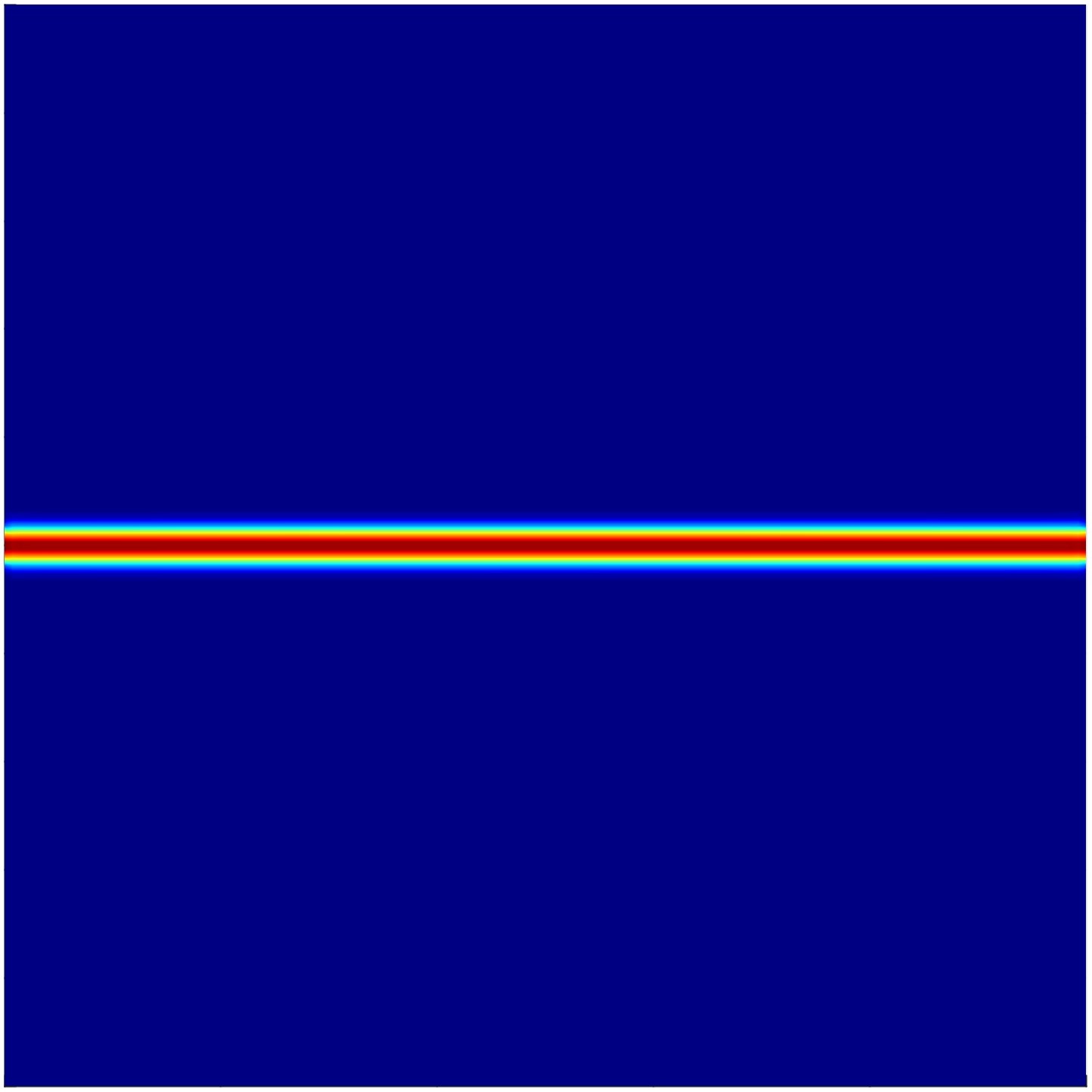}
			\caption{$\v$= 1/4.}
		\end{subfigure}
		\begin{subfigure}[b]{0.18\textwidth}
			\includegraphics[width=2.5cm]{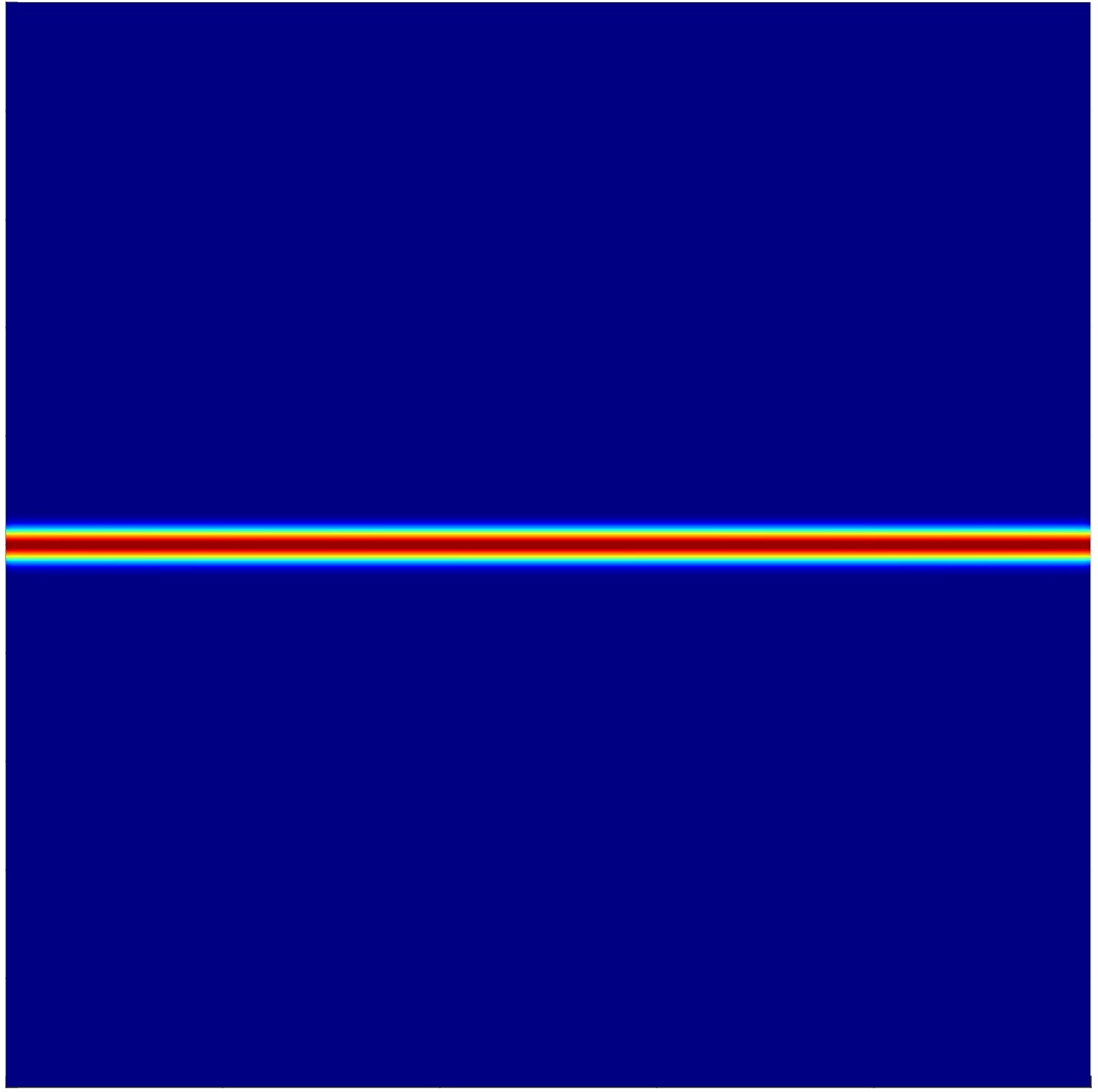}
			\caption{$\v$= 1/16.}
		\end{subfigure} 
	\end{center}
	\vspace{-15pt}
	\caption{SEN tensile test: influence of the variation of coefficient $\v$ on the final crack pattern for mesh size computed {as a function of} the $R_*$ parameter.}
	\label{Crack coeff vary SEN T mv}
\end{figure}

\begin{figure}[h!!]
	\begin{center}
		\begin{subfigure}[b]{0.18\textwidth}
			\includegraphics[width=2.5cm]{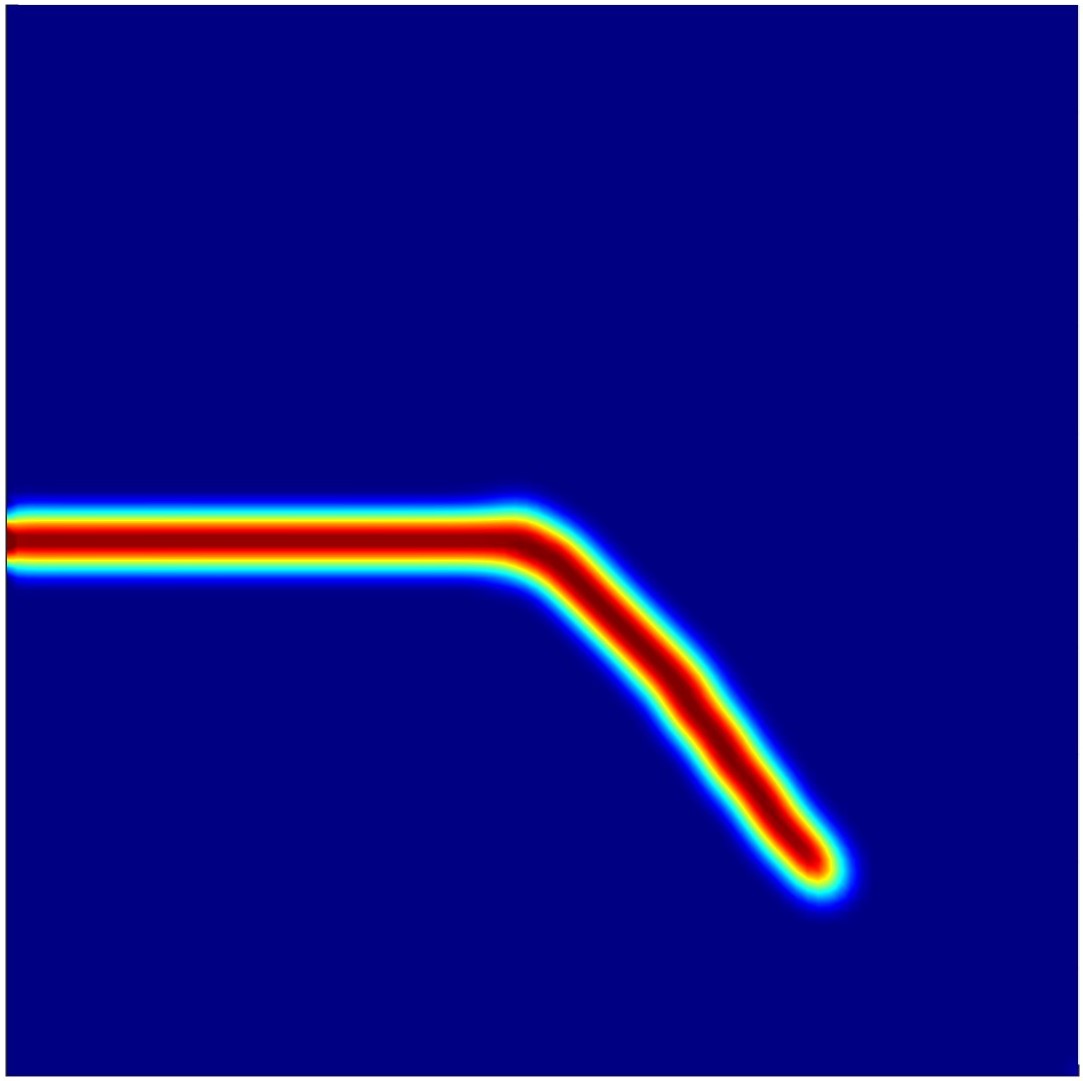}
			\caption{$\v$ = 4.}
		\end{subfigure}
		\begin{subfigure}[b]{0.18\textwidth}
			\includegraphics[width=2.5cm]{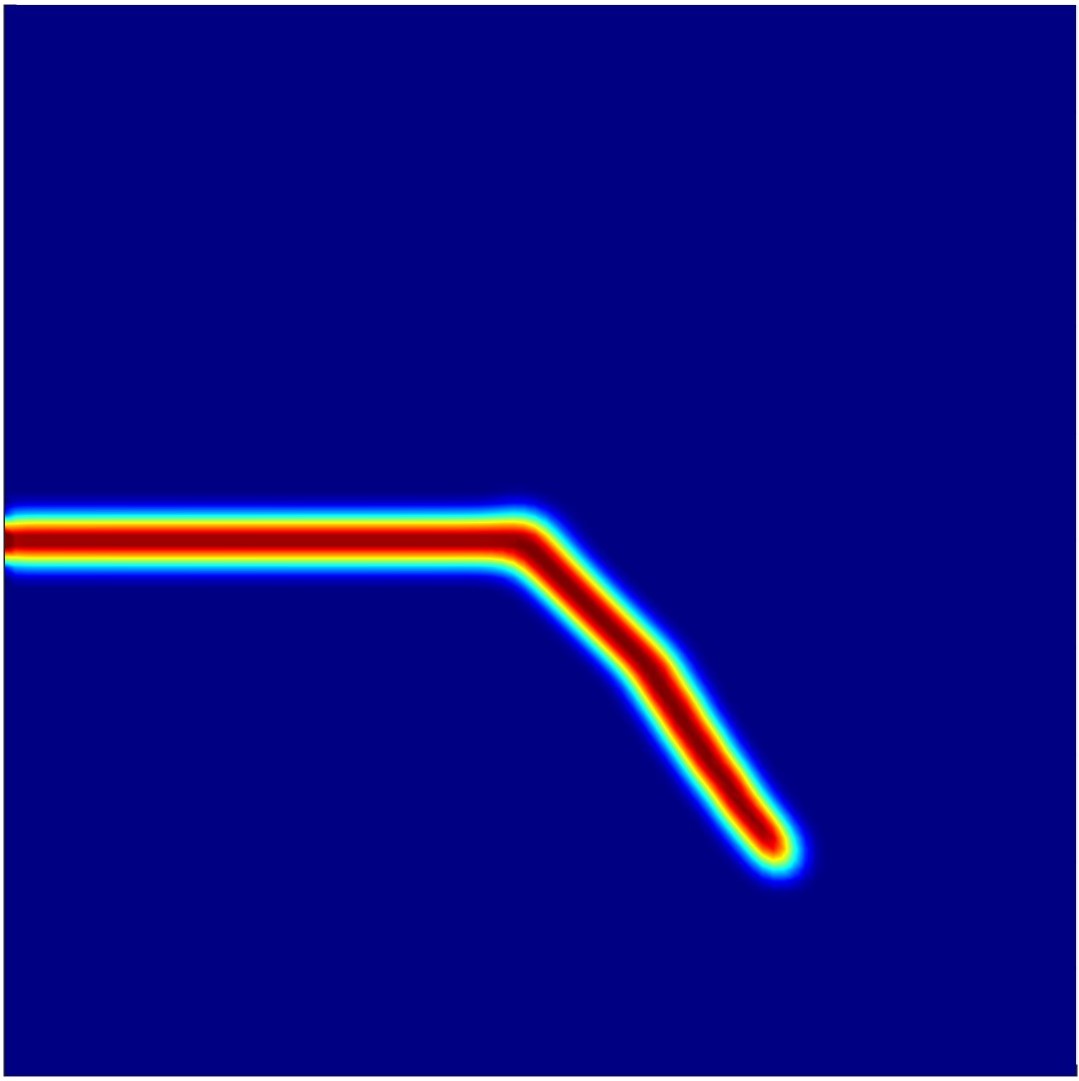}
			\caption{$\v$= 2.}
		\end{subfigure}
		\begin{subfigure}[b]{0.18\textwidth}
			\includegraphics[width=2.5cm]{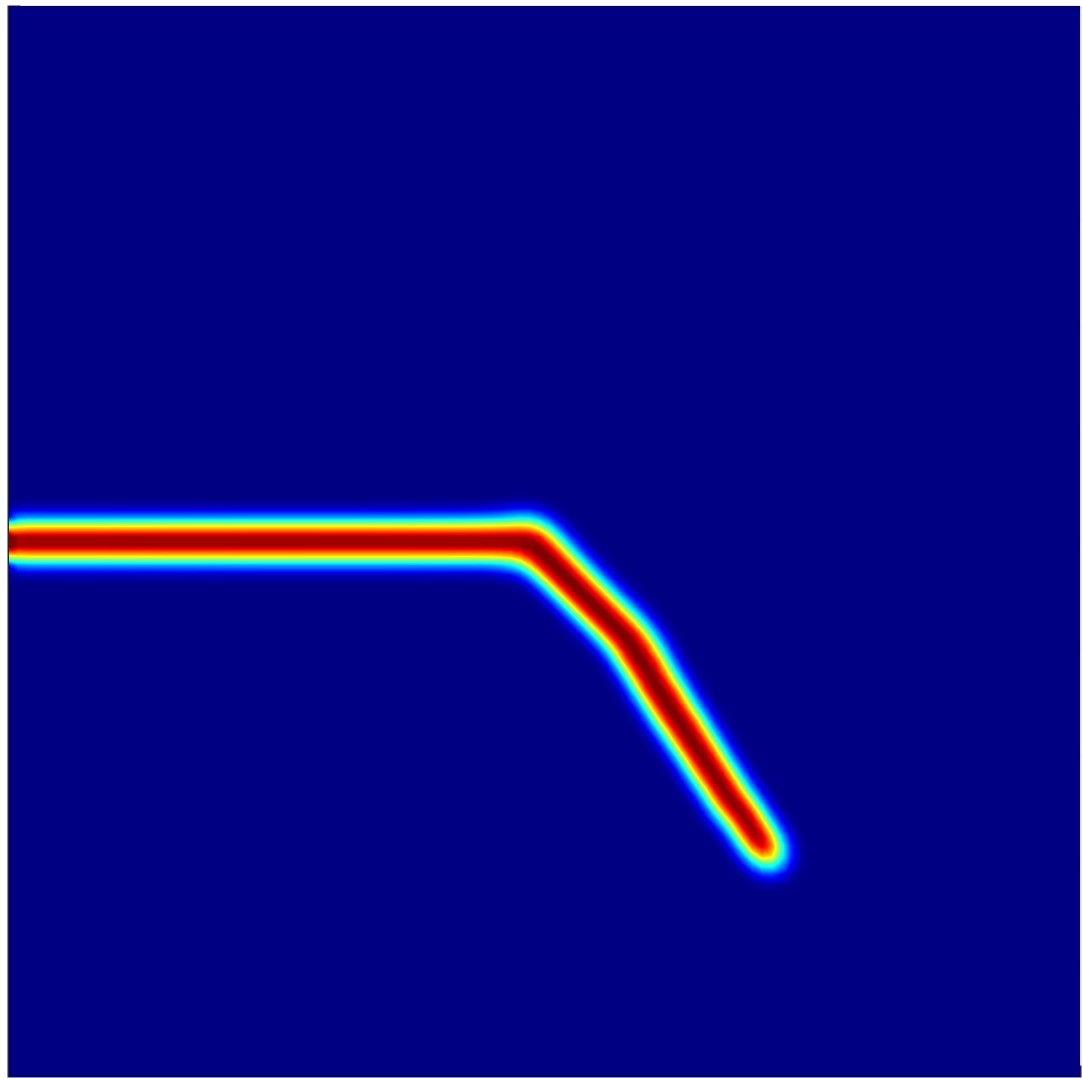}
			\caption{$\v$ = 1.}
		\end{subfigure}
		\begin{subfigure}[b]{0.18\textwidth}
			\includegraphics[width=2.5cm]{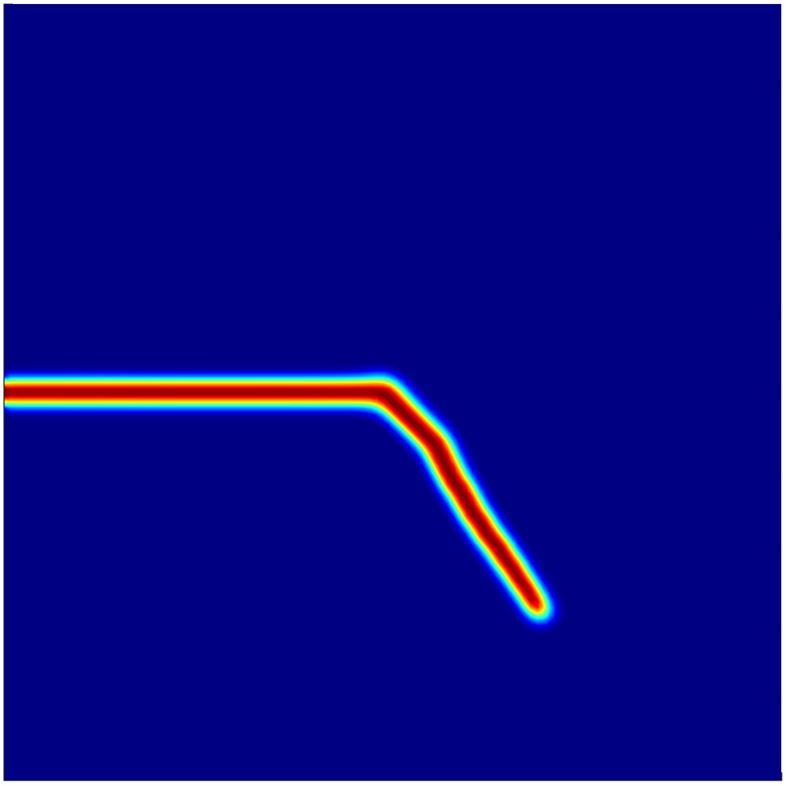}
			\caption{$\v$= 1/2.}
		\end{subfigure}
		\begin{subfigure}[b]{0.18\textwidth}
			\includegraphics[width=2.5cm]{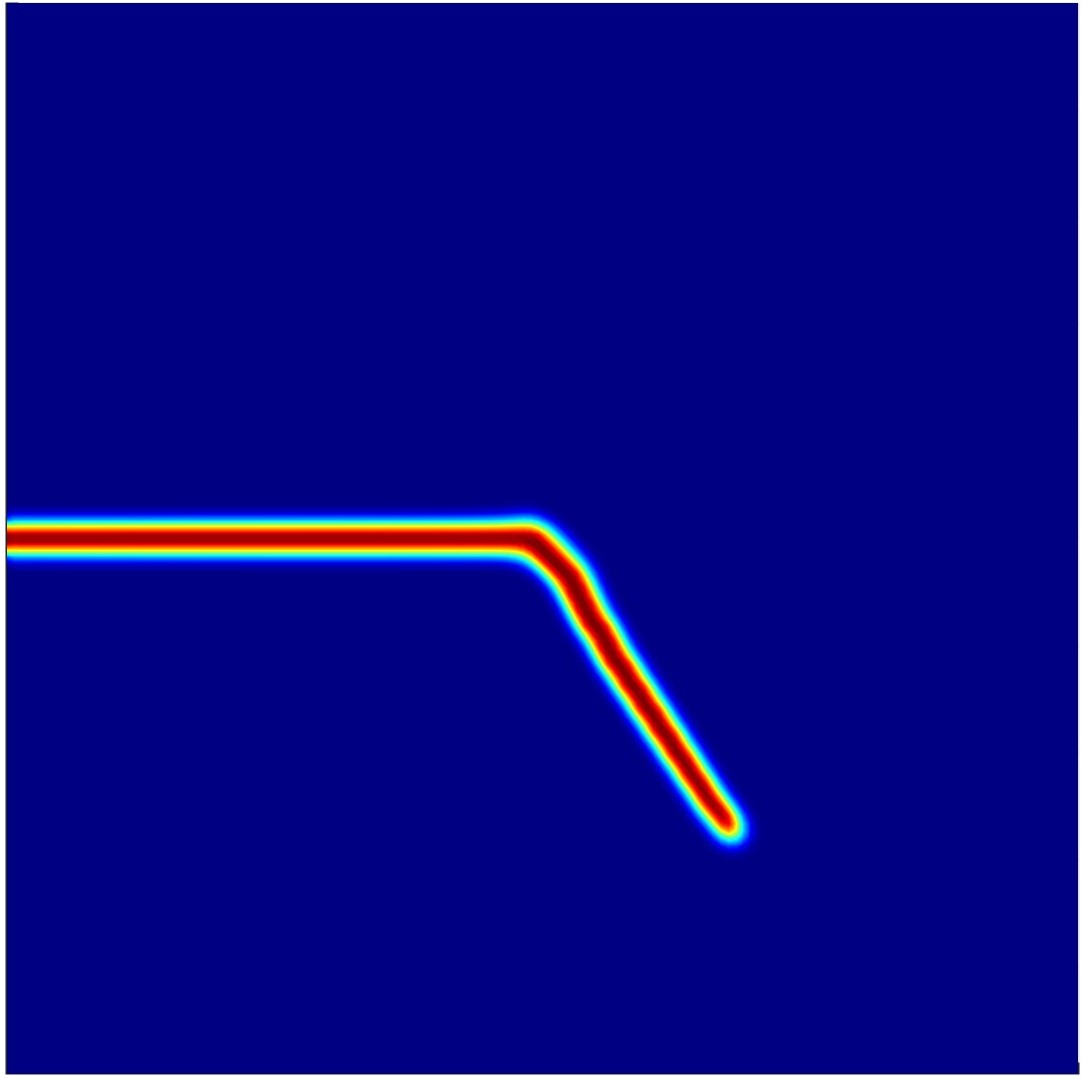}
			\caption{$\v$ = 1/4.}
		\end{subfigure} 
	\end{center}
	\vspace{-15pt}
	\caption{SEN shear test: influence of the variation of coefficient $\v$ on the final crack pattern for mesh size computed {as a function of} the $R_*$ parameter.}
	\label{Crack coeff vary SEN S mv}
\end{figure}
%\NOTE{AP: vambiato caption Figura 14 "SEN tensile test: influence of the variation of coefficient $\v$ on the final crack pattern for fixed mesh size" in attuale. Controllare.}
\begin{figure}[h!!]
	\begin{center}
	\begin{subfigure}[b]{0.5\textwidth}
		\includegraphics[width=9.5cm]{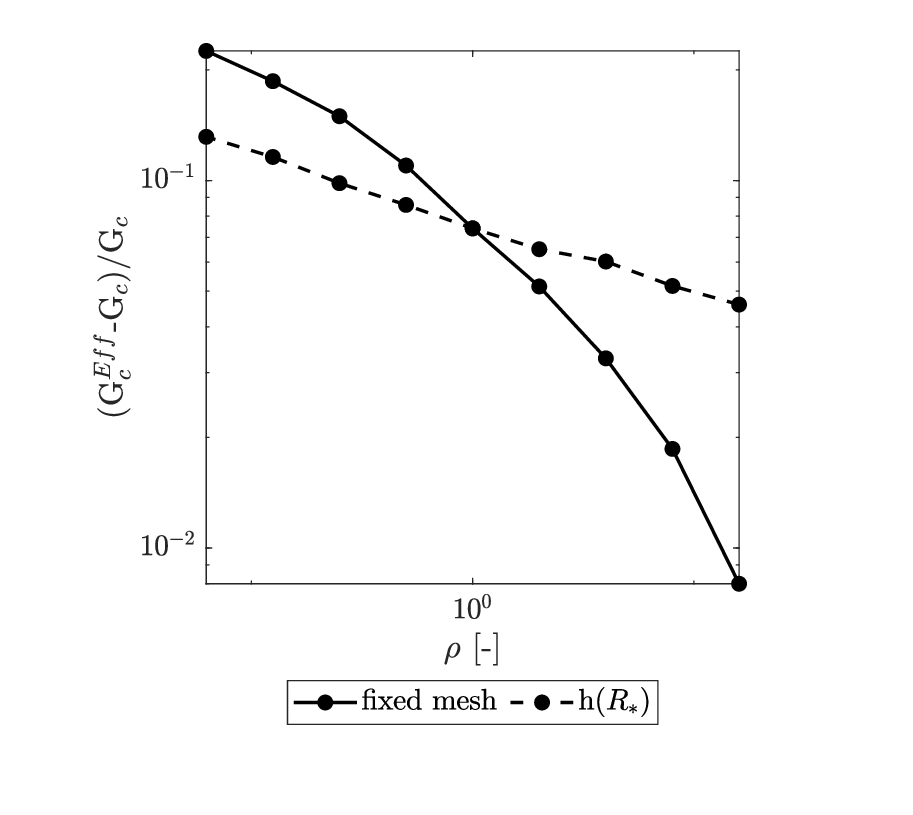}
		\caption{SEN tensile test.}
	\end{subfigure}
	\begin{subfigure}[b]{0.49\textwidth}
	\includegraphics[width=9.7cm]{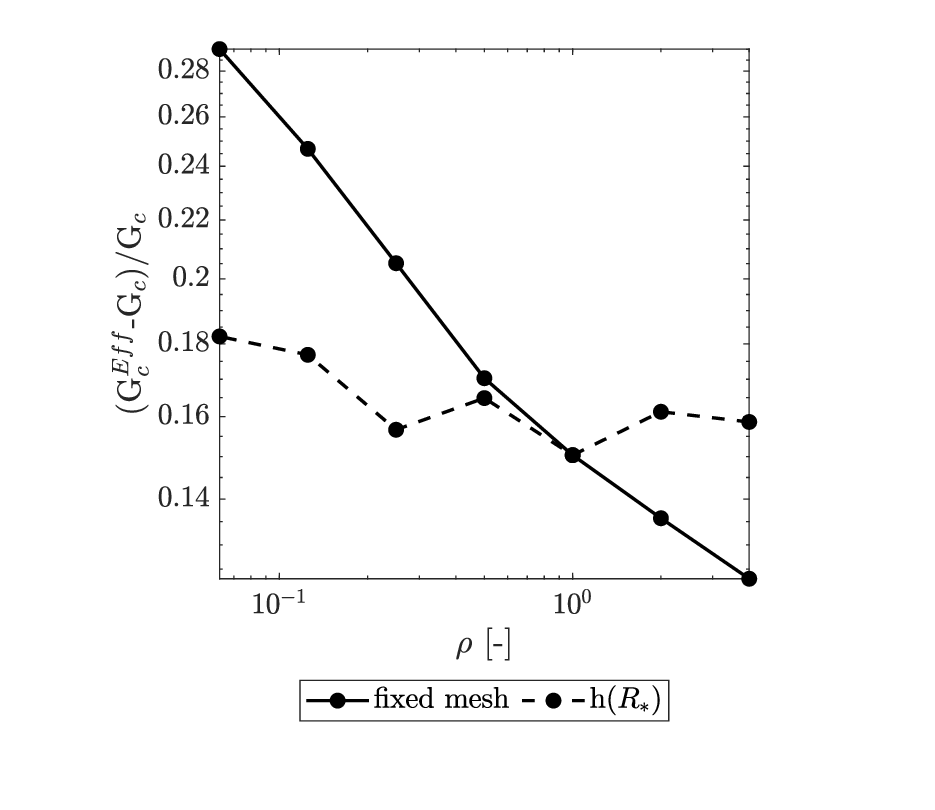}
	\caption{SEN shear test.}
	\end{subfigure}
	\end{center}
	\vspace{-15pt}
	\caption{{Influence of the mesh size on the relative error on the toughness for the SEN tensile and shear tests: comparison for fixed a mesh, {corresponding to the reference case $\v=1$ (i.e., $R_* = 3.83$),} and for a mesh varying with $R_*$.}}
	\label{SENT mesh  influence}
\end{figure}
{In Figure \ref{SENT mesh  influence}, we compare, for varying $\v$, the relative percentage error on the toughness for a fixed mesh, corresponding to the reference case $\v=1$ (i.e., $R_* = 3.83$), and for a mesh varying with $R_*$ (i.e., $h(R_*)$), namely a mesh size computed using Equation \eqref{mesh definition} with $n = 4$.
It can be observed that, for both the SEN tensile and shear tests, the error range for the $h(R_*)$ case is smaller than the one obtained for a fixed mesh, as this latter choice overestimates the required mesh-size value, being approximately twice larger than the one needed to accurately resolve the internal length, leading to a higher error in the toughness. Conversely, for $\v$ grater than 1 the mesh size is underestimated, thereby entailing a lower accuracy.
Then, we focus on the influence of the mesh on the elastic limit in the pure traction case and we highlight in Table \ref{Tab: c and R vary} that the numerical values of the elastic limit computed either for a fixed mesh corresponding to the $\v$ reference case, $\sigma_{c}$, or for a mesh varying with $R_*$, $\sigma_{c} (R_*)$, are comparable, so as, consequently, their corresponding relative errors on the elastic limit. This entails a possible significant reduction in terms of the total number of required control points to reach satisfactory results for this test.}
%\NOTE{AP: Ho fatto swap nel testo attuale dei termini overestimate e underestimate per i valori di mesh size, perche' h fissa e' circa $\epsilon$ per $\v=1$, per $\v=1/16$ sarebbe $h(R_*)=\epsilon/2$ quindi sovrastimo il valore di mesh size, sottostimo il numero di elementi. Controllare  Figure 16 log-log.}

\begin{table}[h!] \centering
	\caption{Pure traction test: comparison between relative percent errors (i) for a fixed mesh size R.err (i.e., corresponding to the reference case $\rho=1$) and (ii) for meshes varying with $R_*$ (R.err($R_*$)).
	} 
	%	\vspace{-5pt}
	\renewcommand\arraystretch{1.25}
	\renewcommand\tabcolsep{8pt}
	\begin{tabular}{c c c c : c c : c c}
		\hline
		\centering
		{{$\v$}} & {{{$R_*$}}} & {c$_\v$} & {$\sigma_{\text{th}}$} & $\sigma_{c}$&  $\sigma_{c}(R_*)$ & {R.}err & {R.}err$(R_*)$\\
		\hline
		{[-]} & {[-]} & [-] & {[kN/mm$^2$]} & {[kN/mm$^2$]} & {[kN/mm$^2$]} & {[\%]} & {[\%]} \\
		\hline
		16 & 7.1041 & 7.7811 & 1.0140 & 1.0199 & 1.0198 & 0.5871 & 0.5776\\ 
		8 & 6.0364 & 6.6769 & 1.0946 & 1.0973 & 1.0973 & 0.2440 & 0.2431\\ 
		4 & 5.1514 & 5.7717 & 1.1773 & 1.1784 & 1.1784 & 0.0929 & 0.0936\\ 
		2 & 4.4230 & 5.0369 & 1.2603 & 1.2605 & 1.2605  & 0.0197 & 0.0198\\ 
		1 & 3.8300 & 4.4485 & 1.3410 & 1.3410 & 1.3410  & 0.0026 & 0.0026\\ 
		1/2& 3.3554 & 3.9852 & 1.4168 & 1.4167 & 1.4167& 0.0115 & 0.0116 \\ 
		1/4 & 2.9847 & 3.6281 & 1.4849 & 1.4848 & 1.4848& 0.0101 & 0.0102\\ 
		1/8 & 2.7045 & 3.3593 & 1.5432 & 1.5430 & 1.5430& 0.0123 & 0.0123\\ 
		1/16 & 2.4998 & 3.1615 & 1.5907 & 1.5900 & 1.5900 & 0.0468 & 0.0468\\ 
		\hline
	\end{tabular}
\label{Tab: c and R vary}
\end{table}
%\NOTE{AP: Table 9: sistemare numero di cifre significative. $R$ e' $R_*$??}

% !Tex root=AT1IV.tex

\section{Conclusions}\label{sec:conclusions}

In this work, an $AT_1$ fourth-order model for phase-field brittle fracture is proposed. First, we prove a $\Gamma$-convergence result with a fine study of the optimal profile, providing (a) the explicit value of the corrector factor for the surface energy and (b) the explicit value of the width of the transition region. 

Numerically, we employ an IgA framework that allows to discretize the high-order differential operator thanks to the smoothness of the isogeometric shape functions, \aprev{highlighting once again the straightforward capability of IgA to model high-order PDEs}. We study the performance of our fourth-order functional in terms of accuracy, as far as the elastic limit and the effective toughness.

%\NOTE{Per quanto riguarda pure tensione mi verrebbe da dire che sono equivalenti, gli errori sono entrambi molto piccoli.}
%\NOTE{AP: Sono d'accordo e in quel caso possiamo proprio parlare di accuratezza.}

%confirms that higher-order formulations anticipate fracture nucleation \bb and are closer to the expected theoretical value

%Then, we propose a workflow for the the fourth-order AT1 profile derivation. This method has general validity and allows to obtain the profile shape and regularization constant c$_\ww$ for any combination of coefficients of the energy dissipation integral. Considering a set of coefficients $\{a, \,b, \,c\} = \{1, \, 1, \, 1\}$, 

 About the accuracy in the evaluation of the toughness, we observed that the proposed $AT_1$ fourth-order model is more accurate than the classical $AT$ ones present in the literature. These results confirm our findings in \cite{greco2024higher}, pointing out that %also the proposed AT1 
(in terms of accuracy) our fourth-order $AT_1$ formulation attains the same level of accuracy as other $AT$-formulations at a lower cost. %one when considering a coarser mesh. 

%As the dissipation functional coefficients considerably influence the shape of the profile, as well as the limit of the interval where the phase-field is positive, 

We also propose defining the mesh size as a function of the width of the phase-field profile, as characterized by $\Gamma$-convergence, leading to a significant reduction in the number of degrees of freedom. 

%\textit{R} parameter obtaining meshes $R$ times larger than the ones obtained considering the classical rule of thumb, i.e.,  four element resolving the internal length parameter.

%\NOTE{Qui starei attento e non imposterei tutto sul risparmio computazionale perche' il SENShear ha dei problemi per $\rho$ grande. In questo, lasciatemi aggiungere una nota: se si calcolasse anche l'energia elastica lungo l'evoluzione, temo che i risultati con $\rho$ grande si discosterebbero abbastanza dal modello sharp perch\`e la zona di transizione \`e piuttosto ampia.} 

Finally, a parametric analysis on the weight of the higher-order term \bl seems to indicate that, even though there is not a specific set of coefficients providing significantly better results, high values of $\rho$ give more regularity to the solution, tend to reduce errors, and allow for larger mesh sizes. 
%Consequently, this entails a significant decrease in terms of the number of degrees of freedom of the problem and the overall computation time. 
\bl{Among future works, we will extend current formulation to simulate fracture in other structural elements such as Kirchhoff-Love shells and study its performance in dynamics. Furthermore, we will possibly include in the present simulation framework adaptive THB-spline refinement and multi-patch geometries \cite{bracco2023adaptive}.} \bl 

\section*{Acknowledgements}
%\NOTE{AP: Da controllare, ho inserito quelli che avevo.}
%This work was partially supported by the Italian Ministry of University and Research (MUR) through the PRIN project XFAST-SIMS (No. 20173C478N). Such a support is gratefully acknowledged.\\\\
L. Greco gratefully acknowledges the ``Erasmus + Traineeship 2023/24'' and ``Bando di mobilità internazionale - 15a edizione'' programmes, that partially supported him within the collaboration between the University of Pavia and the Universit{\"a}t der Bundeswehr M{\"u}nchen.\\\\
%A. Patton was partially supported by the European Research Council (ERC) under the European Union's Horizon 2020 research and innovation programme (grant agreement No. 864482).\\\\
A. Reali is a member of the Gruppo Nazionale Calcolo Scientifico-Istituto Nazionale di Alta Matematica (GNCS-INDAM), and acknowledges the support of the Italian Ministry of University and Research (MUR) through the PRIN project COSMIC (No. 2022A79M75), funded by the European Union - Next Generation EU, as well as the contribution of the National Recovery and Resilience Plan, Mission 4 Component 2 - Investment 1.4 - NATIONAL CENTER FOR HPC, BIG DATA AND QUANTUM COMPUTING, spoke 6.\\\\
\bl E. Maggiorelli and M. Negri are members of the Gruppo Nazionale Analisi Matematica Probabilit\`a Applicazioni-Istituto Nazionale di Alta Matematica (GNAMPA-INDAM) and acknowledge the support of the Italian Ministry of University and Research (MUR) through the PRIN project ``Variational methods for stationary and evolution problems with singularities and interfaces'' (No. 2022J4FYNJ). \bl

\bl
%\NOTE{AP: Ho spostato l'appendice prima delle references.}

\appendix

\section{Optimal profile for the second-order {\boldmath $AT_1$} \aprev{model}} \label{optimalAT1secondo}

In this section, we provide a study of the optimal profile for the second-order $AT_1$ functional. Although the results are well known, here we provide a new proof along the lines of \S \ref{optimal}.

Let $\U = H^1 (a,b)$  be the space of the admissible displacements and let $ \V =  H^1 ( (a,b)  ; [0,1]) $ be the set of the admissible phase-field functions. For $\eps>0$ and $\eta_\eps = o (\eps)$, let $\F_\eps : L^1 (a,b) \times L^1(a,b) \to [0,+\infty]$ be the $AT_1$ (second-order) functional defined by 
\begin{align} \label{e.FepsAP}
	\F_\eps ( u , v) = \begin{cases} 
		{\displaystyle \int_{(a,b)} \left(\kappa \psi_\eps (v)  | u' |_+^2 + \kappa |u'|_-^2  \right) dx +  \frac{G_c}{c} \int_{(a,b)}  \left(\eps^{-1} v + \eps | v' |^2  \right) dx }  & (u,v) \in \U \times \V  , \\ 	+\infty	& \text{otherwise,}
	\end{cases}
\end{align}
where $\psi_\eps (v) = (v-1)^2 + \eta_\eps $ while $c >0 $ will be given later. Let $\F : L^1 (a,b) \times L^1(a,b) \to [0,+\infty]$ be the functional
\begin{align} \label{e.FlimAP}
	\F ( u , v)  = \begin{cases} 
		{\displaystyle \int_{(a,b) \setminus J_u} \kappa | u' |^2 \, dx + G_c \hspace{0.3pt} \# J_u } & \text{if $u \in SBV^2 ( a,b)$ with $\jump{u} > 0$ and $v =0$ a.e.~in $(a,b)$}, \\
		+\infty	& \text{otherwise.}
	\end{cases}
\end{align}
%
%In other terms, the energy $\F$ is finite if $u' \in L^2(a,b)$, the number of discontinuities is finite, and the non-interpenetration condition $\llbracket u \rrbracket > 0$ is satisfied in the jump points.  Before stating our main convergence result, it is necessary to define the optimal profile problem, which ultimately will provide the constant $c_\v$ for the $\Gamma$-convergence result.
%
In \cite{Braides98}, the value $c$ is provided explicitely and is given by 
$$
c = 4 \int_0^1 v^{1/2}= 8 /3 . 
$$
Here, we provide a different proof employing the explicit computation of the optimal profile following the lines of \S \ref{optimal}.  To this end, let $\mathcal{K} : \mathcal{W} \to \R$ be defined by
$$
\mathcal{K} ( w) = \int_{(0,+\infty)} \left(w + | w' |^2 \right) ds 
$$
where $\mathcal{W} = \{ w \in L^1 (\R_+) : w' \in L^2(\R_+), \, w(0)=1 \}$ is the domain of $\mathcal{K}$.

\begin{proposition} \label{p.wopbis} There exists a unique 
	\begin{equation*} %\label{e.optprofbis} 
		w_* \in \argmin \{ \K (w) : w \in \mathcal{W} \text{ such that } 0 \le w \le 1 \} . 
	\end{equation*} 
	Moreover, $\mathrm{supp} (w_*) = [0, R_*]$, where $R_*>0$ is characterized in Proposition \ref{p.eq}.
\end{proposition}
\proof Take $\{w_n\}_{n\in\mathbb N}$ such that $\mathcal{K}(w_n)\to \inf\{\mathcal{K}(w) : w \in \mathcal{W} \text{ such that }  0 \le w \le 1\}$. %This $\{w_n\}_{n\in\mathbf N}$ is customarily called \textit{minimizing sequence}. 
Observe that the functions $w_n$ are bounded in $ H^1(\R_+)$. Indeed, being $0\leq w_n\leq 1$, we have $\|w_n\|^2_{H^1(\R_+)}=\int_{\R_+}\left(|w_n|^2+|w_n'|^2\right)\aprev{dx}\leq \int_{\R_+}\left(|w_n|+|w_n'|^2\right)\aprev{dx}<+\infty$.  Therefore, it exists a (non-relabelled) subsequence of $\{w_n\}_{n\in\mathbb N}$ that weakly converges in $H^1(\R_+)$ to a certain $w_*\in H^1(\R_+)$. 
The set $\{w \in \mathcal{W}\,:\, \,0\leq w\leq1\}$ is weakly closed (being convex and strongly closed) in $H^1(\R_+)$ and hence $w_*$ belongs to this set. 
Now, since the functional $\mathcal{K}$ is weakly lower semicontinuous (being strictly convex), %(being convex and strongly continuous BREZIS), 
we obtain that $w_*$ is indeed the unique minimum in definition \eqref{e.optprof}. \qed 

%As a matter of fact, $w_*|_{[0,R_*]}\in \argmin\{\K_{R_*}(w) : w \in \mathcal{W}_{R_*} \text{ such that } w(0) = 1,  0 \le w \le 1 \} . $

%\begin{remark} \normalfont 
\begin{theorem} \label{t.Gamma-bis} Let $c = 2 \mathcal{K} (w_*) $. Then $\F_\eps$ $\Gamma$-converge to $\F$ as $\eps \to 0^+$ with respect to the (strong) topology of $L^1 ( a,b ) \times L^1( a,b)$.
\end{theorem}

In order the characterize $w_*$, $R_*$, and $c$, it is convenient to introduce for $R\in(0,+\infty)$ the localized energies $\mathcal{K}_R : \mathcal{W}_R \to \R$ given by 
$$
\mathcal{K}_R ( w) = \int_{(0,R)} \left(w + | w' |^2 \right) dx , 
$$
where $\mathcal{W}_R =\{w\in H^1(0,R)  : w(0)=1, w(R)=0 \}$. We will also employ the local unconstrained optimal profiles
$$
w_R \in \argmin \{ \mathcal{K}_R (w) : w \in \mathcal{W}_R \} . 
$$
Clearly, the above minimizer is unique and it is characterized by the ODE
\begin{equation*}\begin{cases} 1-2w''=0\qquad\text{on $(0,R)$} \\w(0)=1\\w(R)=0. \end{cases}
\end{equation*}
Hence $w_R(x)=\frac{1}{4}x^2-\frac{1}{R}(1+\frac{1}{4}R^2)x+1$. In analogy with \S \ref{optimal}, the relationship between the local unconstrained profiles $w_R$ and the optimal profile $w_*$ is given by the following Proposition.

\begin{proposition} \label{p.eq-app} It holds $R_* = \max\{R \in(0,+\infty)\,:\,w_R\in[0,1]\}=2$,
	moreover, $w_*|_{[0,R_*]}(x)= w_{R_*}(x)$, i.e., $w_* (x) = \tfrac14 x^2 - x +1$. 
	%where $w_{R_*}$ is the unique solution of \eqref{EL_R} with $R=R_*$. 
	%
	%The support of the function $w_*$ is $[0,2]$ and $w_*|_{[0,2]}$ is the unique solution of the problem of finding $w$ such that
	%\begin{equation}\begin{cases} 1-2w''=0\qquad\text{on $(0,2)$} \\w(0)=1\\w(2)=0\end{cases}
	%\end{equation}
\end{proposition}

\proof %The Euler-Lagrange equation associated to $ w_* \in \argmin \{ \K (w) : w \in \mathcal{W} \text{ such that } w(0) = 1 \} $ reads:\begin{equation}\begin{cases} 1-2w_*'=0\qquad\text{on $(0,R_*)$} \\w_*(0)=1\end{cases}\end{equation}
% Let's consider the general problem of of finding $w_R$ such that 
%\begin{equation}\label{EL_R}\begin{cases} 1-2w_R''=0\qquad\text{on $(0,R)$} \\w_R(0)=1\\w_R(R)=0\end{cases}
%\end{equation}
% Its solution is $w_R(x)=\frac{1}{4}x^2-\frac{1}{R}(1+\frac{1}{4}R^2)x+1$. 
%
Since $w_R$ is quadratic and convex, with $w_R(0)=1$ and $w_R(R)=0$, requiring that $w_R \in [0,1]$ is equivalent to asking that $w_R'(R)\leq0$, that gives that $R\leq 2$, hence $\max\{R \in(0,+\infty)\,:\,w_R\in[0,1]\}=2$. 

%, i.e. $\max\{R\in(0,+\infty)\,:\,w_R\in[0,1]\}=2$. xcv

%Notice that $\inf\{x\in\R_+\,:\,w_*(x)<1\}=0$.
Notice that $w_*(x)<1$ for all $x>0$. Indeed, if it existed  $R'>0$ such that $w_*|_{[0,R']}=1$, then $\hat w(x)=w_*(x+R')$ would be such that $\mathcal{K} (\hat w) < \mathcal{K} (w_*)$, that is absurd.
We define $R_*:=\inf\{x\,:\,w_*(x)=0\}$ %\in(0,+\infty]$ 
and obviously $w_*|_{[R_*,+\infty)}=0$. 
Now, %$R_*\leq 2$, since 
$w_*|_{[0,R_*]}$ is the minimizer of $\mathcal{K} _{R_*}$ with boundary conditions $w(0)=1$ and $w(R_*)=0$  and with the constraint $w_* \in [0,1]$. However, $0 < w_*(x) < 1$ for $ x \in (0,R_*)$ (by definition of $R_*$), hence $d \mathcal{K} ( w_*) [ \xi] =0$ for every $\xi \in C^\infty_c (0,R_*)$, which means that $w_*$ solves 
the Euler-Lagrange equation
\begin{equation}\label{EL_*}\begin{cases} 1-2w_*''=0\qquad\text{on $(0,R_*)$} \\w_*(0)=1\\w_{R_*}(R_*)=0 . \end{cases}
\end{equation}
%\eqref{EL_R} on $(0,R_*)$ 
%for hypothesis \eqref{e.optprof}. 
In other terms, $w_*|_{[0,R_*]}(x)= w_{R_*}(x)$. As a straightforward consequence, since $w_*\in[0,1]$, $R_*\leq \max\{R\in(0,+\infty):w_R\in[0,1]\}=2$. 

On the other hand, $w_2\in \argmin\{\mathcal{K}_2(w):w\in \mathcal W_2  \}$. Since $w_2\in [0,1]$ we also have  $w_2\in \argmin\{\mathcal{K}_2(w):w\in \mathcal W_2 \,, w \in[0,1]\}$.
Observe that $w_*|_{[0,2]}\in \argmin\{\mathcal{K}_{2}(w) : w \in \mathcal{W}_{2}\,, w\in [0,1] \} $ (indeed $R_*<2$ and $w_*|_{[R_*,+\infty)}=0$) and by the uniqueness of the minimizer $w_*|_{[0,2]}=w_2=\frac{1}{4}(x-2)^2 $. Hence $R_*= 2$. %This means that $w_*|_{[0,2]}$ solves the Euler-Lagrange equation \eqref{EL_*} for $R_*=2$.
\qed

\begin{proposition} The explicit value of the optimal constant is $c=\frac{8}{3}$.
\end{proposition}
\proof The constant $c$ is given by $2\mathcal{K}(w_*)$. By direct calculations we get that $c=R_*-\frac{1}{24}R_*^3+\frac{2}{R_*}$ and since $R_*=2$, then $c=\frac{8}{3}$.
\qed
\bibliographystyle{plain}
\bibliography{bibliografia} 

\end{document}